\documentclass[preprint,1p,11pt,authoryear]{elsarticle}

\usepackage{amssymb}
\usepackage{amsmath}
\usepackage{amsthm}

\usepackage[linesnumbered,ruled]{algorithm2e}
\usepackage{algorithmic}
\usepackage{etoolbox}
\usepackage{booktabs}
\usepackage{multicol}
\usepackage{setspace}
\usepackage{microtype}
\usepackage{wrapfig}
\usepackage{color}

\usepackage{natbib}
 \bibpunct[, ]{(}{)}{,}{a}{}{,}%
 \def\newblock{\ }%

\newcommand{\blue}[1]{{\color{black} #1}}

\def\ssum{\textstyle\sum}
\def\sprod{\textstyle\prod}
\def\sbinom#1#2{\textstyle\left(\! {#1\atop #2}\!\right)}

\def\P{\mathbb{P}}

\newtheorem{thm}{Theorem}

\newtheorem{rem}{Remark}

\newtheorem{lem}{Lemma}

\newtheorem{prop}{Proposition}

\newenvironment{proof}{\paragraph{Proof}}{\hfill $\blacksquare$}

\let\oldnl\nl
\newcommand{\nonl}{\renewcommand{\nl}{\let\nl\oldnl}}

\SetKwInput{KwInput}{Input}                
\SetKwInput{KwOutput}{Output}              

\allowdisplaybreaks[4]

\newcommand{\zerodisplayskips}{%
  \setlength{\abovedisplayskip}{0pt}%
  \setlength{\belowdisplayskip}{0pt}%
  \setlength{\abovedisplayshortskip}{0pt}%
  \setlength{\belowdisplayshortskip}{0pt}}
\appto{\normalsize}{\zerodisplayskips}
\appto{\small}{\zerodisplayskips}
\appto{\footnotesize}{\zerodisplayskips}

\journal{European Journal of Operational Research}

\begin{document}

\begin{frontmatter}



\title{Chance-Constrained Set Multicover Problem} 


\author[1]{Shunyu Yao} 
\author[1]{Neng Fan\corref{cor1}}
\cortext[cor1]{Corresponding author.}
\ead{nfan@arizona.edu}
\author[1]{Pavlo Krokhmal}

\affiliation[1]{organization={Department of Systems and Industrial Engineering},
            addressline={University of Arizona}, 
            city={Tucson},
            state={AZ},
            postcode={85721}, 
            country={USA}}

\begin{abstract}
We consider a variant of the set covering problem with uncertain parameters, which we refer to as the chance-constrained set multicover problem (CC-SMCP). In this problem, we assume that there is uncertainty regarding whether a selected set can cover an item, and the objective is to determine a minimum-cost combination of sets that covers each item $i$ at least $k_i$ times with a prescribed probability. To tackle CC-SMCP, we employ techniques of enumerative combinatorics, discrete probability distributions, and combinatorial optimization to derive exact equivalent deterministic reformulations that feature a hierarchy of bounds, and develop the corresponding outer-approximation (OA) algorithm. Additionally, we consider reducing the number of chance constraints via vector dominance relations and reformulate two special cases of CC-SMCP using the ``log-transformation" method and binomial distribution properties. Theoretical results on sampling-based methods, i.e., the sample average approximation (SAA) method and the importance sampling (IS) method, are also studied to approximate the true optimal value of CC-SMCP under a finite discrete probability space. Our numerical experiments demonstrate the effectiveness of the proposed OA method, particularly in scenarios with sparse probability matrices, outperforming sampling-based approaches in most cases and validating the practical applicability of our solution approaches.
\end{abstract}




\end{frontmatter}




\section{Introduction}\label{sec: introduction}
We consider the following \textit{chance-constrained set multicover problem} (CC-SMCP):
\begin{subequations}
\begin{align}
  \textbf{[CC-SMCP]} \qquad \min _x \quad & \ssum_{j \in [n]} c_j x_j \\
    \text { s.t. } \quad & \mathbb{P}\left[\ssum_{j \in [n]} \tilde{a}_{i j} x_j \geqslant k_i\right] \geqslant 1-\epsilon_i, \quad  \forall i \in [m],  \label{eq: set-k-covering-chance_constraints} \\
    & x \in B \subseteq \{0,1\}^n 
    \end{align}\label{eq: set-k-covering} 
\end{subequations}
where there are $m$ items, indexed by the set $[m] :=\{1, \ldots, m\}$, and $n$ subsets of these $m$ items, indexed by the set $[n]:=\{1, \ldots, n\}$. Here $c \in \mathbb{R}^n$ is the vector of weights, or costs of the subsets in $[n]$, $k_i \in \mathbb{Z}_+$, $\tilde{a}_{i j}\in \{0,1\}$ is a Bernoulli random variable indicating whether item $i\in [m]$ belongs to set $j\in [n]$ or not, and $\epsilon_i \in(0,1)$ is a prespecified allowed failure probability for item $i$.  The set $B$ represents some other deterministic side constraints on the decision variables $x$. For instance, in the facility location problem, $B$ can represent constraints requiring that demands at several points must be serviced by the established facilities, or it can be the budget restriction on the number of facilities to be located such as $\{x \in \{0,1\}^n: \sum_{j \in [n]}x_j \le U \}$ where at most $U$ facilities can be located. 

This problem aims to determine a subset $S \subseteq [n]$ that covers each element $i \in [m]$ at least $k_i$ times with high probability while minimizing the total cost/weight of covering. Note that when $k_i=1$ for all $i \in [m]$, the problem becomes {Probabilistic Set Cover (PSC) problem. The problems PSC and CC-SMCP are NP-hard. Since the standard Set Cover Problem (SCP) and Set Multicover Problem (SMCP) are NP-complete, their probabilistic variations, where elements may be covered with a certain probability, introduce an additional layer of complexity that makes them computationally challenging. Additionally, since different literature sources may have different definitions of PSC, in what follows we will not distinguish between the terms \emph{probabilistic set (multi-)cover} and \emph{chance-constrained set (multi-)cover.}
} 

{Our problem constitutes a variant of the SCP,} a well-known combinatorial optimization problem which serves as a model for a variety of real-world applications \citep{hall1992multicovering}, including marketing \citep{ozener2013allocating,sun2017chance}, security \citep{lessin2018bilevel}, telecommunications \citep{grotschel1992computational}, scheduling \citep{smith1988bus,zaghian2018chance}, production planning \citep{vasko1989set}, facility location \citep{gunawardane1982dynamic,huang2010facility}, and vehicle routing \citep{bramel1997effectiveness,bramel2002set}, etc. {Although many applications can be modeled as SCPs, for reliability and backup purposes, they are often better represented as SMCPs. Consider the critical task of determining optimal locations for emergency service facilities, such as hospitals or fire stations, to ensure comprehensive coverage of $m$ target regions across $n$ potential locations \citep{hall1992multicovering}. In emergency scenarios, such as natural disasters, multiple fire stations within a short distance of at-risk buildings provide a higher level of service compared to having only one nearby. This redundancy is crucial because it ensures that if one station is overwhelmed or inaccessible, others can still respond, thereby enhancing the reliability and resilience of emergency response. In this context, the objective is to strategically identify optimal locations for fire stations that guarantee all target regions receive the necessary level of service while minimizing overall costs. The connection parameter $a_{ij}$ plays a key role, indicating whether the service from candidate location $j$ can cover target $i$. This parameter often depends on distances, where $a_{ij} = 1$ if target $i$ is near location $j$, and $a_{ij} = 0$ otherwise.

In practical emergency scenarios, the connection between target $i$ and facility $j$ may be randomly disrupted due to factors such as road blockages or resource limitations. Therefore, modeling $a_{ij}$ as a Bernoulli random variable, denoted by $\tilde{a}_{ij}$, is appropriate \citep{shen2023chance}. The probability $\P[\tilde{a}_{ij} = 1] := p_{ij}$ represents the likelihood of successfully reaching target $i$ from location $j$ within a limited time frame during an emergency. This probabilistic approach acknowledges the inherent uncertainties and challenges in emergency response logistics. The importance of having multiple facilities (i.e., $k_i > 1$) for the same target becomes apparent in this context. It is not merely about redundancy but also about capacity planning and resource allocation. Multiple facilities ensure that even if one facility is incapacitated or resources are stretched, the remaining facilities can continue to provide critical services. This strategy is vital for maintaining service levels and managing resources effectively during high-demand situations. From the perspective of stochastic optimization, this problem can be formulated as a chance-constrained programming model that incorporates the described uncertainties. By doing so, we aim to identify optimal facility locations that achieve the required level of service for all targets with a high probability. This approach not only addresses the practical needs of emergency planning but also highlights the significance of studying SMCPs in real-world applications where reliability and resilience are paramount.}

In addition to the aforementioned application of CC-SMCP, \cite{beraldi2002probabilistic} considered a game to select arcs on a graph to cover adversaries' path with a high probability. \cite{ahmed2013probabilistic} presented an example to cover $m$ targets by placing sensors at $n$ potential sensor sites. These chance-constrained SCP models can be easily generalized to incorporate backup coverage, i.e., define a parameter $k_i$ to reflect the coverage level of target $i$. 

\subsection{{Literature review on chance-constrained optimization}}

As a powerful paradigm to model optimization problems with uncertainty, chance-constrained programming arises in a wide variety of applications such as finance \citep{lemus1999portfolio}, healthcare \citep{tanner2010iis}, power systems \citep{van2011chance}, transportation \citep{dinh2018exact}, network design \citep{wang2007beta} and wireless communications \citep{soltani2013chance}, etc. From the computational perspective, chance-constrained programs are very challenging to solve primarily for two reasons \citep{kuccukyavuz2022chance}. First, given a candidate solution $x$, it can be computationally demanding to check whether $x$ is feasible or not \citep{ahmed2013probabilistic}; and, second, the feasible region characterized by chance constraints is nonconvex and even disconnected in general \citep{nemirovski2007convex,ahmed2013probabilistic}. There are two main approximation techniques to mitigate the aforementioned challenges: convex conservative approximation \citep{ben2009robust,nemirovski2007convex} and sample average approximation (SAA) \citep{kleywegt2002sample,calafiore2005uncertain,luedtke2008sample}, both of which can efficiently identify feasible solutions with a guarantee of performance. {The strategy of the former is to formulate a convex optimization problem and produce feasible solutions with guarantees. It exploits the convexity of the resulting feasible region; however, these properties are absent in the discrete setting \citep{ahmed2013probabilistic}.} The latter is to solve an approximation problem based on an independent Monte Carlo sample of random data \citep{kleywegt2002sample}. This approach may not require knowledge of the distribution of the random parameters, however, it may need drawing a large number of samples from the true distribution to obtain feasible solutions with a guarantee. {For detailed specific relaxations and conservative approximations for CCO problems} (such as finite scenario approximation, CVaR approximation, Bernstein approximation, Bonferroni approximation, etc), we refer the reader to \citet{ahmed2018relaxations} and \citet{lejeune2025relaxations}. {In addition, \citet{prekopa1971logarithmic, prekopa1972class} showed that, under suitable assumptions, the probability function can exhibit favorable properties such as log-concavity.}

\subsection{{Literature review on chance-constrained combinatorial optimization}}
In many problems of interest, the decision vector $x$ is binary (as is the case in the present paper) and this special structure can be exploited to derive stronger formulations and specialized algorithms. We refer to such chance-constrained problems with purely binary variables as \textit{chance-constrained combinatorial optimization} problems \citep{kuccukyavuz2022chance}. \cite{luedtke2014branch} proposed a branch-and-cut decomposition algorithm for finding exact solutions of  chance-constrained problems having discrete distributions with finite support. For chance-constrained bin packing problems, \cite{song2014chance} proposed an efficient coefficient strengthening method and lifted probabilistic cover inequalities for chance-constrained bin packing problems. Later, \cite{wang2021chance} studied  bilinear formulation of chance-constrained multiple bin packing problem and used the lifting techniques to identify cover, clique, and projection inequalities to strengthen the bilinear formulation. Recently, \cite{zhang2020branch} proposed Dantzig–Wolfe formulations suited to a branch-and-price (B\&P) algorithm to solve two versions of the chance-constrained stochastic bin packing (CCSBP) problem.
For chance-constrained assignment problems, \cite{wang2022solution} considered a chance-constrained assignment problem to develop valid inequalities and derived lifted cover inequalities based on a bilinear reformulation of this problem. For chance-constrained dominating set problems, \cite{sun2019probabilistic} investigated the reliable connected power dominating set problem that met two requirements including the connectivity of the phasor measurement unit (PMU) subgraph and the reliability of its connectivity. For chance-constrained traveling salesman problems (TSPs), \cite{padberg1989branch} proposed a branch-and-cut approach to a TSP with side constraints, and \cite{campbell2008probabilistic} presented two recourse models and a chance-constrained model for probabilistic TSP with deadlines. For chance-constrained knapsack problems, \cite{yoda2016convexity} determined sufficient conditions for the convexity of the formulation under different discrete distributions, and several approximate but more tractable formulations that could provide near-optimal solutions were derived in \citet{klopfenstein2008robust,de2018boolean}, and \citet{han2016robust}. For more references on chance-constrained combinatorial problems and approaches, we refer the reader to Section 2 of survey \citep{kuccukyavuz2022chance}. 

\subsection{{Literature review on chance-constrained set covering problems}}
In the context of chance-constrained SCP models, prior studies have examined two types of uncertainty: right-hand-side (RHS) uncertainty and left-hand-side (LHS) uncertainty. For RHS uncertainty, \cite{beraldi2002probabilistic} developed a specialized branch-and-bound algorithm based on the enumeration of $p$-efficient points. Later, \cite{saxena2010mip} derived the polarity cuts to obtain a stronger formulation and improve the computational performance of this enumeration approach. For LHS uncertainty (as in this paper), \cite{fischetti2012cutting} developed cutting plane approaches to the individual chance-constrained problem where all components of the Bernoulli random vector are independent. Further, \cite{wu2019probabilistic} proposed an exact approach to the probabilistic partial SCP where there existed an efficient probability oracle to retrieve the probability of any events under the true distribution. For distributionally robust chance-constrained SCPs, \cite{ahmed2013probabilistic} studied individual chance constraints under moment-based ambiguity sets, and recently, \cite{shen2023chance} explored joint chance constraints with Wasserstein ambiguity.

As we can see, while there are several works in the literature related to chance-constrained SCPs, it appears that prior studies have paid less attention to chance-constrained set multicover problems, which represent a natural extension of SCPs. {In CC-SMCPs, the requirement that each target must be covered multiple times ($k_i > 1$) adds a layer of complexity not present in PSC problems, where each target needs to be covered only once. This multiplicity requirement significantly increases the problem's combinatorial nature, as we need to consider multiple overlapping sets to ensure sufficient coverage. This can lead to a larger solution space and more intricate interactions between decision variables, complicating both the problem's formulation and its solution process. Consequently, the traditional `log-transformation' technique \citep{haight2000integer,fischetti2012cutting,ahmed2013probabilistic}, commonly used in chance-constrained SCPs, is unsuitable for CC-SMCPs. {In this paper, we derive reformulations using the inclusion--exclusion principle, rather than relying on log-transformation. Beyond this, several commonly used approaches in stochastic optimization include risk-measure surrogates \citep{nemirovski2007convex,alexander2004comparison}, probabilistic envelope constraints (PEC) \citep{xu2012optimization} and sampling-based approaches \citep{luedtke2008sample,song2014chance,zhang2020branch}. The risk-measure surrogates (e.g., CVaR-based approximations) are widely used in distributionally robust CCO to approximate chance constraints; probabilistic envelope constraints generalize chance constraints by imposing probabilistic guarantees across all potential violation levels, thereby bounding the entire tail behavior via an envelope function. Sampling-based reformulations approximate the original problem using an independent Monte Carlo sample of the underlying random data. Our approach is different: rather than replacing the chance constraint with a risk-measure surrogate, an envelope surrogate or a sampling-based formulation, we exploit the discrete multicover structure and derive exact deterministic reformulations, enabling exact (rather than approximate) solution methods.}


\subsection{Contributions}

In this paper, we aim to derive the exact deterministic reformulations and solve the CC-SMCP to optimality. Our focus is on individual chance constraints with LHS uncertainty in set multicover problems,} where both decision and random variables are purely binary. Our contributions may be summarized as follows:
\begin{enumerate}
  \item Using techniques of enumerative combinatorics, discrete probability 
  distributions, and combinatorial optimization, we derive exact deterministic reformulations of CC-SMCP that are based on hierarchies of bounds, which
  are in turn utilized in construction of an outer-approximation (OA) algorithm for CC-SMCP. 
  \item {We investigate how to reduce the number of chance constraints by exploiting vector dominance relations, and we develop reformulations for two special cases of CC-SMCP by applying the ``log-transformation" method and by leveraging properties of the binomial distribution, respectively.}
  \item We present some theoretical results on SAA method to approximate the true optimal value of CC-SMCP under a finite discrete probability space. We also studied the importance sampling (IS) method and obtained a sufficient condition for selecting the optimal IS estimator.
\end{enumerate}
 
The rest of this paper is organized as follows. Section \ref{sec: reformulation}  develops an exact deterministic reformulation for CC-SMCP and proposes an outer-approximation algorithm to address CC-SMCP. Some presolving methods such as reducing the number of chance constraints are introduced as well. Reformulations for two special cases of CC-SMCP are investigated in Section \ref{sec: special_cases}.  Section \ref{sec: sampling} studies two sampling-based methods to approximate the optimal value of CC-SMCP. Section \ref{sec: experiments} presents some computational experiments to evaluate the effectiveness of the proposed model and solution approaches. Finally, Section \ref{sec: conclusion} concludes the paper.

\section{Deterministic Reformulations and Solution Approaches} \label{sec: reformulation}

\subsection{Equivalent reformulation}

In this section, we consider a deterministic reformulation of CC-SMCP using some combinatorial methods. For ease of exposition, we call the probability $\mathbb{P}\left[\sum_{j = 1}^n \tilde{a}_{ij}x_j \ge k_i\right]$ the \textit{cover probability} of the item $i$. Let $X_i$ be the set of all $x$ that satisfy the corresponding cover constraint for the $i$th item, i.e.,
\begin{equation*}
X_i = \left\{x\in \{0,1\}^n:\ \mathbb{P}\left[\ssum_{j = 1}^n \tilde{a}_{ij} x_j \ge k_i \right] \ge 1- \epsilon_i \right\} \label{eq: X},
\end{equation*} 
and thus the problem \eqref{eq: set-k-covering} can be expressed as 
\begin{equation}
\min\left\{\ssum_{j \in [n]} c_j x_j : x \in X_i, i=1, \ldots, m;~~ x \in B \subseteq \{0,1\}^n\right\}. \label{eq: set-k-covering-compact}
\end{equation}
We call the set $X_i$ the \textit{probabilistic covering set} of the item $i$. We use $i$ to index the item and $j$ to index the set throughout the paper. In what follows, for the sake of simplicity, we ignore the index $i$ and study the reformulation of $X_i$ (simplified as $X$) so that an equivalent deterministic reformulation of the form \eqref{eq: set-k-covering-compact} can be obtained for CC-SMCP.

{Throughout the paper we assume that $\tilde{a}_{j}$'s are mutually independent {Bernoulli random variables}} (note that $a_j$'s are not necessarily identically distributed). Let $S$ be a subset of $[n]$, and $S^c: = [n] \backslash S$, and define the events $A_{j}(x) := \{ \tilde{a}_j x_j  =1\}$ and $A_{j}^c(x) := \{ \tilde{a}_j x_j  =0\}$. We also define $A_{S}(x) := \cap_{j \in S}A_j(x)$ and $A_{S}^c(x) := [\cap_{j \in S}A_j(x)]^c$.  {To improve clarity, we summarize the necessary notation and definitions, including sets, events, parameters, decision variables in Table \ref{tab: notations}.}

\begin{table}[h!]
\footnotesize
  \centering
  \setlength{\abovecaptionskip}{0pt}%
  \setlength{\belowcaptionskip}{8pt}%
  {
  \caption{Sets, Parameters, Events and Decision Variables}
    \begin{tabular}{lp{0.75\textwidth}}
    \toprule
     Notation &  Definition \\
    \midrule
    Indices \& Sets  \\
    $[m]$ & Indices of $m$ items; we use $i \in [m]$ to denote the $i$-th item throughout the paper \\
    $[n]$ & Indices of $n$ sets, we use $j \in [n]$ to denote the $j$-th set throughout the paper \\
    $B$ & Set of deterministic constraints on the decision variables $x$ \\
     $X_i$ & Probabilistic covering set related to  item $i$ \\
     $X$ & Same as $X_i$ (ignoring the index $i$) \\
     $\Omega$ & Set of scenarios used in SAA and IS \\
    \hline
    Parameters \\
    $c \in \mathbb{R}^n$ & Cost coefficient vector \\
    $k_i \in \mathbb{Z}_+$ & Cover requirement for item $i$ \\ 
    $\epsilon \in (0, 1)$ & Prespecified allowed failure probability for item $i$ \\
    $\tilde{a}_{ij} \in \{0,1\}$ & Bernoulli random variable indicating whether item $i$ belongs to set $j$ \\
    $\tilde{a}_{j} \in \{0,1\}$ & Same as $\tilde{a}_{ij}$ (ignoring the index $i$)  \\
    \hline
    Variables \\
    $x_{j} \in \{0,1\}$ &  Equals 1 if the $j$-th set is selected, and $0$ otherwise\\
    $z_{i}(\omega) \in \{0,1\}$ & Equals 1 if scenario $\omega \in \Omega$ occurs for item $i$ in the SAA formulation\\
    \hline
    Events \\
    $A_{j}(x)$ & Event $\{ \tilde{a}_j x_j  =1\}$; $A_{j}^c(x)$ denotes its complement \\
    $A_{S}(x)$ & Event $\cap_{j \in S}A_j(x)$; $A^c_{S}(x)$ denotes its complement \\ 
    \bottomrule
    \end{tabular}} %
  \label{tab: notations}%
\end{table}%

A closed-form expression for the cover probability $\mathbb{P}\left[\sum_{j=1}^n \tilde{a}_j x_j  \ge k \right]$ can be developed using the inclusion-exclusion principle, stated as the following lemma:

\begin{lem}\label{lem: 1}
Let $\tilde{a}_{j}$ be {mutually independent} Bernoulli random variables with $\P(\tilde{a}_j =1)=p_j$ for each $j \in [n]$, $x_j$ be binary variables, and {$k \in \{1,\ldots,n-1\}$}. Then
\begin{equation}
\mathbb{P}\left[\ssum_{j=1}^n \tilde{a}_j x_j  \ge k \right]=\ssum_{\ell=k}^n(-1)^{\ell - k}\binom{\ell-1}{\ell -k} h_\ell(x), \label{eq: closed form}
\end{equation}
where
$h_\ell(x) =\sum_{\substack{S \subseteq [n]\\ |S| = \ell}} \mathbb{P}\left[A_S(x)\right] = \sum_{\substack{S \subseteq [n]\\ |S| = \ell}}\prod_{j \in S} x_j p_j, ~~ \ell = k, \ldots, n.$
\end{lem}

Thus, based on Lemma \ref{lem: 1}, the probabilistic covering set $X$ can be reformulated into the following deterministic 0-1 nonlinear set
\begin{equation}
X = \left\{x\in \{0,1\}^n:\ \ssum_{\ell=k}^n(-1)^{\ell - k}\binom{\ell-1}{\ell -k} \sum_{\substack{S \subseteq [n]\\ |S| = \ell}}\prod_{j \in S} x_j p_{j} \ge 1- \epsilon \right\}. \label{eq: X_1}
\end{equation} 
\begin{rem}[Relationship to the Poisson binomial distribution]
If we define $Z := \sum_{j = 1}^n \tilde{a}_j x_j$, then $Z$ is a random variable that follows Poisson binomial distribution with a collection of n independent yes/no experiments with different success probabilities $x_1p_1, x_2p_2,\ldots,x_n p_n$. Actually, from the proof of Lemma \ref{lem: 1}, we can obtain an explicit expression for its probability mass function. 
\end{rem}
\begin{rem}[The probability of the complement of the covering event]\label{rem: complement}
From the proof of Lemma \ref{lem: 1}, we can also obtain the probability 
    \begin{equation*}
      \mathbb{P}\left[\ssum_{j=1}^n \tilde{a}_j x_j  \le k -1 \right]=\ssum_{d = 0}^{k-1}\ssum_{\ell=d}^n(-1)^{\ell - d}\binom{\ell}{d} h_\ell(x).
    \end{equation*}
In Section \ref{sec: OA}, we will explore some lower and upper bounds for this probability and present an outer-approximation algorithm that utilizes these findings.
\end{rem} 

For more references on the inclusion–exclusion principle and combinatorial identities, {we refer the reader to \ref{sec-apx: i-e}, as well as the texts by \citet{Stanley_2011} and \citet{prudnikov1986integrals}.}

Now we consider the linearization of a {general monomial term:
\begin{equation}
y :=\textstyle\prod_{j \in S} x_j, \label{eq: nonlinear}
\end{equation}}
where $S \subseteq [n]$ and $x_j \in\{0,1\}, j \in S$. Since $x_j$  are binary variables, $y$ is also a binary variable. Typically, there are two commonly used approaches to linearize the cross product \eqref{eq: nonlinear}. The first approach establishes the equivalence of the nonlinear equation \eqref{eq: nonlinear} to two linear inequalities. 

\begin{lem}[\citet{glover1974converting}]  \label{lem: 2}
Let $s=|S|$. Equality \eqref{eq: nonlinear} holds if and only if
\begin{align}
\ssum_{j \in S} x_j-y \leq s-1, \quad
-\ssum_{j \in S} x_j+s y \leq 0, \quad
 x_j \in\{0,1\}, \; j \in S, \;\; y \in \{0,1\}.
\label{eq: linearization-I}
\end{align} 
\end{lem} 
If we allow the use of more constraints in exchange for relaxing the binary variable $y$ in \eqref{eq: linearization-I}, then we obtain the so-called \textit{standard linearization} \citep{glover1974converting}:

\begin{equation}
y=\sprod_{j \in S} x_j \Leftrightarrow \left\{
\ssum_{j \in S} x_j-y \leq s-1, \;\;
y \le x_j \; \forall j \in S, \;\;  
x_{j} \in \{0,1\}\; \forall j\in S, \;\; y \ge 0\right\}. \label{eq: linearization-II}
\end{equation}
Indeed, from the experimental results reported below in Section \ref{sec: experiments}, it can be suggested that \eqref{eq: linearization-II} may be preferable to \eqref{eq: linearization-I} in spite of the increased number of constraints.

To obtain a linearization of $X$ given by \eqref{eq: X_1}, let $T_\ell$ be a set of $\ell$ indices (out of $n$ indices),  $\sbinom{[n]}{\ell}$ be a collection of all sets of $\ell$ indices such that $\left|\sbinom{[n]}{\ell}\right| = \binom{n}{\ell}$. {For instance, $\sbinom{[3]}{2} = \{\{1,2\},\{1,3\},\{2,3\}\}$. Let} 
\begin{equation}
y_{_{T_\ell}} := \sprod_{j \in T_{\ell}}x_j, \quad T_\ell \in \sbinom{[n]}{\ell}, \quad\ell=k,\ldots,n. \label{eq: substitution}
\end{equation}
Taking the first linearization approach as an example (combining Lemma \ref{lem: 1} and \ref{lem: 2}), 
i.e., substituting \eqref{eq: substitution} into equation \eqref{eq: closed form} and adding constraints as defined in system \eqref{eq: linearization-I}, we obtain the following result: 
\begin{thm} \label{thm: ref_X}
The chance-constrained set 
\begin{equation*}
X = \left\{x\in \{0,1\}^n:\ \mathbb{P}\left[\ssum_{j = 1}^n \tilde{a}_{j} x_j \ge k \right] \ge 1- \epsilon \right\}
\end{equation*}
admits the following linearized deterministic reformulation:
\begin{subequations} \label{eq: ref_X}
\begin{align}
& \ssum_{\ell = k}^n (-1)^{\ell-k} \binom{\ell-1}{\ell-k}\sum_{T_\ell \in \binom{[n]}{\ell}}\left(\prod_{j \in T_{\ell}}p_j\right)y_{_{T_\ell}} \ge 1- \epsilon \
\\
& \ssum_{j \in T_{\ell}} x_j - y_{_{T_\ell}} \le \ell -1, \quad \forall T_\ell \in \binom{[n]}{\ell}, \; \ell = k,\ldots,n 
\\
& -\ssum_{j \in T_{\ell}} x_j + \ell y_{_{T_\ell}} \le 0, \quad \forall T_\ell \in \binom{[n]}{\ell}, \; \ell = k,\ldots,n 
\\
&x_j \in \{0,1\},  \ \forall j \in [n]; \; y_{_{T_\ell}} \in \{0,1\} \ \forall T_\ell \in \sbinom{[n]}{\ell}, \; \ell = k,\ldots,n 
\end{align}
\end{subequations} 
\end{thm} 

For small values of $n$ or when $k$ is very close to $n$, the original CC-SMCP problem \eqref{eq: set-k-covering-compact} can be directly solved using the aforementioned linearization techniques (an equivalent BIP). However, the above reformulation \eqref{eq: ref_X} has some significant drawbacks. First, the number of variables and constraints included in this reformulation may grow exponentially as the number of sets increases, making it very time-consuming to compute feasible solutions even for moderately sized problems. Second, the product of probabilities can be very small, and the summation of these products may lead to a substantial \textit{cumulative rounding error} (also known as \textit{accumulation of roundoff error}). Third, an optimal solution often assigns a value of one to numerous decision variables, making the traditional column generation or branch-and-price approach unsuitable for this reformulation.
Next, we will discuss several ways to alleviate these computational challenges. We first propose an outer-approximation algorithm in Section \ref{sec: OA}, and discuss how to reduce the number of chance constraints by using vector domination relations in Section \ref{sec: reduced_problem}. Then we investigate two special cases of CC-SMCP in Section \ref{sec: special_cases}, which allow for additional simplifications.

\subsection{Outer approximations} \label{sec: OA}

In this section, we propose an outer-approximation (OA) algorithm for CC-SMCP with small $k$. We first notice that the probabilistic covering set $X$ is equivalent to 
\begin{equation*}
X = \left\{x\in \{0,1\}^n:\ \mathbb{P}\left[\ssum_{j = 1}^n \tilde{a}_jx_j \le k -1\right] \le \epsilon \right\}, 
\end{equation*} 
which can be further reformulated into the following  0-1 nonlinear set (as mentioned in Remark \ref{rem: complement})
\begin{equation*}
X = \biggl\{x\in \{0,1\}^n:\ \ssum_{d = 0}^{k-1}\ssum_{\ell=d}^n(-1)^{\ell - d}\binom{\ell}{d}  \ssum_{\substack{S \subseteq [n]\\ |S| = \ell}}\prod_{j \in S} x_j p_j \le \epsilon \biggr\}.
\end{equation*} 
Next we derive some lower and upper bounds for the probability $\mathbb{P}\left[\sum_{j = 1}^n \tilde{a}_jx_j \le k -1\right]$. 

\begin{thm}\label{thm: bounds}
Let $\tilde{a}_{j}$ be {mutually independent Bernoulli random variables} with $\P(\tilde{a}_j =1)=p_j$ for each $j \in [n]$, $x_j$'s be binary variables, and {$k \in \{1,\ldots,n-1\}$.} Then
\begin{equation}
\mathbb{P}\left[\ssum_{j=1}^n \tilde{a}_j x_j  \le k -1 \right]=\left\{ 
\begin{aligned}
g_n(x) &\le g_t(x), \  \text{if $t$ is even} \\
g_n(x) &\ge g_t(x), \  \text{if $t$ is odd}
\end{aligned}
\right.
\end{equation}
where
\begin{equation*}
g_t(x):=\ssum\limits_{d = 0}^{k-1}{\ssum\limits_{\ell=d}^{\min\{t+d, n\}}}(-1)^{\ell - d}\binom{\ell}{d} \ssum\limits_{\substack{S \subseteq [n]\\ |S| = \ell}}\prod_{j \in S} x_j p_j, \ {t \in \mathbb{Z}_{\ge 0}}.
\end{equation*}
Moreover, if for any fixed $S \subseteq [n]$ with $|S| = d \le \ell < n$, we have 

\begin{equation}\label{condition}
\small
\ssum\limits_{\substack{S \subseteq T \subseteq [n]\\ |S| = d, |T| = \ell}} \mathbb{P}\left[A_T(x) \right] \ge \ssum\limits_{\substack{S \subseteq T' \subseteq [n]\\ |S| = d, |T'| = \ell + 1}} \mathbb{P}\left[A_{T'}(x)\right],
 \end{equation}
then the following hierarchy of inequalities holds:

\begin{align}
& \mathbb{P}\left[\ssum_{j=1}^n \tilde{a}_j x_j  \le k -1 \right]=g_n(x) \le \cdots \le g_4(x) \le g_2(x) \le g_0(x),  \label{eq: ub}
\\
& \mathbb{P}\left[\ssum_{j=1}^n \tilde{a}_j x_j  \le k -1 \right]=g_n(x) \ge \cdots \ge g_5(x) \ge g_3(x) \ge g_1(x). \label{eq: lb}
\end{align}
\end{thm}

\begin{figure}[htp]
\centering
\begin{tabular}{cc}
\includegraphics[width=.4\textwidth]{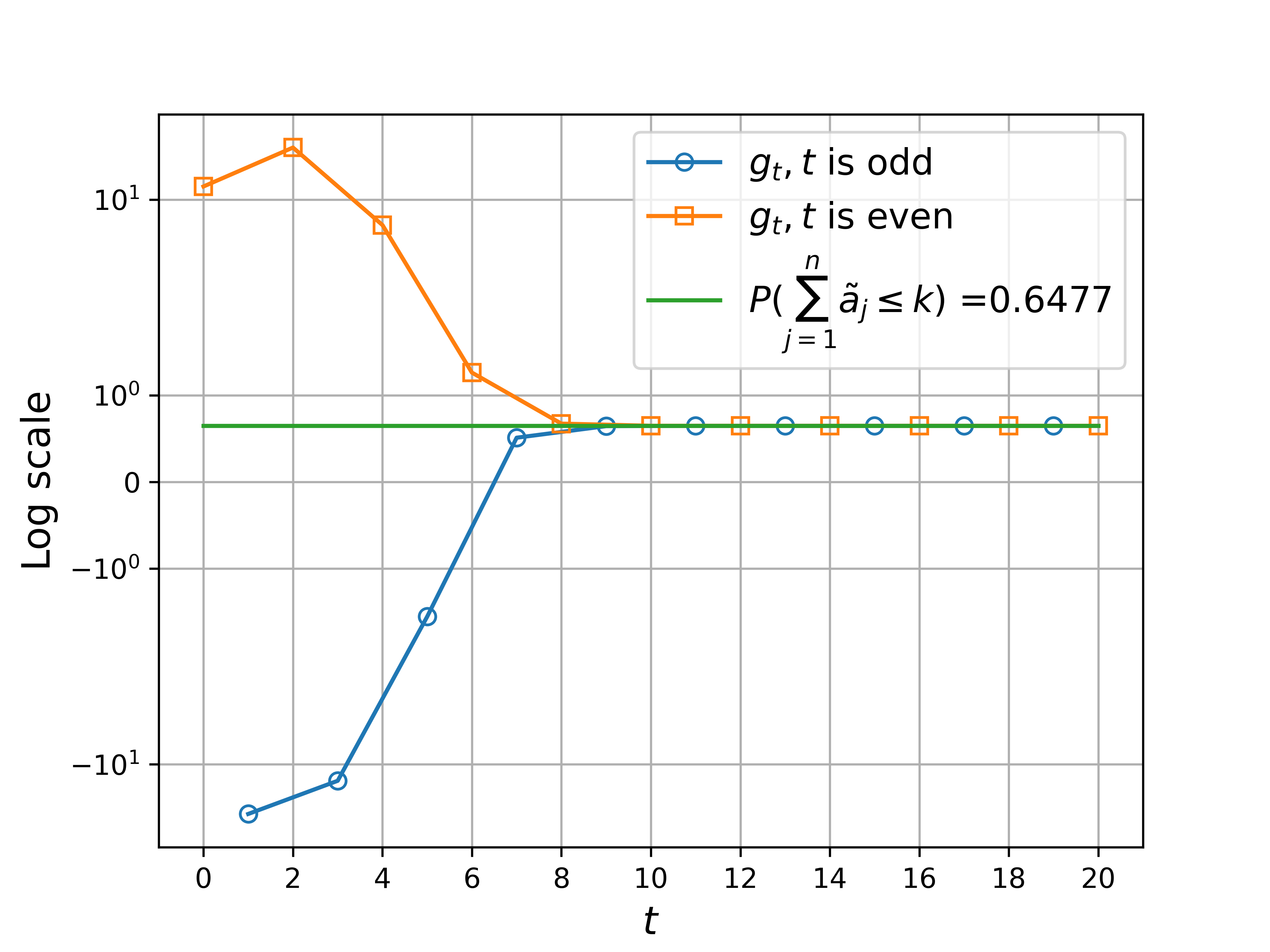} & \quad \includegraphics[width=.4\textwidth]{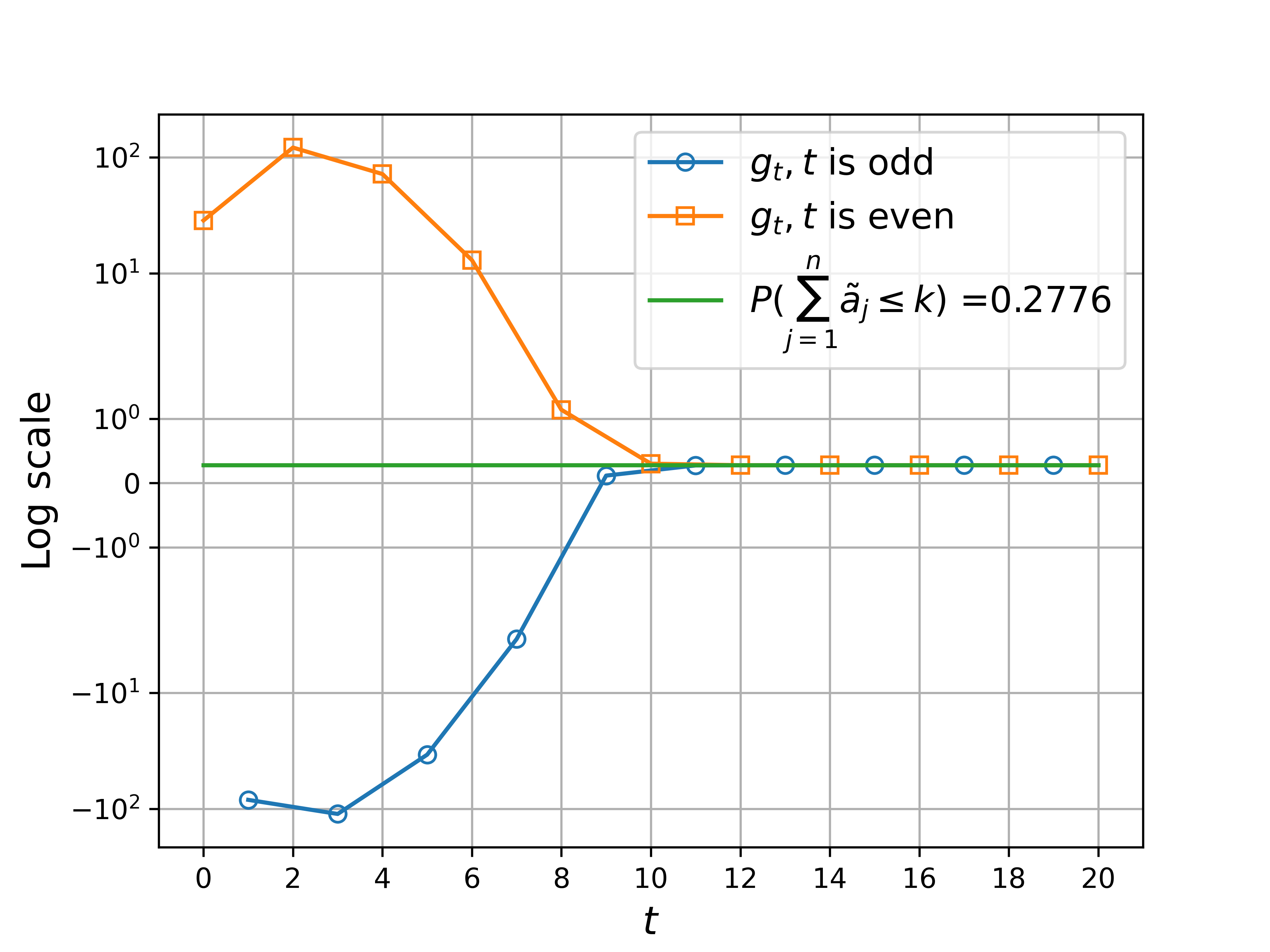}  \\
{\small (a) $n = 20, k = 3, \P(\tilde{a}_j=1)=0.15, \ \forall j \in [n].$}  & {\small (b) $n = 20, k = 3, \P(\tilde{a}_j=1)\in \{.1,.2,.3,.5\}$} \\
\end{tabular}
\caption{Convergence of bounds for the cover probability $\mathbb{P}\left[\sum_{j = 1}^n \tilde{a}_j \le k\right]$. }
\label{fig: convergence}
\end{figure}

Figure \ref{fig: convergence} illustrates how upper and lower bounds converge to the probability $\P[\sum_{j = 1}^n \tilde{a}_j \le k]$ as $t$ increases. In these experiments, we set $n = 20$, $k = 3$ and let $x_j = 1 $ for each $j \in [n]$. Specifically, in Figure \ref{fig: convergence}(a), we assume that random variables $\tilde{a}_j$ are i.i.d. with Bernoulli distribution with $\P(\tilde{a}_j=1)=0.15$, for each $j \in [n]$. In Figure \ref{fig: convergence}(b), we assume that random variables follow different Bernoulli distributions, such that $\P(\tilde{a}_j=1)=0.1, j \in \{1,\ldots,7\}$, $\P(\tilde{a}_j=1)=0.2, j \in \{8,\ldots,13\}$, $\P(\tilde{a}_j=1)=0.3, j \in \{14,\ldots,17\}$ and $\P(\tilde{a}_j=1)=0.5, j \in \{18,\ldots,20\}$. As can be seen in these two figures, the upper and lower bounds converge rapidly to the final probability, within about $4-5$ iterations. 

{Please note that condition \eqref{condition}, which guarantees the monotonicity of the bounds in \eqref{eq: ub} and \eqref{eq: lb}, is quite strong. As illustrated in Figure \ref{fig: convergence}, the sequence $g_t(x)$ typically does not exhibit monotonic behavior for smaller values of $t$.} However, it can be demonstrated that, after a sufficient number of iterations, the sequence $g_t(x)$ will eventually become monotonically increasing or decreasing for odd and even $t$, respectively.

We need additional notations in the following. Let $g_t(x) = \ssum\limits_{d=0}^{k-1} g_{t,d}(x)$, where $g_{t,d}(x) = \ssum\limits_{\ell=d}^{{\min\{t+d,n\}}}(-1)^{\ell - d}\binom{\ell}{d}  \ssum\limits_{\substack{S \subseteq [n]\\ |S| = \ell}}\prod_{j \in S} x_j p_j$. From the proof of Theorem \ref{thm: bounds}, for any $0 \le d\le k-1$, we have 
\begin{equation}
\mathbb{P}\left[\ssum_{j=1}^n \tilde{a}_j x_j  = d \right]=\left\{ 
\begin{aligned}
g_{n,d}(x) &\le g_{t,d}(x), \  \text{if $t$ is even} \\
g_{n,d}(x) &\ge g_{t,d}(x), \  \text{if $t$ is odd}
\end{aligned}
\right.
\end{equation}

\begin{lem}\label{lem:monotone}
  Given $d$ where $0\le d \le k-1$, if for some odd $t$, $g_{t,d}(x) \le g_{t+2,d}(x)$, then $g_{\tau,d}(x) \le g_{\tau+2,d}(x),$ for all odd $\tau$ such that $t \le \tau \le n-2$;
  Similarly, if for some even $t$, $g_{t,d}(x) \ge g_{t+2,d}(x)$, then $g_{\tau,d}(x) \ge g_{\tau+2,d}(x),$ for all even $\tau$ such that $t \le \tau \le n-2$.
\end{lem}

\begin{prop}\label{thm:monotone}
Let $t_{\min}^{odd} : =\min\{t\ |\ t\text{ is odd}, 1 \le t \le n-2, g_{t,d}(x) \le g_{t+2,d}(x), 0 \le d \le k-1\}$. Then $g_{\tau}(x) \le g_{\tau+2}(x)$ for all odd $\tau$ such that $t_{\min}^{odd} \le \tau \le n-2$;  
Similarly, let $t_{\min}^{even} : =\min\{t\ |\ t\text{ is even}, 0 \le t \le n-2, g_{t,d}(x) \ge g_{t+2,d}(x), 0 \le d \le k-1\}$. Then $g_{\tau}(x) \ge g_{\tau+2}(x),$ for all even $\tau$ such that $t_{\min}^{even} \le \tau \le n-2$.
\end{prop}

One of the popular methods of global optimization consists in approximating the given problem by a sequence of easier problems, such that the solutions of the approximating problems converge  
to an optimal solution of the given problem. If the feasible sets of the approximating problems contain the feasible set of the original problem, this method is known as \textit{outer approximation}. According to the above theorem, we can propose an outer-approximation algorithm by using the lower bounds to solve CC-SMCP to optimality.

Let $g_{t_i}^{(i)}(x):=\sum_{d = 0}^{k-1}\sum_{\ell^=d}^{t_i+d}(-1)^{\ell - d}\binom{\ell}{d}  \sum_{\substack{S \subseteq [n]\\ |S| = \ell}}\prod_{j \in S} x_{ij} p_{ij} $ and initialize $t_i: =1$ for each $i \in [m]$. We start by solving the following relaxed problem: 
\begin{equation} 
\nu:= \min\left\{\ssum_{j \in [n]} c_j x_j : g_{t_i}^{(i)}(x) \le \epsilon_i, \forall i \in [m], x \in B \subseteq \{0,1\}^n\right\} \label{eq: relaxed problem}.
\end{equation}
Clearly, $\nu$ provides a lower bound for the optimal value of problem \eqref{eq: set-k-covering-compact} by Theorem \ref{thm: bounds}. Let $\bar{x}$ be an optimal solution to problem \eqref{eq: relaxed problem}, and check whether $\bar{x} \in X_i$ for every $i \in [m]$ by using Lemma \ref{lem: 1}. If $\bar{x} \in X_i$ for each $i \in [m]$, then $\bar{x}$ is also an optimal solution to the chance-constrained problem \eqref{eq: set-k-covering-compact}. Otherwise, we have $\bar{x} \notin X_i$ for some $i \in [m]$, replace $t_i$ by $t_i + 2$ in $g_{t_i}^{(i)}(x)$ and solve an updated relaxed problem \eqref{eq: relaxed problem}, obtaining an updated solution $\bar{x}$. The process is repeated until $\bar{x} \in X_i$ for each $i \in [m]$. Algorithm \ref{alg: oa} summarizes this outer-approximation method for solving CC-SMCP.


\begin{algorithm}[H]
Initialize $t_i := 1$ for each $i \in [m]$\; 


Solve the relaxed  problem \eqref{eq: relaxed problem} to obtain a solution $\bar{x}$ and a lower bound $\nu$ for the original problem \eqref{eq: set-k-covering-compact}\; 


\uIf {$\bar{x} \in X_i$ for each $i \in [m]$}{\Return{the optimal solution $\bar{x}$ and the optimal value $\nu$}}
\ElseIf{$\bar{x} \notin X_i$ for some $i \in [m]$}{set $t_i : = t_i + 2$}

Go to line 2 and repeat this process until $\bar{x} \in X_i$ for every $i \in [m]$.

\caption{Outer Approximation Algorithm} \label{alg: oa}
\end{algorithm} 


It is noteworthy that in each iteration, new constraints as well as decision variables are added to the relaxed problem \eqref{eq: relaxed problem}. That is, both column and constraint generation procedures are included in our algorithm. {At each iteration, we use Lemma 1 to verify whether the current solution satisfies the chance constraints. If it does not, we refine $g_t(x)$ to better approximate the true probability function and repeat. We continue this process until we obtain a solution that satisfies all chance constraints, and the algorithm terminates after a finite number of iterations.}

{
\begin{rem}
  \textit{Inner approximations.} Similar to the outer-approximation procedure, using the upper bounds that we have derived in Theorem \ref{thm: bounds}, we can propose an inner-approximation algorithm considering that $t_{i}$ is an even number in problem \eqref{eq: relaxed problem}. However, this method has some limitations. Two of the main challenges with this method are that, the upper bounds it provides are often trivial (larger than one) in the first few iterations, which may not lead to a feasible solution. Additionally, for inner approximations, checking the feasibility or optimality of the original problem is very challenging. In other words, it is difficult to determine a stopping criterion when the original problem is infeasible.
\end{rem}}

\subsection{Reducing the number of chance constraints} \label{sec: reduced_problem}

In this section, we propose an approach to reducing the number of constraints in the chance-constrained problem \eqref{eq: set-k-covering-compact}, using vector dominance relations defined in a partially ordered set (also called poset). We first consider problem \eqref{eq: set-k-covering-compact} under the assumption that for any fixed $j \in [n]$, $P(\tilde{a}_{ij} =1) = p_j\ \forall i \in [m]$, i.e., {the set indexed by $j$} covers each item $i \in [m]$ with equal probability $p_j$. 
Then, the following lemma provides an approach to reduce the number of chance constraints in \eqref{eq: set-k-covering-compact}. 

\begin{lem} \label{lem: k-e}
    If {for any fixed $j \in [n]$}, $\P(\tilde{a}_{ij} =1) = p_j,\ \forall i \in [m]$
    and $i_1, i_2 \in [m]$ are two rows with $k_{i_1} \ge k_{i_2}$ and $\epsilon_{i_1} \le \epsilon_{i_2}$, then $X_{i_1} \cap X_{i_2} = X_{i_1}$, where the probabilistic covering sets $X_{i_1}$ and $X_{i_2}$ correspond to the $i_1$-th and ${i_2}$-th items, respectively.
\end{lem}

For a more general setting, we define a partially ordered set by considering the relation between the coefficients of the constraints in the chance-constrained problem. 

Let $P(\tilde{a}_{ij} =1) =: p_{ij}$ for each $i \in [m]$ and $j \in [n]$. Given two vectors $v^{i_1}:= (p_{i_11},p_{i_12},\ldots,p_{i_1n},-k_{i_1},\epsilon_{i_1})$
and $v^{i_2}:= (p_{i_21}, p_{i_22},\ldots,p_{i_2n},-k_{i_2},\epsilon_{i_2})$ of length $n+2$ corresponding to the $i_1$-th and $i_2$-th rows of the chance constraints \eqref{eq: set-k-covering-chance_constraints}, the relation $v^{i_1} \le v^{i_2}$ means that  $v^{i_1}_t \le v^{i_2}_t$ for every $t = 1,\ldots, n+2$, whereas the relation $v^{i_1} \nleq v^{i_2}$ means that there exists at least one index $t$ such that $v^{i_1}_t > v^{i_2}_t$. The relations ``$\ge$'' and ``$\ngeq$'' can also be defined similarly. We say $v^{i_1}$ and $v^{i_2}$ are \textit{comparable} if $v^{i_1} \le v^{i_2}$ or $v^{i_2} \le v^{i_1}$. Otherwise they are incomparable. Let $\mathcal{C} : = \{v^1,v^2,\ldots,v^m\}$ be the set of all vectors corresponding to the associated chance constraints. Clearly, the relation ``$\le$'' is a \textit{partial order} on the set $\mathcal{C}$, as it is reflexive, antisymmetric, and transitive. Then the following result generalizes Lemma \ref{lem: k-e}: 
\begin{lem} \label{lem: relation}
    If there exist two items $i_1,i_2 \in [m]$ such that $v^{i_1} \le v^{i_2}$, then $X_{i_1} \cap X_{i_2} = X_{i_1}$, where the probabilistic covering sets $X_{i_1}$ and $X_{i_2}$ correspond to the $i_1$-th and ${i_2}$-th items, respectively.
\end{lem} 
We consider minimal and maximal elements of the partially ordered set  $(\mathcal{C}, \le)$. An element $v \in \mathcal{C}$ is a \textit{minimal} element if there is no element $u \in \mathcal{C}$ and $u \neq v$ such that $u \le v$. An element $v \in \mathcal{C}$ is a \textit{maximal} element if there is no element $u \in \mathcal{C}$ and $u \neq v$  such that $v \le u$. Based on Lemma \ref{lem: relation}, the chance-constrained problem \eqref{eq: set-k-covering-compact} can be further simplified by only focusing on the chance constraints whose rows correspond to the minimal elements of $\mathcal{C}$.

\begin{prop}
    Let $\mathcal{C}' \subseteq \mathcal{C}$ be the set of all minimal elements of $\mathcal{C}$, and $I \subseteq [m]$ be the rows corresponding to the set $\mathcal{C}'$. Then the chance-constrained problem \eqref{eq: set-k-covering-compact} is equivalent to the following problem
\begin{equation*}
\min\left\{\ssum_{j \in [n]} c_j x_j : x \in X_i, \forall i \in I, x \in B \subseteq \{0,1\}^n\right\} 
\end{equation*}
\end{prop}

It can be readily seen that the worst-case complexity of finding all minimal elements of the poset $(\mathcal{C}, \le)$ is $O(m^2n)$, {since the number of comparisons between elements is $O(m^2)$, and comparing every two elements takes $O(n+2) = O(n)$.}   


\subsection{Presolving method}

Based on the analyses conducted in the previous and following sections, we can propose a presolving method to simplify chance-constrained problem and improve the efficiency of our algorithm. First, we use vector dominance relations to reduce the problem's complexity by identifying and eliminating redundant chance constraints in the problem \eqref{eq: set-k-covering-compact} (see details in Section \ref{sec: reduced_problem}). Denote the collection of remaining items by $I$. If $k_i = 1$ for some item $i \in I$, then use the ``log-transformation'' technique to reformulate the probabilistic covering set $X_i$ (see details in Section \ref{sec: case 1}). If the item $i$ is covered by each set $j \in [n]$ with equal probability $p_i$ for some $i \in I$, then reformulate the probabilistic covering set $X_i$ using the method mentioned in Theorem \ref{thm: ref} (see details in Section \ref{sec: case 2}). Finally, we solve the simplified problem using the outer-approximation algorithm (see details in Section \ref{sec: OA}).

To sum up, we first attempt to reduce the number of chance constraints by only considering the set of all minimal vectors corresponding to the associated chance constraints. Further, we employ two distinct reformulation techniques (as stated in Section \ref{sec: special_cases}) to handle two specific cases. These reformulation techniques will be more computationally efficient, as they involve significantly fewer variables and constraints compared to those in Theorem \ref{thm: ref_X}. Additionally, it is worth mentioning that at the final step, in the process of checking feasibility of our solution, we can utilize the discrete Fourier transform (DFT) method to compute the PMF and CDF of Poisson binomial distribution instead of using Lemma \ref{lem: 1} directly. The DFT method provides us with an effective way to mitigate the cumulative rounding error. For more information about the DFT method, please refer to \cite{hong2013computing}.

\section{Reformulation for Two Special Cases of CC-SMCP} \label{sec: special_cases}

\subsection{Reformulation for probabilistic set covering problem} \label{sec: case 1}

Suppose $k_i = 1$ for $i=1,\ldots,m.$ As mentioned before, when $k_i=1$ for each $i \in [m]$, the problem \eqref{eq: set-k-covering} becomes the probabilistic set covering (PSC) problem, which has already been well studied by \citet{haight2000integer}, \citet{fischetti2012cutting}, and \citet{ahmed2013probabilistic}. In the independent case, the probabilistic covering set $X$ can be further simplified using ``log-transformation''. 
We include it here for the sake of completeness. 

\begin{prop}[\cite{haight2000integer,fischetti2012cutting,ahmed2013probabilistic}]
    Let $\tilde{a}_{j}$ be independent random variable with $\P(\tilde{a}_j =1)=p_j$ for each $j \in [n]$, $x_j$'s be binary variables, and $k =1$ in the probabilistic covering set $X$. Then
\begin{equation}
X = \left\{x\in \{0,1\}^n:\ \ssum_{j \in [n]}\log (1-p_j) x_j \le \log \epsilon \right\}. \label{eq: log-X}
\end{equation}
\end{prop}

It follows that, when $k_i = 1$ for each $i \in [m]$, CC-SMCP is equivalent to the following deterministic binary linear program (BIP): 
\begin{equation*}
\min \left\{\ssum_{j \in [n]} c_j x_j : \ssum_{j \in [n]}\log (1-p_{ij}) x_{j} \le \log \epsilon_i,  \ \forall i \in [m],\ x \in B \subseteq \{0,1\}^n \right\}.
\end{equation*}

\subsection{Case when cover probabilities of each item are equal} \label{sec: case 2}

Assume that for any fixed item $i \in [m]$, we have $\P(\tilde{a}_{ij} =1) = p_i,\ \forall j \in [n]$. That is, the item $i$ is covered by each set $j \in [n]$ with equal probability $p_i$. We first notice that in this case, the value of the probability $\mathbb{P}\left[\sum_{j=1}^n \tilde{a}_{ij} x_j  \ge k_i \right]$ is only related to the number of selected sets in the solution. In this case, given a candidate solution $\bar{x} \in \{0,1\}^n$, the following lemma provides an efficient way to compute the probability $\mathbb{P}\left[\sum_{j=1}^n \tilde{a}_{j} \bar{x}_j  \ge k \right]$ rather than using Lemma \ref{lem: 1} directly: 
\begin{lem}\label{lem: covering-1}
Let $\tilde{a}_{j}$ be {mutually independent Bernoulli random variables} with $\P(\tilde{a}_j =1)=p$ for each $j \in [n]$, $\bar{x} \in \{0,1\}^n$ with $\sum_{j \in [n]} \bar{x}_j = d$, and {$k \in \{1,\ldots,n-1\}$.}  Then
\begin{equation}\label{eq: case-2.0}
\mathbb{P}\left[\ssum_{j=1}^n \tilde{a}_j \bar{x}_j  \ge k \right]=\ssum_{\ell = k}^d\binom{d}{\ell}p^{\ell}(1-p)^{d-\ell}.
\end{equation}
\end{lem} 

The above result provides an efficient way to compute the cover probability $\mathbb{P}\left[\sum_{j=1}^n \tilde{a}_j \bar{x}_j  \ge k \right]$ without knowing the exact value of $\bar{x}_j$, and it only depends on the value of the summation $\sum_{j\in [n]}{\bar{x}_j}$. However, it has a severe drawback:  due to the computer's inability to represent some numbers exactly, the calculation may result in a large cumulative rounding error that cannot be ignored, especially when $d$ is a large number. Indeed, the calculation \eqref{eq: case-2.0} including a series of summation operations of very small numbers can lead the cumulative error to be especially problematic. Even if the rounding error from a single operation (addition, product) is small, the cumulative error from many operations may be significant. 
Taking the above considerations into account, we develop an alternative equivalent explicit formula to compute the cover probability. This formula can be readily implemented in practice, as numerous software packages offer efficient methods for calculating the values of certain special functions. 
Specifically, we utilize the 
beta function $\mathrm{B}(a,b)$  and 
the incomplete beta function $\mathrm{B}(q;a,b)$:
\begin{equation*}
\mathrm {B} (a,b)=\textstyle\int _{0}^{1}t^{a-1}\,(1-t)^{b-1}\,dt \;\text{ and }\; \mathrm {B} (q;\,a,b)=\int _{0}^{q}t^{a-1}\,(1-t)^{b-1}\,dt.
\end{equation*}
Then, the following lemma provides another computational method that reduces the cumulative rounding error to obtain the cover probability: 

\begin{lem}\label{lem: covering-2}
Let $\tilde{a}_{j}$ be {mutually independent Bernoulli random variables} with $\P(\tilde{a}_j =1)=p$ for each $j \in [n]$, $\bar{x} \in \{0,1\}^n$ with $\sum_{j \in [n]} \bar{x}_j = d$, and {$k \in \{1,\ldots,n-1\}$.}  Then
\begin{equation}
\mathbb{P}\left[\ssum_{j=1}^n \tilde{a}_j \bar{x}_j  \ge k \right]=k \sbinom{d}{k} \ssum_{\ell^=k}^d(-1)^{\ell - k}\binom{d-k}{\ell -k} \frac{p^\ell}{\ell}=k\binom{d}{k}\int_{0}^pu^{k-1}(1-u)^{d-k}du, \label{eq: case-2}
\end{equation}
or equivalently,
\begin{equation}
\mathbb{P}\left[\ssum_{j=1}^n \tilde{a}_j \bar{x}_j  \ge k \right]=\sbinom{d}{k} p^k {}_{2}F_{1}(k-d,k;k+1;p) = I_{p}(k,d-k+1),
\end{equation}  
where 
${}_{2}F_{1}(a,b;c;z)=\sum _{n=0}^{\infty }{\frac {a^{(n)}b^{(n)}}{c^{(n)}}}{\frac {z^{n}}{n!}}=1+{\frac {ab}{c}}{\frac {z}{1!}}+{\frac {a(a+1)b(b+1)}{c(c+1)}}{\frac {z^{2}}{2!}}+\cdots,$ $|z| < 1$, is the \textit{hypergeometric function}  
and $I_q(a,b) ={\frac {\mathrm {B} (q;\,a,b)}{\mathrm {B} (a,b)}}$ is the \textit{regularized incomplete beta function.} 
\end{lem}

Next we derive an equivalent reformulation for the probabilistic covering set $X$ based on the above results. Considering CC-SMCP with a single chance-constraint: 
\begin{equation}
\min\left\{\ssum_{j \in [n]} c_j x_j : x \in X,\  x \in B \subseteq \{0,1\}^n\right\},  \label{eq: set-k-covering-compact-single}
\end{equation} 
we explore properties of a solution to problem \eqref{eq: set-k-covering-compact-single} when $P(\tilde{a}_{j} =1) = p$, for each $j \in [n]$: 

\begin{lem}\label{lem: feasible-sol}
    Let $\tilde{a}_{j}$ be {mutually independent Bernoulli random variables} with $\P(\tilde{a}_j =1)=p$ for each $j \in [n]$. If $\bar{x}$ is a feasible solution to the problem \eqref{eq: set-k-covering-compact-single} with $\sum_{j= 1}^n\bar{x}_j=d$, then for any $\hat{x} \in B \subseteq \{0,1\}^n$ with $\sum_{j= 1}^n\hat{x}_j=\hat{d} \ge d$, $\hat{x}$ is also a feasible solution to the problem \eqref{eq: set-k-covering-compact-single}.
\end{lem}

Lemmas \ref{lem: covering-2} and \ref{lem: feasible-sol} imply the following result: 
\begin{thm} \label{thm: ref}
Let $\tilde{a}_{j}$ be {mutually independent Bernoulli random variables} with $\P(\tilde{a}_j =1)=p$ for each $j \in [n]$. Then $x$ is a feasible solution to problem \eqref{eq: set-k-covering-compact-single} if and only if $x$ is a feasible solution to the following deterministic covering problem:
\begin{equation*}
\min\left\{\ssum_{j \in [n]} c_j x_j : \ssum_{j \in [n]}x_j \ge \bar{d}, \; x \in B \subseteq \{0,1\}^n\right\}, 
\end{equation*}
where 
\begin{equation} \label{eq: subproblem}
\bar{d} := \min\left\{d \in \mathbb{Z}_+ : I_{p}(k,d-k+1) \ge 1 - \epsilon, k \le d \le n \right\}.
\end{equation}
That is, the probabilistic covering set $X$ admits the deterministic equivalent reformulation
\begin{equation*}
X = \left\{x \in \{0,1\}^n : \ssum_{j \in [n]}x_j \ge \bar{d}\right\}.
\end{equation*}
\end{thm} 
\begin{proof}
The necessity follows directly from Lemma \ref{lem: covering-2}, and the sufficiency from Lemma \ref{lem: feasible-sol}.
\end{proof}

For any feasible solution $\bar{x}$ to problem \eqref{eq: set-k-covering-compact-single}, a simple (but potentially weak) lower bound for the sum $\sum_{j \in [n]}\bar{x}_j$ can be obtained using Markov's inequality. 
\begin{prop}\label{prop: markov}
    Let $\tilde{a}_{j}$ be {mutually independent Bernoulli random variables} with $\P(\tilde{a}_j =1)=p$ for each $j \in [n]$, and $\bar{x}$ be a feasible solution to the chance-constrained  problem \eqref{eq: set-k-covering-compact-single} with $\sum_{j= 1}^n\bar{x}_j=d$, then 
    \begin{equation*}
    d : = \ssum_{j \in [n]}\bar{x}_j  \ge \max\left\{k, \left\lceil\frac{k(1-\epsilon)}{p}\right\rceil \right\}.
    \end{equation*}
\end{prop}
Actually, from Lemma \ref{lem: covering-2} and \ref{lem: feasible-sol}, we can notice that for any fixed $k$ and $p$, the function $I_{p}(k,d-k+1)$ is monotone non-decreasing with respect to $d$. Thus, the tightest lower bound for $d := \sum_{j \in [n]} \bar{x}_j$ (i.e., the optimal value of problem \eqref{eq: subproblem}) can be obtained by the binary search. Based on the above analysis, we can provide the following algorithm to solve CC-SMCP when cover probabilities of each item are equal: we first use a binary-search algorithm to obtain $\bar{d}_i$ for each item $i \in [m]$ where 
\begin{equation*} 
\bar{d}_i := \min\left\{d_i \in \mathbb{Z}_+ : I_{p_i}(k_i,d_i-k_i+1) \ge 1 - \epsilon_i, k_i \le d_i \le n \right\},
\end{equation*}
and choose the maximum $d_i$ to satisfy all the chance constraints, i.e., let $\bar{d}:= \max\{\bar{d}_i,\ i \in [m]\}$. Finally, we solve the following deterministic covering problem to optimality:
\begin{equation*}
\min\left\{\ssum_{j \in [n]} c_j x_j : \ssum_{j \in [n]}x_j \ge \bar{d} , x \in B \subseteq \{0,1\}^n\right\}.\;
\end{equation*}

\begin{rem}
    Note that if $B = \{0,1\}^n$, then we can obtain an explicit formula for the optimal value of the chance-constrained problem \eqref{eq: set-k-covering-compact} when cover probabilities of each item are equal to each other. Suppose that the sequence $c'_1,\ldots, c'_n$ is an order of $c_1,\ldots,c_n$ satisfying $c'_1 \le c'_2 \le \cdots \le c'_n$. Based on the above algorithm, the optimal value of the chance-constrained problem \eqref{eq: set-k-covering-compact} can be written as $\sum_{j =1}^{\bar{d}}c'_j.$
\end{rem}


\section{Sampling-Based Approaches}\label{sec: sampling}

In the following sections, we discuss two sampling-based methods, sample average approximation (SAA) and importance sampling (IS), to solve CC-SMCP on a finite discrete probability space. {The high-level ideas are inspired by previous works \citep{luedtke2008sample} and \citep{barrera2016chance}, which establish statistical relationships between chance-constrained problems and their reformulations, and prove the asymptotic convergence of the SAA and IS approaches. However, our setting focuses on chance-constrained combinatorial optimization problems with individual constraints over a finite discrete probability space. As is often the case with sampling-based techniques, analyzing these approaches relies on exploiting the specific structure of the problem. For CC-SMCP, we provide a convergence guarantee for the SAA approach and a sufficient condition for selecting the optimal IS estimator.
}



\subsection{Sample average approximation} \label{sec: SAA}

In this section, we consider CC-SMCP on a finite discrete probability space $(\Omega, 2^{\Omega}, \mathbb{P}_N)$ where $\Omega = \{ \omega_1,\ldots,\omega_N \}$ and $p(\omega_i) = \mathbb{P}_N(\omega = \omega_i)$. The sample average approximation (SAA) method approximates the true distribution via a finite empirical distribution, $\mathbb{P}_N$. We assume that the samples are independent and identically distributed (i.i.d.) with $p(\omega) = 1/N, \omega \in \Omega$. The SAA formulation of CC-SMCP is 
\begin{subequations}
        \begin{align}
         \min \quad &  \mathbf{c}^T \mathbf{x} \label{eq: SAA-big-M-0}\\
        \text { s.t. }\quad  & \tilde{A}_{i}(\omega) \mathbf{x} + M_i(\omega)(1-z_i(\omega)) \ge k_i, \ \forall \omega \in \Omega, \forall i \in [m] \label{eq: SAA-big-M-1}\\
        &{\textstyle\frac{1}{N}}\ssum_{\omega \in \Omega} z_i(\omega) \ge 1 - \alpha_i, \ \forall i \in [m] \label{eq: SAA-big-M-2}\\
        & \mathbf{z}_i \in \{0,1\}^{N}, \ \forall i \in [m] \label{eq: SAA-big-M-3}\\
        &\mathbf{x} \in B \subseteq \{0,1\}^n \label{eq: SAA-big-M-4}
        \end{align} \label{eq: SAA-big-M}
\end{subequations} 
where $\tilde{A}_{i}(\omega)$ is the $i$-th row of the coefficient matrix $\tilde{A}$ in scenario $\omega$, $M_i(\omega)$ is the big-M coefficient that guarantees feasibility whenever any element in $z_i(\omega)$ is equal to zero, and $\alpha_i \in [0,1)$ is the risk level that may be different from $\epsilon_i$ in CC-SMCP. 

\begin{rem}
    Generally, SAA is a sampling-based method commonly employed to approximate the function/distribution $f(\mathbf{x},\mathbf{\xi})$ that cannot be observed or computed directly. Here we choose to use the term ``SAA'' since the idea of our method is the same as SAA. 
\end{rem}

Often, the linear programming (LP) relaxation of the SAA reformulation \eqref{eq: SAA-big-M} is very weak because of the big-M coefficients that are introduced to model the chance constraints. Fortunately, note that in CC-SMCP problem, we have $\tilde{A} \in \{0,1\}^{m \times n}$, $\mathbf{x} \in \{0,1\}^n$ and $k_i \in \mathbb{Z}_+$. {In this case, we can choose the big-M coefficient in the SAA reformulation as $M_i(\omega) : = k_i$ for each $\omega \in \Omega$:} 
\begin{equation}
\textbf{[SAA]} ~~~ \min \left\{\mathbf{c}^T \mathbf{x}: \tilde{A}_{i}(\omega) \mathbf{x} \ge k_iz_i(\omega), \forall \omega \in \Omega, \forall i \in [m], \eqref{eq: SAA-big-M-2} - \eqref{eq: SAA-big-M-4} \right\}.     \label{eq: SAA-big-M-free}
\end{equation} 

In addition, we provide a convergence guarantee for the SAA approach. Details are provided in \ref{sec-apx: SAA}.



\subsection{Importance sampling}

Importance sampling (IS) is probably the most popular approach to reduce the variance of an estimator. In the case of rare event estimation, this also means increasing the occurrence of rare events.  In this paper, we consider an IS formulation based on CC-SMCP as an example and propose a sufficient condition for selecting the optimal IS estimator. We replace the chance-constraint \eqref{eq: set-k-covering-chance_constraints} by 
\begin{subequations}
\begin{align}
&\hat{A}_{i}(\omega) \mathbf{x} \ge k_iz_i(\omega), \forall \omega \in \Omega, \forall i \in [m] \\
&\ssum_{\omega \in \Omega} L(\hat{A}_i(\omega))\left(1-z_i(\omega)\right) \le N\epsilon_i, \forall i \in [m] \label{eq: IS_likelihood_cons}\\
&z_i(\omega) \in \{0,1\}, \forall \omega \in \Omega
\end{align} \label{eq: IS_reformulation}
\end{subequations}
where the definitions of $\hat{A}_i (\omega)$ and $z_i(\omega)$ are similar to those in the SAA reformulation \eqref{eq: SAA-big-M-free}, and $L(\cdot)$ is the likelihood ratio $L(\hat{a}) = \prod_{j = 1}^n (\frac{p_j}{\hat{p}_j})^{\hat{a}_j}(\frac{1- p_j}{1-\hat{p}_j})^{1- \hat{a}_j}$ with $p_j = P(\tilde{a}_j =1)$ and $\hat{p}_j = P(\hat{a}_j = 1)$ for each $j \in [n]$. In this case, $L$ is the ratio between the respective probability mass functions, since both $\tilde{a}$ and $\hat{a}$ have discrete support. 

We provide some theoretical results based on the above reformulation. The details of this part can be found in \ref{sec-apx: IS}.



\section{Numerical Experiments}\label{sec: experiments}

In this section, we study the computational performance of our proposed methods for solving the CC-SMCP. We compare the outer-approximation (OA) algorithm with two different linearization techniques with the sampling-based algorithms, i.e., the SAA and IS approaches, in terms of both time and effectiveness. The detailed experimental settings are listed at the beginning of the following subsections.

{All experiments are conducted using Python 3.9 with the optimization solver Gurobi 10.0.2 on a single node from the University of Arizona's HPC puma cluster \citep{hpc_cluster}. This node is equipped with one AMD EPYC 7642 48-core processor (2.4GHz) and {240 GB RAM (5GB per core)} and runs on Linux CentOS 7.} For all experiments, we use a time limit of 3,600 seconds and report the running time in CPU seconds. All other control parameters remain at their default settings within the branch-and-cut framework of Gurobi. Overall, in our experiments, we compare the performance of the following four solution approaches:

\begin{itemize}
    \item OA-I: outer approximation approach with the first linearization technique (Lemma \ref{lem: 2}) and presolving process applied (i.e., presolving + OA + \blue{linearization}-I)
    \item OA-II: outer approximation approach with the second linearization technique (standard linearization \eqref{eq: linearization-II}) and presolving process applied (i.e., presolving + OA + \blue{linearization}-II)
    \item SAA: sample average approximation approach (with reformulation \eqref{eq: SAA-big-M-free})
    \item IS: importance sampling approach (with reformulation \eqref{eq: IS_reformulation})
\end{itemize}

{
  \subsection{Case study: a real-world facility location problem}
  We evaluate the proposed methods on a real-world facility location problem derived from the \texttt{spopt} facility-location dataset describing an area of census tract 205 in San Francisco \citep{spopt2022}. The dataset contains 205 demand points, 16 candidate facility sites, and a network-distance matrix between each demand--facility pair, together with geographic information used for visualization. In our setting, the objective is to determine the minimum number of facilities such that each demand point \(i\) is covered at least \(k_i\) times with high probability. We set $k_i:=2$ for all $i$. To model uncertain coverage, we construct the probability matrix  $(p_{ij})$ based on the distance $d_{ij}$ between demand point $i$ and candidate facility $j$. Specifically, we define 
  \begin{equation*}
    p_{ij} = \left\{ 
\begin{aligned}
&\frac{1}{1+\exp(\frac{d_{ij}-10000}{s})}, \ &d_{ij} \le 15000\text{m} \\
&0, \  &d_{ij} >15000\text{m}
\end{aligned}
\right.
\end{equation*}  
where $s = 5000/\ln 9$. With this choice of \(s\), the resulting probabilities satisfy $p_{ij}(d_{ij}=5000\text{m}) = 0.9, p_{ij}(d_{ij}=10000\text{m}) = 0.5$ and $p_{ij}(d_{ij}=15000\text{m}) = 0.1$. This definition captures the desired property that the covering probability decreases as the distance increases.

\blue{We emphasize that this case study uses real spatial and network-distance data with artificially defined covering/coverage probabilities, rather than probabilities calibrated from historical coverage data. The probability function is chosen because it provides a smooth and monotone decrease in coverage reliability with distance. The parameter values are selected so that the coverage probabilities equal $0.9$, $0.5$, and $0.1$ at distances of 5,000, 10,000, and 15,000 meters, respectively, representing high, moderate, and low coverage reliability. Coverage is set to zero beyond 15,000 meters to impose a maximum effective service distance.}

Figure \ref{fig: real-world} illustrates the optimal facility locations returned by OA and SAA methods for $\epsilon_i = 0.1$.  While the two methods may identify different optimal facilities, both yield geographically balanced solutions. A detailed comparison of the four methods on different risk parameters is presented in Table \ref{tab:real-world}. As $\epsilon_i$ increases, the number of optimal facilities decreases from $7$ to $4$, reflecting the less coverage requirement under a higher allowable  risk level.

\begin{figure}[htp]
\centering
\includegraphics[width=.8\textwidth]{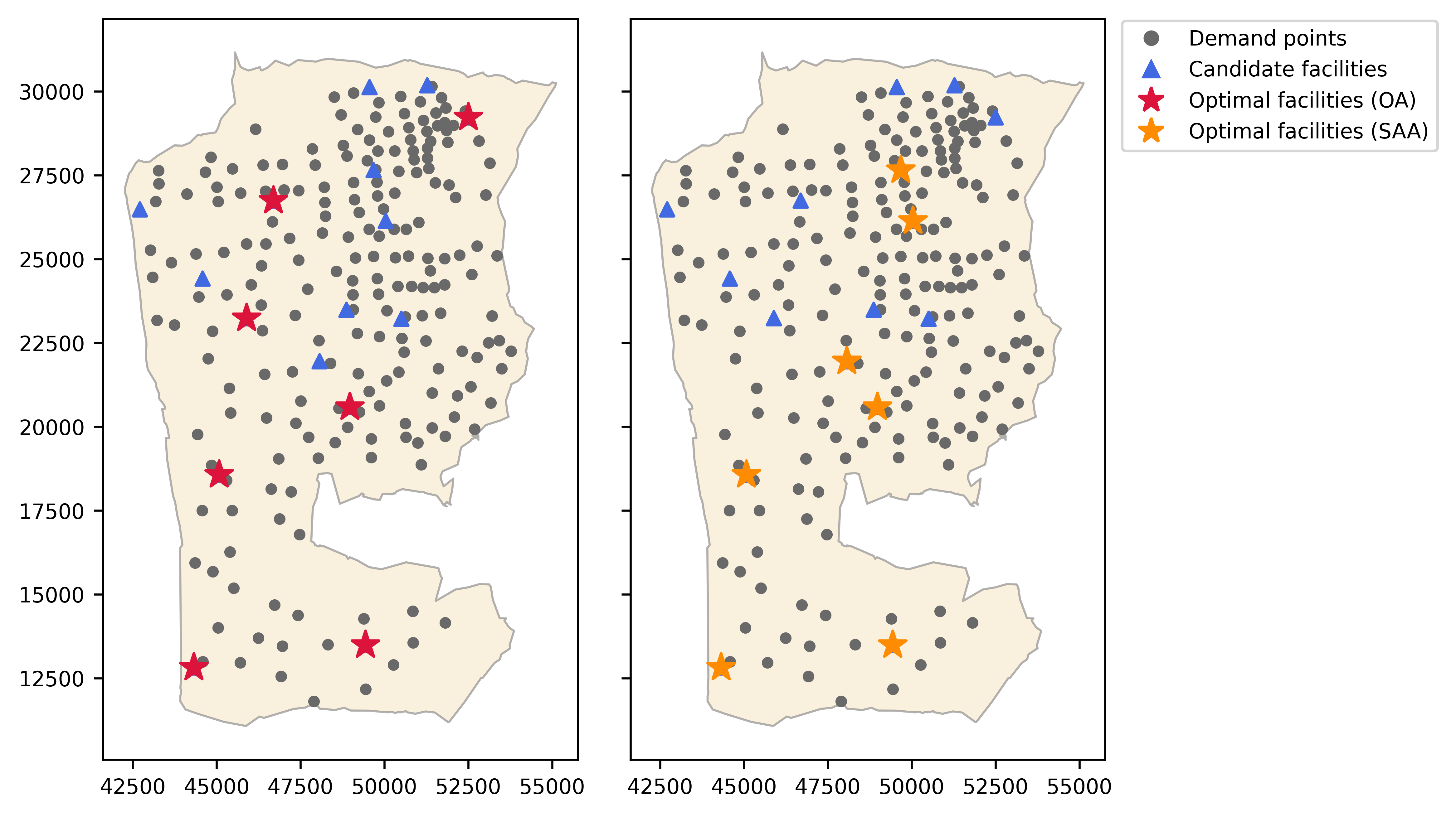} 
\caption{Optimal facility locations returned by OA and SAA methods for $\epsilon_i = 0.1$}
\label{fig: real-world}
\end{figure}

\begin{table}[htbp]
\centering
\caption{Computational results under different $\epsilon_i$ in the facility location problem ($N = 200$)}
\label{tab:real-world}
\resizebox{\textwidth}{!}{%
\begin{tabular}{ccccccccccc}
\toprule
& \multicolumn{3}{c}{OA-I} & \multicolumn{3}{c}{OA-II} & \multicolumn{2}{c}{SAA} & \multicolumn{2}{c}{IS} \\
\midrule
$\epsilon_i$ & Time & Iterations & Val & Time & Iterations & Val & Time & Val & Time & Val \\
\midrule
0.1 & 1361.90 & 3 & 7 & 2514.20 & 3 & 7 & 10.26 & 7 & 17.54 & 7 \\
0.2 & 1022.72 & 4 & 6 & 510.57  & 3 & 6 & 34.19 & 6 & 30.51 & 6 \\
0.3 & 372.12  & 6 & 5 & 228.97  & 4 & 5 & 50.39 & 5 & 38.51 & 5 \\
0.4 & 283.60  & 2 & 5 & 141.68  & 2 & 5 & 58.34 & 5 & 41.94 & 5 \\
0.5 & 217.69  & 2 & 4 & 288.11  & 3 & 4 & 19.00 & 4 & 13.10 & 4 \\
\bottomrule
\end{tabular}%
}
\end{table}

}

\subsection{Computational performance of methods for CC-SMCP} \label{sec: overall_performance}

{To further evaluate the performance of solution methods}, we randomly generate a \textit{sparse} probability matrix $P = (p_{ij})_{m \times n}$: we assume that the $i$-th row can only be covered by at most $12$ columns when $k_i \ge 2$. That is, for each row $i$ where $k_i \ge 2 $, we randomly choose $n'\le 12$ columns and sample the probability of success $p_{ij} = \mathbb{P}\left[\tilde{a}_{ij} = 1 \right]$ uniformly from the interval $[0.9, 1]$. For the other $n-n'$ columns in the row $i$, we set $p_{ij} = 0$. Each covering number $k_i$ is randomly selected from $\{1,2,3\}$, and costs $c_{ij}$ in the objective function are all set to 1. Further, we set the deterministic set $B = \{0,1\}^n$, and the risk parameter $\epsilon_i = 0.05$ and $0.1$ for each $i \in [m]$ in these experiments, reflecting the natural assumption that we want to cover the item $i$ with high probability. We implement the presolving process before we apply the OA method to solve each instance. For the SAA method, we let the risk parameter $\alpha_i := \epsilon_i$ in problem \eqref{eq: SAA-big-M-free} for each $i \in [m]$, and generate each scenario $\omega \in \Omega$ independently according to the probability matrix $P$. For the IS approach, we set $\hat{p}_{ij} := 0$ for any column $j$ where $p_{ij} = 0$. For columns $j \in [n]$ where $p_{ij} \neq 0$, we choose $u = \lceil (n' + k_i - 1) / 2 \rceil$ and compute the IS estimator $\hat{p}_{ij}$ using equations \eqref{eq: IS_est_calc}.

\begin{table}[h]
  \centering
  \caption{Comparison between OA and sampling-based methods ($p_{ij} \in [0.9,1], N = 200$)}
  \resizebox{\textwidth}{!}{%
      \begin{tabular}{ccccccccccccccccccccc}
    \toprule
    \multicolumn{3}{c}{Parameters} & \multicolumn{3}{c}{OA-I(w/o presolving)} & \multicolumn{3}{c}{OA-I} & \multicolumn{3}{c}{OA-II(w/o presolving)} & \multicolumn{3}{c}{OA-II} & \multicolumn{3}{c}{SAA} & \multicolumn{3}{c}{IS} \\
    $n$ & $m$ & $\epsilon$ & Time  & Iter. & Val   & Time  & Iter. & Val   & Time  & Iter. & Val   & Time  & Iter. & Val   & Time  & Val   & Gap(\%) & Time  & Val   & Gap(\%) \\
    \midrule
    30    & 10    & 0.05  & 3.918216 & 3     & 7     & 3.86  & 3     & 7     & 0.676478 & 2     & 7     & 0.52  & 2     & 7     & 0.15  & 7     & 0.00  & 0.26  & 7     & 0.00 \\
    30    & 10    & 0.1   & 1.8321 & 2     & 7     & 1.39  & 2     & 7     & 0.896609 & 2     & 7     & 0.58  & 2     & 7     & 0.42  & 7     & 0.00  & 0.21  & 7     & 0.00 \\
    30    & 20    & 0.05  & 3.045077 & 2     & 11    & 3.55  & 2     & 11    & 0.608429 & 2     & 11    & 1.24  & 2     & 11    & 0.44  & 11    & 0.00  & 0.63  & 11    & 0.00 \\
    30    & 20    & 0.1   & 13.3727 & 5     & 9     & 7.92  & 4     & 9     & 2.625688 & 3     & 9     & 1.16  & 2     & 9     & 1.33  & 9     & 0.00  & 0.51  & 8     & -11.11 \\
    30    & 30    & 0.05  & 7.165127 & 3     & 13    & 5.90  & 2     & 13    & 3.509374 & 3     & 13    & 4.05  & 2     & 13    & 8.37  & 13    & 0.00  & 2.82  & 13    & 0.00 \\
    30    & 30    & 0.1   & 2.060083 & 2     & 11    & 2.62  & 2     & 11    & 1.867163 & 2     & 11    & 1.13  & 2     & 11    & 12.12 & 11    & 0.00  & 1.16  & 11    & 0.00 \\
    30    & 50    & 0.05  & 10.83688 & 2     & 13    & 17.79 & 3     & 13    & 2.245444 & 2     & 13    & 2.28  & 2     & 13    & 2.30  & 13    & 0.00  & 2.22  & 13    & 0.00 \\
    30    & 50    & 0.1   & 14.17505 & 3     & 13    & 13.70 & 3     & 13    & 2.481468 & 2     & 13    & 2.24  & 2     & 13    & 3.67  & 13    & 0.00  & 4.19  & 13    & 0.00 \\
    30    & 100   & 0.05  & 11.99792 & 2     & 19    & 19.71 & 2     & 19    & 5.413821 & 3     & 19    & 3.29  & 3     & 19    & 4.26  & 19    & 0.00  & 3.86  & 19    & 0.00 \\
    30    & 100   & 0.1   & 23.77158 & 3     & 18    & 21.63 & 3     & 18    & 3.707056 & 2     & 18    & 5.52  & 2     & 18    & 6.36  & 18    & 0.00  & 3.42  & 17    & -5.56 \\
    30    & 150   & 0.05  & 16.21416 & 3     & 20    & 25.69 & 2     & 20    & 10.24752 & 3     & 20    & 5.80  & 3     & 20    & 6.20  & 20    & 0.00  & 7.22  & 20    & 0.00 \\
    30    & 150   & 0.1   & 17.53648 & 3     & 19    & 16.82 & 3     & 19    & 7.955401 & 2     & 19    & 9.42  & 3     & 19    & 8.09  & 18    & -5.26 & 9.65  & 19    & 0.00 \\
    50    & 30    & 0.05  & 26.29923 & 2     & 14    & 21.94 & 2     & 14    & 10.26118 & 2     & 14    & 7.44  & 2     & 14    & 3.21  & 14    & 0.00  & 1.59  & 13    & -7.14 \\
    50    & 30    & 0.1   & 19.14469 & 2     & 12    & 32.04 & 4     & 12    & 8.787309 & 3     & 12    & 6.63  & 3     & 12    & 10.37 & 13    & 8.33  & 2.54  & 11    & -8.33 \\
    50    & 50    & 0.05  & 27.06219 & 2     & 17    & 50.87 & 3     & 17    & 13.26956 & 2     & 17    & 14.39 & 3     & 17    & 2.78  & 17    & 0.00  & 1.97  & 16    & -5.88 \\
    50    & 50    & 0.1   & 127.7748 & 4     & 14    & 195.57 & 4     & 14    & 27.52517 & 3     & 14    & 50.22 & 4     & 14    & 6.42  & 14    & 0.00  & 2.00  & 12    & -14.29 \\
    50    & 100   & 0.05  & 44.16067 & 2     & 24    & 83.83 & 3     & 24    & 20.54323 & 2     & 24    & 16.98 & 2     & 24    & 8.13  & 24    & 0.00  & 6.19  & 24    & 0.00 \\
    50    & 100   & 0.1   & 116.9192 & 4     & 23    & 26.59 & 2     & 23    & 30.39361 & 3     & 23    & 22.74 & 3     & 23    & 12.66 & 22    & -4.35 & 6.92  & 22    & -4.35 \\
    50    & 150   & 0.05  & 64.4578 & 3     & 33    & 26.93 & 2     & 33    & 27.30152 & 4     & 33    & 20.35 & 3     & 33    & 9.76  & 33    & 0.00  & 10.99 & 33    & 0.00 \\
    50    & 150   & 0.1   & 17.55522 & 3     & 29    & 21.88 & 3     & 29    & 10.83044 & 3     & 29    & 13.02 & 3     & 29    & 11.31 & 29    & 0.00  & 12.08 & 29    & 0.00 \\
    100   & 50    & 0.05  & 2561.773 & 2     & 28    & -     & 3     & 28    & 551.1082 & 4     & 28    & 148.56 & 2     & 28    & 15.74 & 27    & -3.57 & 4.64  & 27    & -3.57 \\
    100   & 50    & 0.1   & 504.1497 & 2     & 23    & 2019.72 & 4     & 23    & 138.642 & 2     & 23    & 42.88 & 2     & 23    & 147.91 & 25    & 8.70  & 9.56  & 24    & 4.35 \\
    100   & 100   & 0.05  & 1866.637 & 3     & 40    & 592.19 & 2     & 40    & 185.2476 & 2     & 40    & 122.96 & 2     & 40    & 15.64 & 40    & 0.00  & 9.80  & 39    & -2.50 \\
    100   & 100   & 0.1   & 2010.724 & 4     & 35    & 1310.52 & 3     & 35    & 209.3176 & 5     & 35    & 108.59 & 3     & 35    & 30.20 & 35    & 0.00  & 49.53 & 36    & 2.86 \\
    100   & 150   & 0.05  & 961.1219 & 5     & 46    & 439.65 & 3     & 46    & 251.3717 & 5     & 46    & 183.10 & 4     & 46    & 0.19  & INF   & INF   & 18.07 & 46    & 0.00 \\
    100   & 150   & 0.1   & 1138.104 & 4     & 44    & 1660.78 & 5     & 44    & 277.3671 & 4     & 44    & 200.79 & 3     & 44    & 49.03 & 42    & -4.55 & 65.74 & 43    & -2.27 \\
    300   & 50    & 0.05  & 842.8555 & 3     & 43    & 198.26 & 3     & 43    & 653.3265 & 2     & 43    & 71.47 & 2     & 43    & 28.08 & 43    & 0.00  & 29.08 & 41    & -4.65 \\
    300   & 50    & 0.1   & 1071.726 & 4     & 38    & 248.80 & 4     & 38    & 956.1633 & 3     & 38    & 93.53 & 3     & 38    & 31.69 & 35    & -7.89 & 37.32 & 34    & -10.53 \\
    300   & 100   & 0.05  & -     & 2     & 39    & -     & 2     & 39    & 1111.11 & 2     & 61    & 1437.51 & 4     & 61    & 102.19 & 61    & 0.00  & 72.79 & 57    & -6.56 \\
    300   & 100   & 0.1   & -     & 2     & 36    & -     & 2     & 36    & 1937.254 & 4     & 54    & 796.72 & 3     & 54    & 89.33 & 52    & -3.70 & 92.51 & 52    & -3.70 \\
    300   & 150   & 0.05  & -     & 2     & 51    & -     & 2     & 51    & -     & 3     & 78    & -     & 3     & 78    & 192.23 & 78    & 0.00  & 135.22 & 76    & -2.56 \\
    300   & 150   & 0.1   & -     & 2     & 48    & -     & 2     & 48    & -     & 5     & 70    & 720.15 & 2     & 70    & 366.51 & 67    & -4.29 & 145.31 & 69    & -1.43 \\
    300   & 200   & 0.05  & -     & 2     & 56    & -     & 2     & 56    & 2283.844 & 3     & 88    & 1349.81 & 3     & 88    & 212.94 & 86    & -2.27 & 180.20 & 87    & -1.14 \\
    300   & 200   & 0.1   & -     & 2     & 49    & -     & 2     & 49    & 2770.346 & 4     & 77    & 1210.56 & 3     & 77    & 685.08 & 73    & -5.19 & 276.36 & 74    & -3.90 \\
    300   & 250   & 0.05  & 2298.834 & 2     & 114   & 2650.57 & 3     & 114   & 2078.619 & 3     & 114   & 787.82 & 3     & 114   & 283.98 & 113   & -0.88 & 262.12 & 112   & -1.75 \\
    300   & 250   & 0.1   & -     & 2     & 65    & -     & 2     & 65    & -     & 4     & 102   & 1249.90 & 3     & 102   & 843.70 & 101   & -0.98 & 474.56 & 100   & -1.96 \\
    300   & 300   & 0.05  & -     & 3     & 120   & 3434.26 & 3     & 121   & 2144.307 & 3     & 121   & 1358.17 & 4     & 121   & 453.30 & 120   & -0.83 & 336.46 & 120   & -0.83 \\
    300   & 300   & 0.1   & -     & 4     & 112   & -     & 4     & 112   & 2524.173 & 3     & 112   & 1989.73 & 6     & 112   & 717.31 & 108   & -3.57 & 561.76 & 110   & -1.79 \\
    \bottomrule
    \end{tabular}
    }
  \label{tab: comparison}%
\end{table}%

Table \ref{tab: comparison} lists some comparison results between four different methods: OA-I, OA-II, SAA and IS. \blue{We also report additional experiments comparing OA-I and OA-II with and without the presolving procedure.} {Similar as in Table \ref{tab:real-world}, we} use ``Val'' to denote the best achievable value obtained among four different methods within one hour for the corresponding instance. For OA methods, if an instance has been solved within one hour, then the ``Val'' of OA methods should be the optimal value of the instance. The symbol `-' in `Time' columns indicates that the method did not obtain an optimal solution within the time limit of 3,600 seconds. \blue{We use “Iter.” to denote the number of iterations performed by the OA methods.}  The gaps (Gap) reported are the percent by which the output value of SAA or IS is below/above the optimal value generated by the OA methods, defined as Gap = (Val $-$ opt)/opt, where opt denotes the optimal value of an instance and ``Val" represents the objective value obtained from the SAA or IS approach. Specifically, a negative (positive) gap indicates the output value of SAA or IS is below (above) the optimal value.
Since the output of SAA or IS depends on sampling results, an  ``INF" is placed in the column ``Gap (\%)'' if the corresponding SAA or IS problem is infeasible.

As can be seen in Table \ref{tab: comparison}, the OA approaches perform well on these instances with a sparse probability matrix, and the OA-II method outperforms the OA-I method in all $38$ cases. \blue{Since the OA algorithm solves a sequence of relaxed problems of CC-SMCP, the  objective value of the current relaxation provides a valid lower bound for the original minimization problem if the algorithm does not terminate within the time limit. If the current relaxation cannot be solved within the time limit, a valid lower bound can still be obtained from the objective value of the relaxation solved in the previous OA iteration. Moreover, if a candidate solution is verified to satisfy all original chance constraints during the OA procedure, then it provides an optimal solution. As shown in Table \ref{tab: comparison}, in our experiments, OA methods always return either an optimal value or a valid lower bound within the one-hour time limit across all 38 instances.} \blue{In addition, the OA approaches with presolving outperform their counterparts without presolving in most instances, particularly as the instance size increases. The improvement is mainly attributable to the reduction in the number of chance constraints through vector dominance and the use of equivalent deterministic reformulations when $k_i=1$. For a few small instances, however, the approaches without presolving perform slightly better. This may be because, in these cases, the vector-dominance procedure eliminates only a small number of chance constraints, while the equivalent reformulation in \eqref{eq: log-X} may not provide sufficient computational benefits for small instances.} Sampling-based methods (SAA and IS) provide good approximations for CC-SMCP on these instances, with output values close to the optimal ones. The solution time for SAA and IS may be less than that of OA methods when $n$ or $m$ increases. Note that most gaps in SAA and all gaps in IS are negative, indicating that sampling-based methods often struggle to find feasible solutions. 
Hence, developing methods and techniques for faster identification of feasible solutions in sampling-based approaches remains an area of interest.

\begin{figure}[htp]
\centering
\begin{tabular}{cc}
\includegraphics[width=0.38\textwidth]{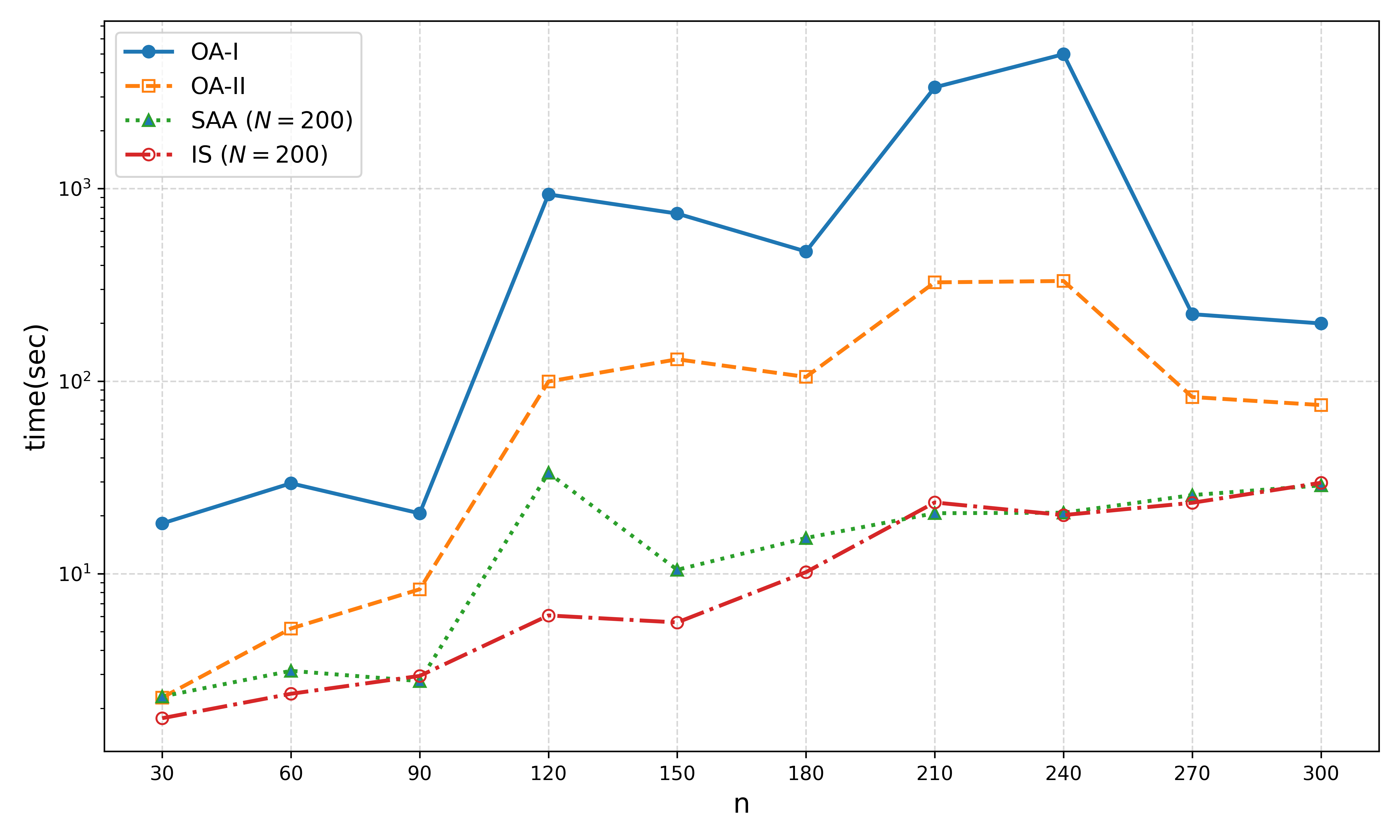}  &  \quad \includegraphics[width=0.38\textwidth]{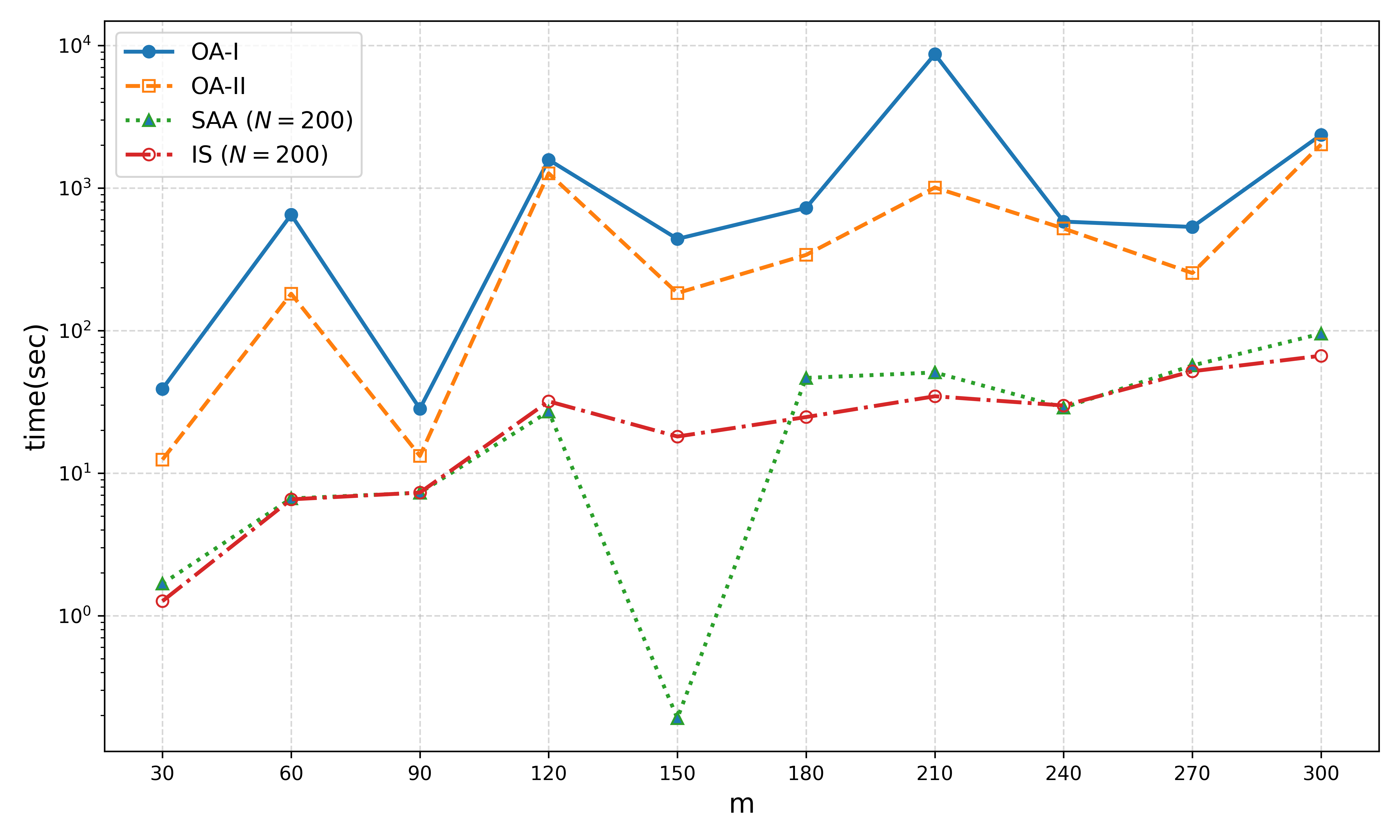}  \\
{\small (a) $m=50$ and $\epsilon = 0.05$} & {\small (b) $n=100$ and $\epsilon = 0.05$} \\
\end{tabular}
\caption{Runtime performance as a function of problem size}
\label{fig: time_vs_size}
\end{figure}

{
Figure~\ref{fig: time_vs_size} illustrates how solution times evolve as the problem size ($n$ or $m$) increases, with the vertical axis shown on a logarithmic scale. Overall, SAA and IS are much less sensitive to increases in $n$ or $m$ than the OA methods. One exception appears in Figure~\ref{fig: time_vs_size}(b) for the instance with \(n=100\) and \(m=150\), where SAA solves the problem very quickly because it returns an infeasible solution. In general, the runtimes of all four methods increase as the problem size grows; however, the OA methods are substantially more sensitive to this growth. This behavior suggests that their computational performance may also heavily depend on other factors, such as the sparsity structure of the probabilistic matrix and the effectiveness of presolving.

}

\begin{figure}[htp]
\centering
\begin{tabular}{cc}
\includegraphics[width=.35\textwidth]{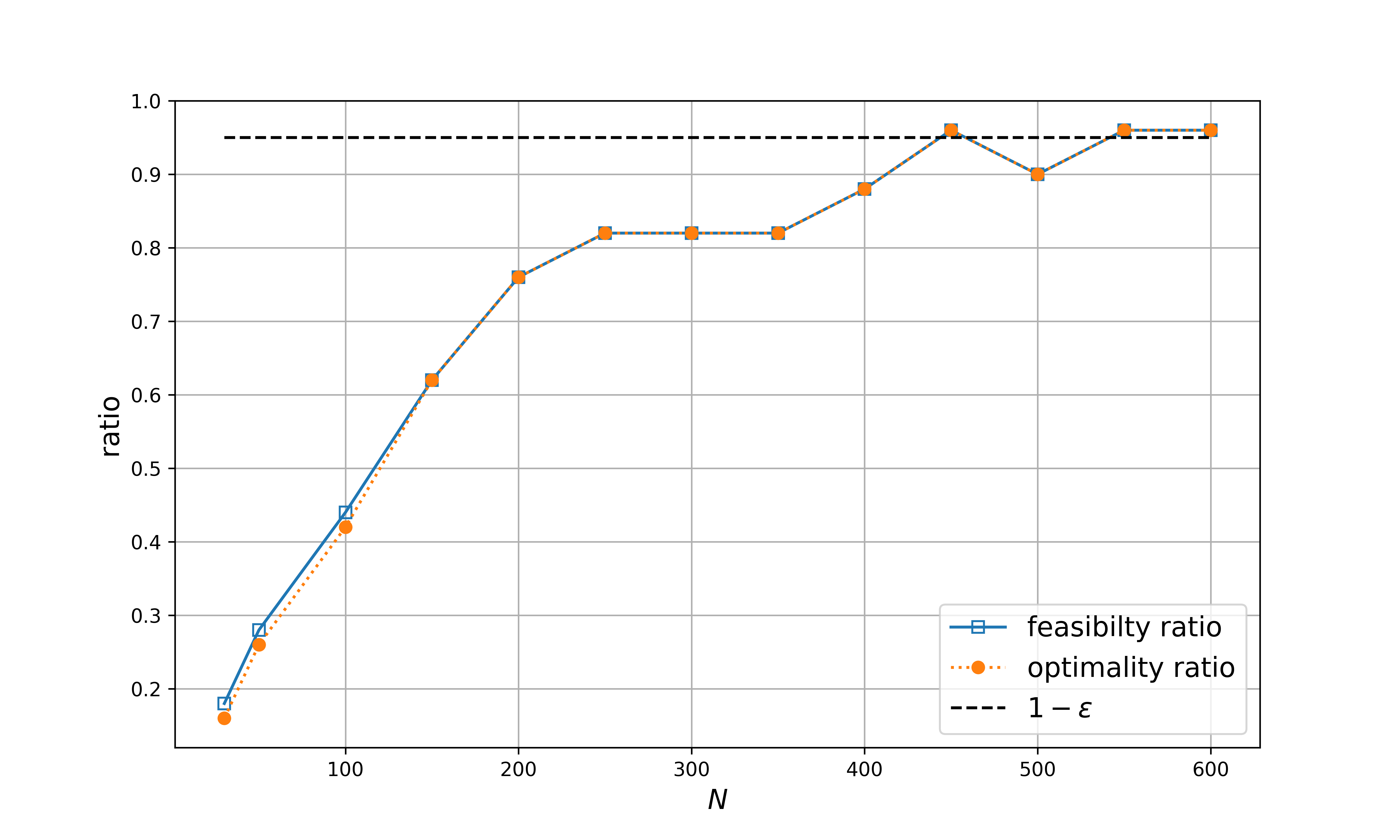} &  \includegraphics[width=.35\textwidth]{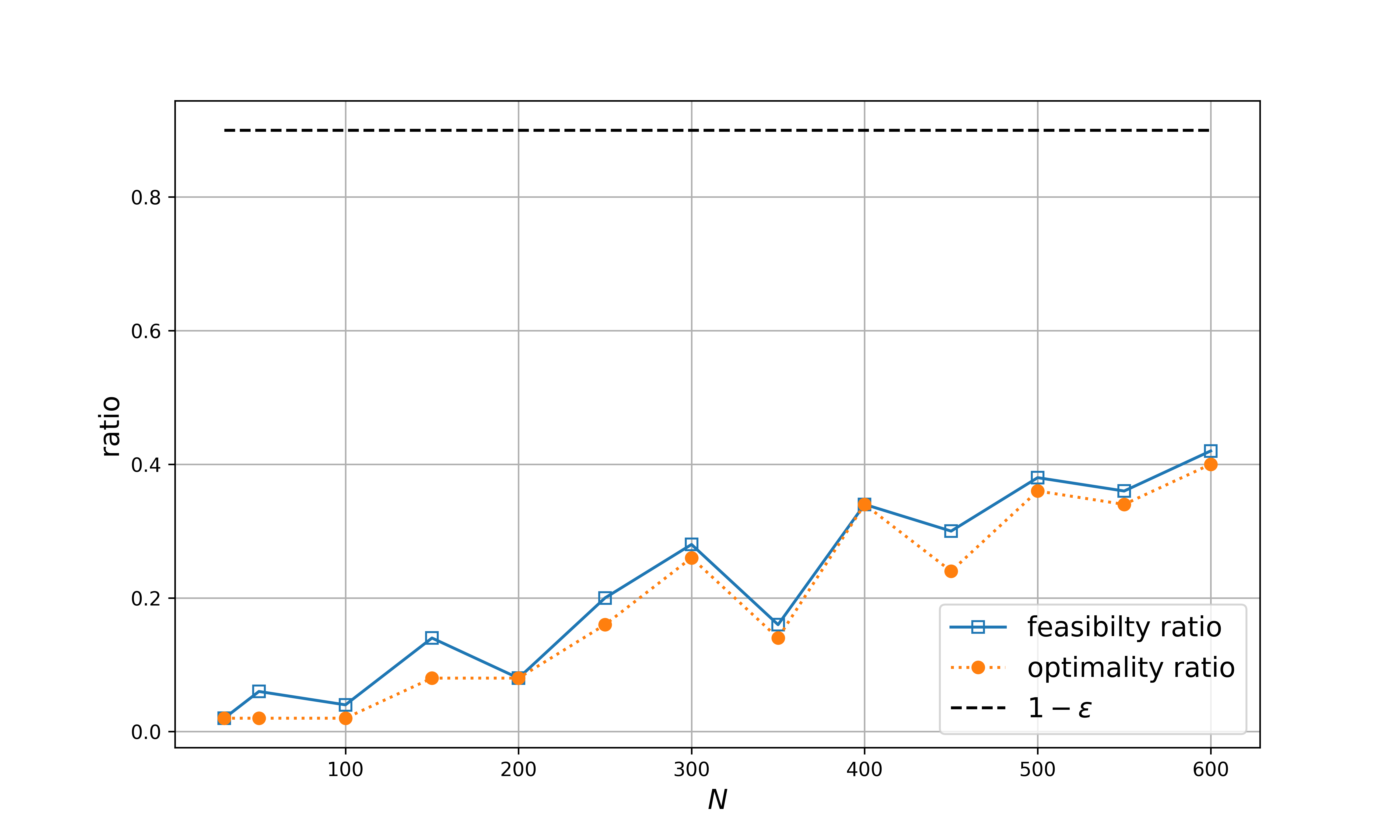} \\[-0.5em]
{\small (a) $n=50, m = 30$ and $\epsilon = 0.05$} & {\small (b) $n=50, m = 30$ and $\epsilon = 0.1$} \\
\includegraphics[width=.35\textwidth]{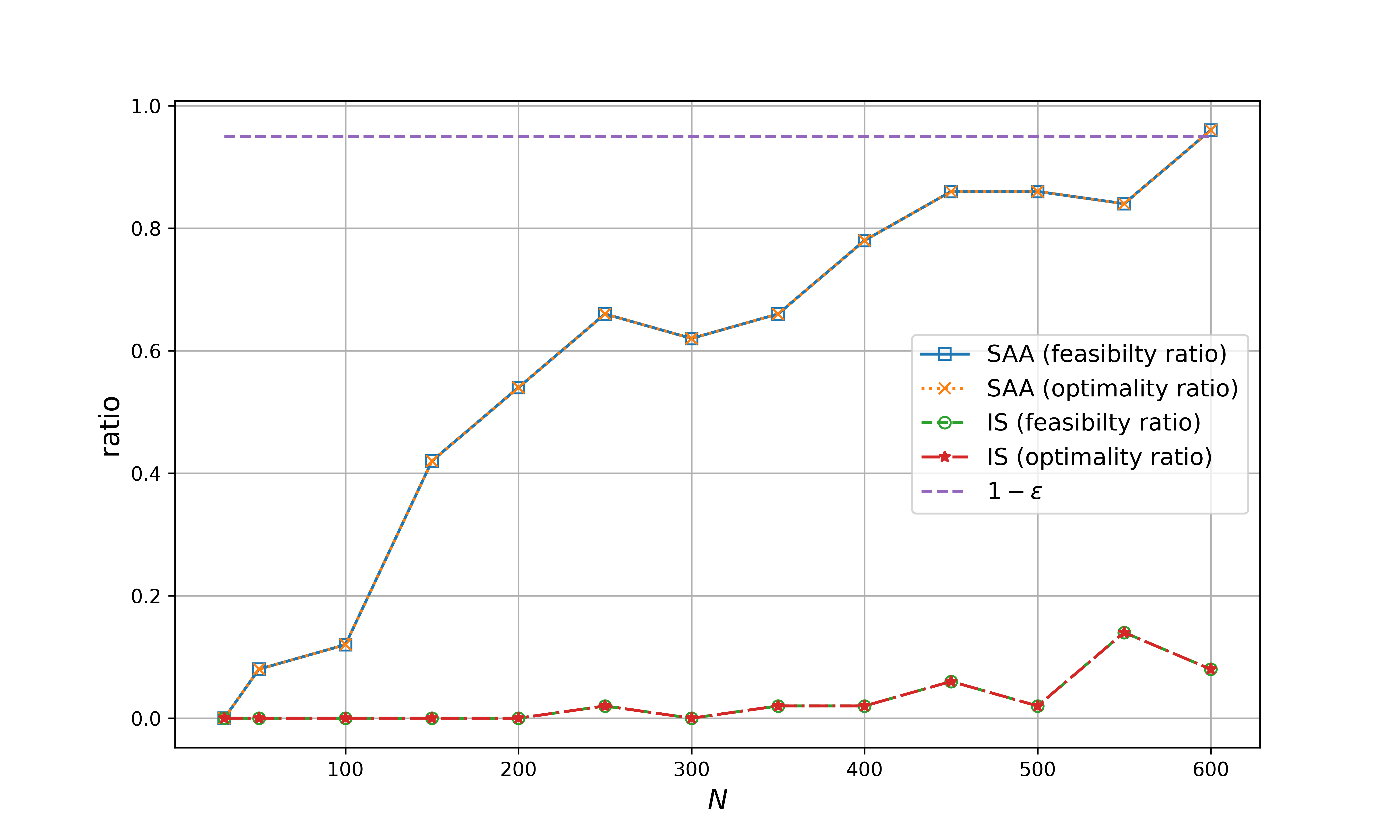} &  \includegraphics[width=.35\textwidth]{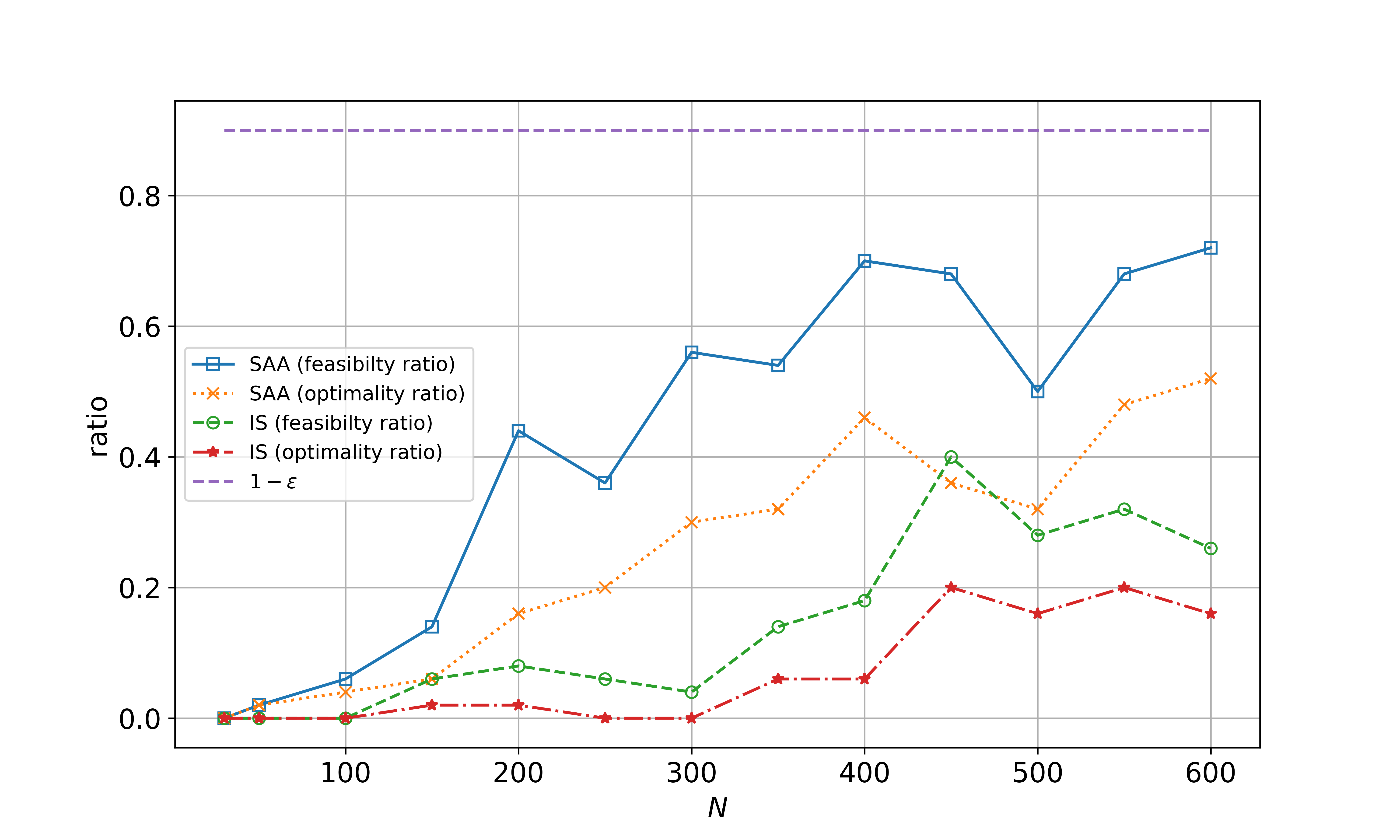} \\[-0.5em]
{\small (c) $n=100, m = 50$ and $\epsilon = 0.05$} & {\small (d) $n=100, m = 50$ and $\epsilon = 0.1$} \\
\end{tabular}
\caption{Feasibility-ratio and optimality-ratio curves of SAA and IS as functions of sample size $N$}
\label{fig: ratio}
\end{figure}

Note that even the presence of a zero gap in the SAA/IS results of Table \ref{tab: comparison} does not guarantee that the SAA/IS has generated a feasible solution for CC-SMCP. \blue{In this work, we do not report out-of-sample performance under independently generated testing scenarios because, in our setting, the true distribution is known and the exact probability expression in Lemma 1 can be used to evaluate any solution returned by SAA or IS. Therefore, empirical violation probabilities or coverage success rates based on finite test samples/scenarios would only approximate quantities that can already be computed exactly and would introduce additional sampling error. To assess the reliability and variability of the sampling-based methods, we instead conduct the following additional experiments with independently generated samples and report their feasibility ratios, optimality ratios, and average solution times.} We test the sampling-based approaches on two instances with $n = 50, m=30$ and $n = 100, m =50$ generated in the aforementioned
way. For the experiments, we consider $\epsilon_i = \epsilon  \in \{0.05, 0.1\}$ and take independent samples of size $N \in \{30, 50, 100, 150, 200, 250, 300, 350,$ $ 400, 450, 500, 550, 600\}$, and for each sample size we take $50$ different replications. 
We study how the risk parameter $\epsilon$ and the sample size $N$ affect the performance of SAA/IS. To better evaluate the performance of the sampling-based approaches, we define two performance measures: \textit{feasibility ratio} and \textit{optimality ratio}, and consider these two measures as functions of the sample size $N$. Specifically, for each $N$, the feasibility ratio is defined as the number of replications that output feasible solutions divided by the total number of replications, and the optimality ratio is defined as the number of replications that generate optimal solutions divided by the total number of replications. We can use Lemma \ref{lem: 1} to verify the feasibility of the solution provided by SAA/IS, and check the optimality by comparing the objective value generated by SAA/IS with the optimal value produced by the OA approaches. Figure \ref{fig: ratio} illustrates feasibility-ratio and optimality-ratio curves of SAA/IS as functions of $N$.

Notice that the feasibility ratio and optimality ratio of SAA increase as the sample size $N$ becomes larger. This makes sense as the larger the sample size $N$, the more likely SAA is to approximate the true distribution. Consequently, it becomes more likely for SAA to generate a feasible or optimal solution for CC-SMCP. This observation aligns with the theoretical findings previously presented in Theorem \ref{thm: feasible_region} and \ref{thm: SAA_opt}. However, practical experiments show that a smaller $N$ can also yield a feasible or optimal solution. For IS curves, the feasibility ratio and optimality ratio also increase with larger sample sizes $N$, but this effect is less pronounced compared to SAA curves. Moreover, it is important to note that the performance of IS is generally worse than that of SAA in almost all cases. This may be because IS is more sensitive to the choice of the importance sampling estimator $\hat{p}_{ij}$, which is crucial for the performance of IS.

\begin{figure}[htp]
\centering
\begin{tabular}{cc}
\includegraphics[width=0.38\textwidth]{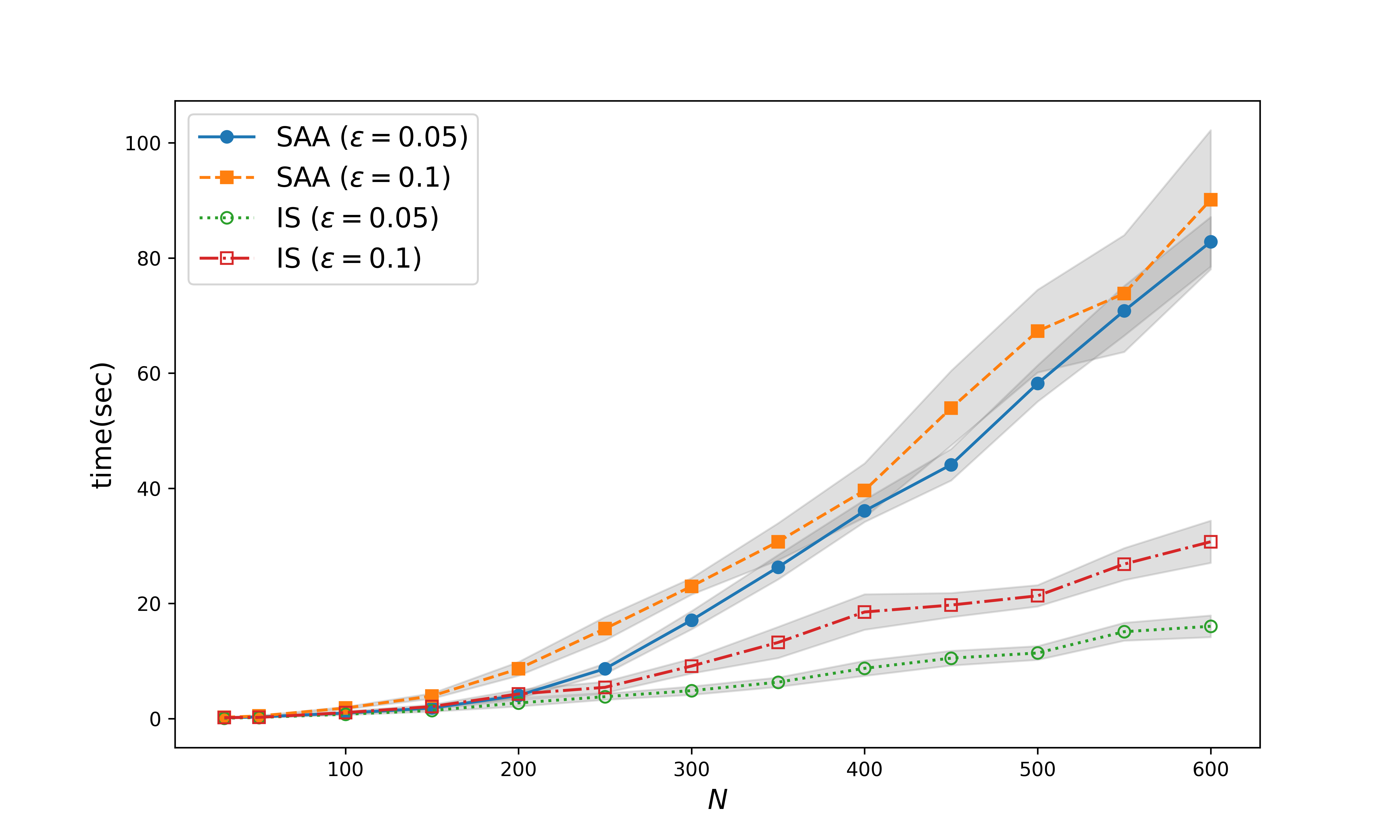}  &  \quad \includegraphics[width=0.38\textwidth]{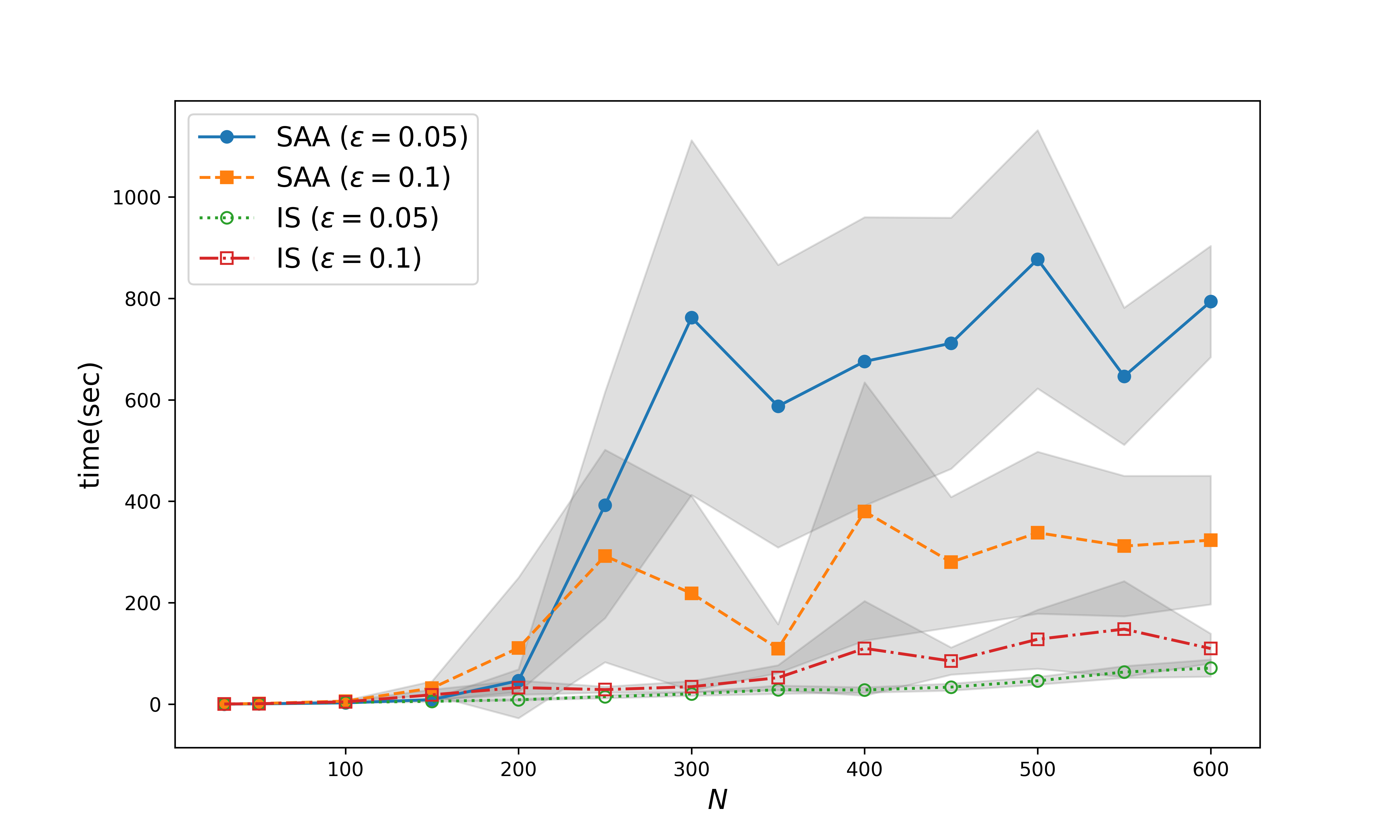}  \\
{\small (a) $n=50$ and $m = 30$} & {\small (b) $n=100$ and $m = 50$} \\
\end{tabular}
\caption{The average solution time of the SAA and IS methods as a function of $N$, where the 95\% confidence intervals are indicated as the shaded area}
\label{fig: avg_time}
\end{figure}

Furthermore, we compute the average solution time of SAA and IS as well as their $95$\% confidence intervals, as illustrated in Figure \ref{fig: avg_time}. As we observe from Figure \ref{fig: avg_time}, the solution time for
the IS at beginning iterations is similar to the solution time for SAA, but when the sample size $N$ becomes large, IS becomes significantly faster than SAA. This may be because the IS method is more efficient in generating feasible solutions satisfying constraints \eqref{eq: IS_likelihood_cons} rather than constraints \eqref{eq: SAA-big-M-2} in the SAA reformulation.

\subsection{Additional Experiments}

In addition to the above experiments, we also conduct several experiments to investigate the effect of risk parameter $\epsilon$ and the performance of methods for checking infeasibility and special cases. Please find more details in \ref{sec: additional_experiments} of this paper.


\section{Conclusion} \label{sec: conclusion}

In this paper, we have investigated CC-SMCP with individual chance constraints under LHS uncertainty, where both decision and random variables are purely binary. We derived exact deterministic reformulations and proposed an outer-approximation (OA) algorithm for CC-SMCP using some combinatorial methods. This OA method has significantly alleviated computational challenges and reduced cumulative rounding errors. Additionally, we established statistical relationships between CC-SMCP and its sample average approximation (SAA) reformulation and theoretically demonstrated the effectiveness of the SAA method in generating optimality bounds and feasible solutions for CC-SMCP. We also studied the importance sampling (IS) method and obtained a sufficient condition for selecting the optimal IS estimator. Finally, several computational experiments were conducted to validate the effectiveness of the OA methods (with two different linearization techniques) and sampling-based approaches (SAA and IS).

{
  In this work, the mutual independence assumption imposed on the underlying random variables enables the exact probability calculations and deterministic reformulations we develop. However, in many real applications,  disruptions and failures may be positively correlated. An important direction for future research is to incorporate structured dependence and to study whether our approach can be extended using tractable evaluations or provable bounds for joint event probabilities. A second promising direction is to extend this approach to distributionally robust variants (e.g., moment-based or Wasserstein-based chance-constrained formulations), leading to solutions that provide theoretical guarantees among possible probability distributions.

}

\bibliographystyle{elsarticle-harv} 
\bibliography{CC-SMCP-ref}

\clearpage
\appendix


\setcounter{page}{1}
\setcounter{lem}{0}
\setcounter{thm}{0}
\setcounter{equation}{0}
\counterwithin{thm}{section}
\counterwithin{lem}{section}
\renewcommand{\thelem}{\Alph{section}.\arabic{lem}}
\renewcommand{\thethm}{\Alph{section}.\arabic{thm}}

\begin{center}
{\bf \Large Supplementary Materials for\\
``Chance-Constrained Set Multicover Problem"}
\end{center}

\section{Preliminaries} 

\subsection{Inclusion–exclusion principle} \label{sec-apx: i-e}

\begin{lem}[Combinatorial version]
    For finite sets $A_1, \ldots, A_n$, one has the identity
    \begin{equation}
    \left|\bigcup_{i=1}^{n}A_{i}\right|=\sum_{\emptyset \neq S\subseteq [n]}(-1)^{|S|+1}\left|\bigcap_{j\in S}A_{j}\right|.
    \end{equation}
    
\end{lem}

\begin{lem}[Probabilistic version]
    For events $A_1, \ldots, A_n$ in a probability space $(\Omega,{\mathcal{F}},\mathbb{P})$, one has the identity
    \begin{equation}
    \mathbb {P} \left(\bigcup _{i=1}^{n}A_{i}\right)=\sum _{k=1}^{n}\left((-1)^{k-1}\sum _{I\subseteq [n] \atop |I|=k}\mathbb {P} \left(\bigcap _{i\in I}A_{i}\right)\right)
    \end{equation}
    
\end{lem}

\subsection{Bonferroni inequalities}

For events $A_1, \ldots, A_n$ in a probability space $(\Omega,{\mathcal{F}},\mathbb{P})$, we define 
\begin{equation*}
S_{1}:=\sum _{i=1}^{n}{\mathbb {P} }(A_{i}),\ S_{2}:=\sum _{1\leq i_{1}<i_{2}\leq n}{\mathbb {P} }(A_{i_{1}}\cap A_{i_{2}}),\ldots ,\ S_{k}:=\sum _{1\leq i_{1}<\cdots <i_{k}\leq n}{\mathbb {P} }(A_{i_{1}}\cap \cdots \cap A_{i_{k}})
\end{equation*}
for all integers $k$ in $\{1, ..., n\}$. Then we have the following inequalities:

\begin{lem}
    For odd $k \in \{1,\ldots,n\}$,
    \begin{equation}
      \sum _{j=1}^{k}(-1)^{j-1}S_{j}\geq \mathbb {P} \left(\bigcup _{i=1}^{n}A_{i}\right)=\sum _{j=1}^{n}(-1)^{j-1}S_{j}
      \end{equation}
\end{lem}

\begin{lem}
  For even $k \in \{2,\ldots,n\}$,
    \begin{equation}
      \sum _{j=1}^{k}(-1)^{j-1}S_{j}\leq \mathbb {P} \left(\bigcup _{i=1}^{n}A_{i}\right)=\sum _{j=1}^{n}(-1)^{j-1}S_{j}
    \end{equation}
\end{lem}

\subsection{Multiplicative Chernoff bound}
\begin{lem}[\citet{mitzenmacher2005probability}] \label{lem: chernoff}
Let $Y \sim Binomial(n; p)$ and $\mu = \mathbb{E}[Y]$. For any $\delta \ge 0$, 
\begin{equation*}
   {\mathbb{P}\left[Y\geq (1+\delta )\mu \right] \le \left({\frac {e^{\delta }}{(1+\delta )^{1+\delta }}}\right)^{\mu }} \le \exp \left(\frac{-\delta ^{2}\mu}{2+\delta} \right).
\end{equation*}
\end{lem}

\subsection{Hoeffding's inequality}
\begin{lem}[\citet{hoeffding1963probability}]
    Let $Y_1, \ldots, Y_N$ be independent random variables such that $a_i \le Y_i \le b_i$ almost surely. Then for all $t >0$,
    \begin{equation*}
    \mathbb{P}\left[\sum_{i=1}^N (Y - \mathbb{E}\left[Y_i\right]) \ge t\right] \le \exp \left\{- \frac{2t^2}{\sum_{i = 1}^N (b_i - a_i)^2} \right\}.
    \end{equation*}
\end{lem}

{
\subsection{Newton's inequalities} \label{ap:newton}
 \begin{lem}[\citet{hardy1952inequalities}]
\label{lem:newton}
Let $a_1,\ldots,a_n\ge 0$ be real numbers. For $k=0,1,\ldots,n$, let
\[
e_k
:=\sum_{\substack{T\subseteq [n]\\ |T|=k}}\ \prod_{i\in T} a_i
\]
denote the $k$-th elementary symmetric polynomial (with $e_0:=1$). Define the elementary symmetric means
\[
S_k:=\frac{e_k}{\binom{n}{k}},\qquad k=0,1,\ldots,n.
\]
Then, for every $k=1,2,\ldots,n-1$,
\[
S_{k-1}\,S_{k+1}\le S_k^{\,2}.
\]
Moreover, equality holds if and only if $a_1=\cdots=a_n$.
\end{lem} 

}

\section{Supplementary Materials for Sampling-Based Approaches}

\subsection{Sample average approximation} \label{sec-apx: SAA}

For ease of analysis, we rewrite the SAA reformulation as follows: 
\begin{equation}
 \nu_\alpha^N : = \min\left\{ c^T x :  x \in X_\alpha^N       \right\}, \label{eq: SAA-compact}
\end{equation}
where 
\begin{equation*}
X_\alpha^N :=\left\{x \in B: \textstyle\frac{1}{N} \ssum_{\omega \in \Omega} \mathbb{I}\left( \tilde{A}_{i}(\omega) x \ge k_i \right) \geq 1-\alpha_i, \ \forall i \in [m] \right\},
\end{equation*}
where $\mathbb{I}(\cdot)$ is the indicator function. 
Motivated by theoretical results in \citet{luedtke2008sample} regarding SAA for chance-constrained
problems with joint probabilistic constraints, we study approximations of chance-constrained combinatorial optimization problems with individual probabilistic constraints in this paper. Our goal is to establish statistical relationships between CC-SMCP and its SAA reformulation \eqref{eq: SAA-compact} for $\alpha_i \ge 0, \forall i \in [m]$. We assume that CC-SMCP has an optimal solution $x^*$ and a finite optimal value $\nu^*_\epsilon$. For each $i \in [m]$, we use $q_i(x^*)$ to denote the probability of the event $\{\tilde{A}_{i}(\omega) x^* \le k_i - 1\}$. 
By Lemma \ref{lem: 1}, we have $q_i(x^*) = \sum_{d = 0}^{k_i-1}\sum_{\ell^=d}^n(-1)^{\ell - d}\binom{\ell}{d} h_\ell(x^*)$ and $q_i(x^*) \le \epsilon_i$ since $x^* \in X_i$ is a solution to CC-SMCP. We now establish a bound on the probability that $\nu_\alpha^N$ yields a lower bound for $\nu^*_\epsilon$. 
\begin{thm} \label{thm: lower_bound}
Assume that $\nu^*_\epsilon$ and $\nu_\alpha^N$ are the optimal values to CC-SMCP and the SAA reformulation \eqref{eq: SAA-compact}, respectively. Then 
\begin{equation}
\mathbb{P}\left[\nu_\alpha^N \le \nu^*_\epsilon \right] \ge 1 - \ssum_{i = 1}^m I_{\epsilon_i}(\lfloor \alpha_i N \rfloor + 1, N - \lfloor \alpha_i N \rfloor),
\end{equation}
where $I_x(a,b)$ is the regularized incomplete beta function. Moreover,  if $\epsilon_i \le \alpha_i$ for each $i \in [m]$, then
\begin{equation}
\mathbb{P}\left[\nu_\alpha^N \le \nu^*_\epsilon \right] \ge 1-m\exp\left(-\kappa_1 N\right),
\end{equation}
where $\kappa_1 := \min_{i \in [m]} \left\{(\alpha_i -\epsilon_i)^2 /(\alpha_i + \epsilon_i)\right\} $.
\end{thm} 

Theorem \ref{thm: lower_bound} states that, when the risk parameter $\alpha_i > \epsilon_i$  for each $i \in [m]$, the SAA method yields a lower bound to the true optimal value of CC-SMCP with probability approaching one exponentially fast as $N$ increases. Further, given a confidence $1 - \delta$, where $0<\delta <1$, Theorem \ref{thm: lower_bound} ensures that $\nu_\alpha^N \le \nu^*_\epsilon$ with probability of at least $1 - \delta$, if we choose $\epsilon_i<\alpha_i$ for any $i \in [m]$ and the sample size 
\begin{equation*}
N \ge \textstyle\frac{1}{\kappa_1} \ln \frac{m}{\delta}.
\end{equation*}
In fact, with this choice of $\alpha_i$ and $N$, we have
\begin{equation*}
\mathbb{P}\left[\nu_\alpha^N \le \nu^*_\epsilon \right] \ge  1-m\exp\left(-\kappa_1 N\right) \ge 1 - \delta.
\end{equation*}

Next we investigate conditions under which an optimal solution of SAA problem 
\eqref{eq: SAA-compact} 
is a feasible solution to CC-SMCP. Let $X_\epsilon := \cap_{i \in [m]} X_{\epsilon, i}$  where $X_{\epsilon, i}:= \{x \in B: \mathbb{P}[\tilde{A}_{i} x \ge k_i] \ge 1 - \epsilon_i\}$. We assume that $\alpha_i < \epsilon_i$ for each $i \in [m]$. The idea is that if the risk level $\alpha_i$ is less than $\epsilon_i$ for each $i \in [m]$, then,  given a sufficiently large sample size $N$, the feasible region $X_\alpha^N$ will be a subset of $X_\epsilon$. Consequently, any optimal solution to the SAA reformulation \eqref{eq: SAA-compact} must be feasible for CC-SMCP. We have the following theorem: 
\begin{thm} \label{thm: feasible_region}
Assume that $\alpha_i < \epsilon_i$ for each $i = 1,\ldots, m,$ and let 
$\kappa_2 :=2 \min_{i \in [m]} (\epsilon_i - \alpha_i)^2$. Then
    \begin{equation*}
      \mathbb{P}\left[X_\alpha^N \subseteq X_\epsilon \right] \ge 1 - m|B \backslash X_\epsilon| \exp \left\{ -\kappa_2 N\right\}.
    \end{equation*}
\end{thm} 
    
The above theorem also provides us a way to estimate the sample size $N$ such that the feasible solutions of SAA reformulation \eqref{eq: SAA-compact} are feasible to CC-SMCP with a high probability (confidence) $1 -\delta$, if $\alpha_i < \epsilon_i$ for each $i \in [m]$ and we choose the sample size 
\begin{equation*}
N \ge \textstyle \kappa_2^{-1} 
\ln \frac{m |B \backslash X_\epsilon|}{\delta}.
\end{equation*}
Moreover, note that $B \subseteq \{0,1\}^n$ in CC-SMCP, we can take 
\begin{equation*}
N \ge\textstyle \frac{1}{\kappa_2} \ln \frac{m}{\delta} + \frac{n}{\kappa_2}\ln(2)
\end{equation*}
to attain the required confidence.

We are interested in determining the conditions under which the optimal value $\nu_\alpha^N$ of the SAA problem will converge to the true optimal value $\nu_\epsilon^*$ of CC-SMCP with probability one as $N$ approaches infinity for $\alpha = \epsilon$. Notice that Theorem \ref{thm: lower_bound} provides a trivial lower bound for the probability $\mathbb{P}\left[\nu_\epsilon^N \le \nu^*_\epsilon \right]$ when $\alpha_i = \epsilon_i$ for any $i \in [m]$. In the following, we will delve further by combining Theorem \ref{thm: lower_bound} and Theorem \ref{thm: feasible_region} to establish a lower bound for the probability $\mathbb{P}\left[\nu_\epsilon^N = \nu^*_\epsilon \right]$. Let $X_\epsilon^*$ be the set of optimal solutions to CC-SMCP, and define 
\begin{equation*}
\underline{\alpha}_i := \max \left\{ \mathbb{P}\left[\tilde{A}_i x \le k_i -1 \right] \ :\ x \in X_\epsilon^* \right\}, \ \forall i \in [m],
\end{equation*}
and $\underline{\alpha} := (\underline{\alpha}_1, \ldots, \underline{\alpha}_m)$. By definition, we have $\underline{\alpha}_i \le \epsilon_i$ for any $i \in [m]$ and $\nu_{\underline{\alpha}}^* = \nu^*_\epsilon$. 
Let  $X_\alpha^N = \cap_{i \in [m]} X_{\alpha, i}^N,$ where $X_{\alpha, i}^N := \left\{x \in B: \frac{1}{N} \sum_{\omega \in \Omega} \mathbb{I}\left( \tilde{A}_{i}(\omega) x \ge k_i \right) \geq 1-\alpha_i\right\}$,  
as well as 
\begin{equation*}
\overline{\alpha}_i := \min \left\{ \mathbb{P}\left[\tilde{A}_i x \le k_i -1 \right] \ :\ x \in B \backslash X_{\epsilon,i} \right\}, \ \forall i \in [m],
\end{equation*}
and $\overline{\alpha} := (\overline{\alpha}_1, \ldots, \overline{\alpha}_m)$. Without loss of generality, we assume that $B\backslash X_{\epsilon,i} \neq \varnothing$. Otherwise, since $B\backslash X_{\epsilon,i} = \varnothing$ implies $B \subseteq X_{\epsilon,i}$, we can remove the $i$th chance-constraint in CC-SMCP.  By definition, we have $\overline{\alpha}_i > \epsilon_i$ for any $i \in [m]$. Then we have the following theorem: 

\begin{thm} \label{thm: SAA_opt}
    Assume that $\underline{\alpha}_i < \epsilon_i$ for each $i \in [m]$, $\nu^*_\epsilon$ and $\nu_\alpha^N$ are the optimal values to CC-SMCP and the SAA reformulation \eqref{eq: SAA-compact}, respectively. Then
    \begin{equation*}
     \mathbb{P}\left[\nu_\epsilon^N = \nu^*_\epsilon \right] \ge 1 - m\left(|B \backslash X_\epsilon| + 1\right)\exp\{-\kappa_3 N\},
    \end{equation*}
where $\kappa_3 := \min \left\{\min_{i \in [m]} \left\{(\underline{\alpha}_i -\epsilon_i)^2 /(\underline{\alpha}_i + \epsilon_i)\right\},  2 \min_{i \in [m]} \{(\overline{\alpha}_i -\epsilon_i)^2\} \right\}.$
\end{thm} 

Note that the assumption that $\underline{\alpha}_i < \epsilon_i$ for each $i \in [m]$ is mild, because $B \subseteq \{0,1\}^n$ is finite and there are only a finite number of values of $\epsilon_i$ such that $\underline{\alpha}_i = \epsilon_i$. This fact inspires us to add a random perturbation uniformly distributed in $[-\gamma_i, \gamma_i]$ to $\epsilon_i$ when $\underline{\alpha}_i = \epsilon_i$, where $\gamma_i$ can be arbitrarily small, then the assumption will hold with probability one.

Theorem \ref{thm: SAA_opt} proves that solving a sample approximation with $\alpha = \epsilon$ will yield an exact optimal solution with probability approaching one exponentially fast with $N$. However, the sample size $N$ required to ensure a reasonably high probability of obtaining the optimal solution will be at least proportional to $\kappa_3^{-1}$ and therefore may be very large. Hence, Theorem \ref{thm: SAA_opt}, which reflects the qualitative behavior of the sample approximation with $\alpha = \epsilon$, may not be suitable for estimating the required sample size. In Section \ref{sec: experiments}, we will observe that SAA exhibits good performance even with a small sample size $N$.

\subsection{Importance Sampling} \label{sec-apx: IS}

We first briefly review the basic ideas of IS based on CC-SMCP as an example. For more details about IS methods and applications, we refer the interested reader to \citet{barrera2016chance} and \citet{rubino2009rare}. 

Let 
\begin{equation}
q_i(x) : = \mathbb{P}[\tilde{A}_ix < k] = \mathbb{E}_{\tilde{a}}[\mathbb{I}(\tilde{A}_ix < k_i)], ~~\forall i \in [m], \label{eq: qx}
\end{equation}
where $\tilde{A}_i := (\tilde{a}_{i1},\ldots,\tilde{a}_{in})$ is a row random vector. In the following we will assume $m=1$, and drop the row subscript $i$ like what we have done in Section \ref{sec: reformulation}. We want to estimate $q(x)$ for all $x \in X$. For a given sample $\tilde{a}(\omega_1), \ldots, \tilde{a}(\omega_N)$ of size $N$ from the distribution of $\tilde{a}$, a natural approximation of $q(x)$ in \eqref{eq: qx} is the SAA estimator:
\begin{equation}
  \hat{q}^{SAA}(x) : = \textstyle\frac{1}{N} \ssum_{\ell = 1}^N \mathbb{I}\left(\tilde{a}(\omega_\ell)x <k\right). \label{eq: SAA_estimator}
\end{equation}
Note that $\hat{q}^{SAA}(x)$ is an unbiased estimator of $q(x)$, since 
\[
\mathbb{E}_{\tilde{a}}\left[\hat{q}^{SAA}(x)\right] 
= \textstyle\frac{1}{N}\ssum_{\ell = 1}^N \mathbb{E}_{\tilde{a}}\left[\mathbb{I}\left(\tilde{a}(\omega_\ell)x <k\right)\right] = q(x).
\] 
Now let us consider $\hat{a}$ a new random vector and let $\hat{a}(\omega_1), \ldots, \hat{a}(\omega_N)$ be i.i.d. copies of $\hat{a}$. Define 
\begin{equation}
\hat{q}^{IS}(x) : = \textstyle\frac{1}{N} \ssum_{\ell = 1}^N \mathbb{I}\left(\hat{a}(\omega_\ell)x <k\right) L(\hat{a}(\omega_\ell)), \label{eq: IS_estimator}
\end{equation}
where $L(\cdot)$ is the likelihood ratio $L(\hat{a}) = \prod_{j = 1}^n (\frac{p_j}{\hat{p}_j})^{\hat{a}_j}(\frac{1- p_j}{1-\hat{p}_j})^{1- \hat{a}_j}$ with $p_j = P(\tilde{a}_j =1)$ and $\hat{p}_j = P(\hat{a}_j = 1)$ for each $j \in [n]$. In this case, $L$ is the ratio between the respective probability mass functions, since both $\tilde{a}$ and $\hat{a}$ have discrete support. Notice that for any function $f(\cdot): \mathbb{R}^n \mapsto \mathbb{R}$, we have 
\begin{equation}
\mathbb{E}_{\tilde{a}}[f(\tilde{a})] = \mathbb{E}_{\hat{a}}[f(\hat{a})L(\hat{a})]. \label{eq: IS_expectation}
\end{equation}
Based on the above observation, we obtain $\hat{q}^{IS}(x)$ is also an unbiased estimator of $q(x)$, since 
\begin{equation*}
\mathbb{E}_{\hat{a}}\left[\hat{q}^{IS}(x)\right] = \textstyle\frac{1}{N} \ssum_{\ell = 1}^N \mathbb{E}_{\hat{a}}\left[\mathbb{I}\left(\hat{a}(\omega_\ell)x <k\right) L(\hat{a}(\omega_\ell))\right]
\end{equation*}
and 
\begin{equation*}
\mathbb{E}_{\hat{a}}\left[\mathbb{I}\left(\hat{a}x <k\right) L(\hat{a})\right] = \mathbb{E}_{\tilde{a}}\left[\mathbb{I}(\tilde{a}x < k)\right] = q(x).
\end{equation*}
Further, the variance of the SAA estimator $\hat{q}^{SAA}(x)$ in \eqref{eq: SAA_estimator} is 
\begin{equation*}
\text{Var}\left[\hat{q}^{SAA}(x)\right] = \textstyle\frac{1}{N}\big(q(x) -q(x)^2\big)=\textstyle\frac{1}{N}\mathbb{E}_{\tilde{a}}\left[\mathbb{I}(\tilde{a}x < k)\right] - \textstyle\frac{1}{N}q(x)^2 
\end{equation*}
whereas the variance of the IS estimator $\hat{q}^{IS}(x)$ in \eqref{eq: IS_estimator} is given by
\begin{equation}
\text{Var}\left[\hat{q}^{IS}(x)\right] = \textstyle\frac{1}{N}
\mathbb{E}_{\hat{a}}\left[\mathbb{I}(\hat{a}x < k)^2L(\hat{a})^2\right] - \textstyle\frac{1}{N} q(x)^2 = \textstyle\frac{1}{N} \mathbb{E}_{\tilde{a}}\left[\mathbb{I}(\tilde{a}x < k)L(\tilde{a})\right] - \textstyle\frac{1}{N} q(x)^2, \label{eq: IS_var}
\end{equation}
where the second equality follows from \eqref{eq: IS_expectation}. Note that if we choose a ``good'' IS distribution in such a way that the event $\mathbb{I}{(\tilde{a}x<k)}$ becomes more likely under that distribution provided that $L(\tilde{a}) \le 1$, the variance of the IS estimator will be smaller than that of the standard SAA estimator.

To use the IS estimator $\hat{q}^{IS}(x)$ in the formulation of CC-SMCP, we only need to replace the chance-constraint \eqref{eq: set-k-covering-chance_constraints} by 
\begin{subequations}
\begin{align}
&\hat{A}_{i}(\omega) \mathbf{x} \ge k_iz_i(\omega), \forall \omega \in \Omega, \forall i \in [m] \\
&\ssum_{\omega \in \Omega} L(\hat{A}_i(\omega))\left(1-z_i(\omega)\right) \le N\epsilon_i, \forall i \in [m] \\
&z_i(\omega) \in \{0,1\}, \forall \omega \in \Omega
\end{align} 
\end{subequations}
where the definitions of $\hat{A}_i (\omega)$ and $z_i(\omega)$ are similar to those in the SAA reformulation \eqref{eq: SAA-big-M-free}. Note that if we choose $L(\hat{A}_i(\omega)) \equiv 1$ for each $i \in [m]$, then the obtained reformulation will exactly match the SAA reformulation as indicated in equation \eqref{eq: SAA-big-M-free}.

The remaining problem becomes finding a ``good'' IS estimator by minimizing the variance of $\hat{q}^{IS}(x)$. Suppose all components of $\tilde{a}$ and $\hat{a}$ are independent, from the expression \eqref{eq: IS_var}, we consider minimizing the term 
\begin{align}
\mathbb{E}_{\tilde{a}}\left[\mathbb{I}(\tilde{a}x < k)L(\tilde{a})\right] &= \mathbb{E}_{\tilde{a}}\left[\mathbb{I}\left(\ssum_{j = 1}^n\tilde{a}_jx_j < k\right)\sprod_{j = 1}^n \left(\frac{p_j}{\hat{p}_j}\right)^{\hat{a}_j}\left(\frac{1- p_j}{1-\hat{p}_j}\right)^{1- \hat{a}_j}\right] \notag\\
&= \mathbb{E}_{\tilde{a}}\left[\mathbb{I}\left(\ssum_{j = 1}^n\tilde{a}_jx_j < k\right)\sprod_{j = 1}^n \left(\frac{p_j(1-\hat{p}_j)}{\hat{p}_j(1-p_j)}\right)^{\hat{a}_j}\sprod_{j = 1}^n\left(\frac{1- p_j}{1-\hat{p}_j}\right)\right] \label{eq: IS_Item} 
\end{align}
For each $j \in [n]$, let 
\begin{equation*}\textstyle
\lambda_j := \log \left(\frac{p_j(1-\hat{p}_j)}{\hat{p}_j(1-p_j)}\right) = \log \left(\frac{1/\hat{p}_j -1}{1/p_j -1}\right)
\end{equation*}
Note that $\lambda_j \ge 0$ when $\hat{p}_j \le p_j$ which is  the case that IS distribution works by reducing the cover probability of each set $j$ such that the event $\sum_{j = 1}^n\tilde{a}_jx_j < k$ happens more often. By plugging $\lambda_j$ into equation \eqref{eq: IS_Item}, it follows that 
\begin{equation}
\mathbb{E}_{\tilde{a}}\left[\mathbb{I}(\tilde{a}x < k)L(\tilde{a})\right] = \mathbb{E}_{\tilde{a}}\left[\mathbb{I}\left(\ssum_{j = 1}^n\tilde{a}_jx_j < k\right) \exp\left(\ssum_{j = 1}^n \lambda_j \tilde{a}_j \right)\sprod_{j =1}^n \left(e^{-\lambda_j}p_j + (1 - p_j) \right)\right],
\end{equation}
where the equation follows from $\frac{1- p_j}{1- \hat{p}_j} = e^{-\lambda_j}p_j + (1- p_j)$.  The task of minimizing $\mathbb{E}_{\tilde{a}}\left[\mathbb{I}(\tilde{a}x < k)L(\tilde{a})\right]$ requires solving a multidimensional stochastic nonlinear problem, which is quite challenging. Alternatively, our approach focuses on minimizing the largest term within the expectation, that is
\begin{equation}
B_x(\lambda) := \max_{ \substack{\tilde{a}: \sum_{j = 1}^n \tilde{a}_jx_j < k \\ \tilde{a} \in \{0,1\}^n}} \exp\left(\ssum_{j = 1}^n \lambda_j \tilde{a}_j \right)\sprod_{j =1}^n \left(e^{-\lambda_j}p_j + (1 - p_j) \right).
\end{equation}
Note that $B_x(\lambda) \ge 0$. We consider minimizing $B_x(\lambda)$ over $\lambda \ge 0$. Then we have the following theorem that reduces the minimization problem to a one-dimensional problem, which can be efficiently solved by using some numerical methods. 

\begin{thm}\label{thm: IS_OPT}
    Suppose that $0<p_j <1$ for all $j = 1,\ldots, n$. Let $x$ be the solution that satisfies $\sum_{j = 1}^n x_j =u$ and $n-u+k-1 < \sum_{j = 1}^n p_j x_j$. Then the function $B_x(\lambda)$ is convex and there exists $\lambda_x^* \in \mathbb{R}_+$ such that the vector $\lambda$ defined as $\lambda_j = \lambda_x^*$ for each $j = 1, \ldots, n$ minimizes $B_x(\lambda)$. Moreover, $\lambda_x^*$ and $\hat{p}_j(\lambda_x^*)$ satisfy
    \begin{equation}
    \ssum_{j = 1}^n \hat{p}_j(\lambda_x^*) = n-u+k-1 \text{ and } \hat{p}_j(\lambda_x^*) = \frac{e^{-\lambda_x^*}p_j}
    {e^{-\lambda_x^*}p_j+(1- p_j)}.
    \label{eq: IS_est_calc}
    \end{equation}
\end{thm}

\begin{rem}
    Note that if $\sum_{j=1}^n x_j <k$, then all components of vector $\tilde{a}$ in $B_x(\lambda)$ should be one, which implies the optimal solutions to minimizing $B_x(\lambda)$ should be $\lambda_1 = \lambda_2 = \cdots = \lambda_n = 0$, since $\nabla \log B_x(\lambda) \ge 0$ for all $\lambda \ge 0$. In this case, due to the definition of $\lambda_j$, we have $\hat{p}_j = p_j$ for each $j \in [n]$.
\end{rem} 

Theorem \ref{thm: IS_OPT} provides a method for selecting good IS parameters for a given solution $x \in X$. To choose the importance sampling estimator for the whole set $X$, based on the conditions in Theorem \ref{thm: IS_OPT}, we notice that a reasonable lower bound for the sum $\sum_{j =1} ^n x_{ij}$ can be $\lceil (n+k_i-1)/2 \rceil$. This is because $n-u + k_i-1 <\sum_{j =1} ^n p_{ij}x_{ij} \le \sum_{j =1} ^n x_{ij} = u$. In the following numerical experiments, we will use $\lceil (n+k_i-1)/2 \rceil$ as an estimator for the sum $\sum_{j =1} ^n x_{ij}$, and then calculate the value for $\hat{p}_{ij}$ using equations $\eqref{eq: IS_est_calc}$. However, as we do not know the exact value of the sum $\sum_{j}x_{ij}$, this choice can sometimes be too radical to obtain a feasible solution to the CC-SMCP, as observed in the experimental section.

\section{Proofs}

\subsection{Proof of Lemma \ref{lem: 1}}

Note that $\mathbb{P}[A_{j}(x)] =\mathbb{P}[  \tilde{a}_j x_j  =1]= x_jp_{j}$ because $x_j \in \{0,1\}$. We first derive a closed-form expression for the probability $\mathbb{P}\left[\sum_{j=1}^n \tilde{a}_j x_j  = k \right]$, which can be regarded as the probability of exactly $k$ out of $n$ events $A_{1}(x), \ldots, A_{n}(x)$ occurring. We have
\begin{align*}
\mathbb{P}[\sum_{j=1}^n \tilde{a}_j x_j  = k ] & = \sum_{\substack{S \subseteq [n]\\ |S| = k}}\mathbb{P}[A_S(x)\bigcap_{j \in S^c} A_{j}^c(x) ] 
 = \sum_{\substack{S \subseteq [n]\\ |S| = k}}\mathbb{P}\left[A_S(x)\right]  (1 - \mathbb{P}[\bigcup_{j \in S^c} A_j(x)] )\\
& = \sum_{\substack{S \subseteq [n]\\ |S| = k}} \left(\mathbb{P}\left[A_S(x)\right] - \mathbb{P}[\bigcup_{j \in S^c} A_{S \cup \{j\}}(x)] \right),
\end{align*}
where the last two equalities follow from independence. Moreover, by the inclusion-exclusion principle, the probability of finite unions of events
\begin{small}
\begin{align*}
\mathbb{P}[\bigcup_{j \in S^c} A_{S \cup \{j\}}(x)] & = \sum_{j_1 \in S^c}\mathbb{P}\left[A_{S \cup \{j_1\}}(x)\right] - \sum_{\substack{j_1<j_2 \\ j_1,j_2\in S^c}} \mathbb{P}\left[A_{S \cup \{j_1, j_2\}}(x)\right] + \cdots + (-1)^{n- |S|-1} \mathbb{P}\left[A_{[n]}(x)\right]  \\
& = \sum_{S \subsetneq T \subseteq [n]} (-1)^{|T| - |S|-1}\mathbb{P}\left[A_{T}(x)\right].
\end{align*}
\end{small}

Therefore,
\begin{small}
\begin{align}
\mathbb{P} [\sum_{j=1}^n \tilde{a}_j x_j  = k ] & = \sum_{\substack{S \subseteq [n]\\ |S| = k}}\left(\mathbb{P} [A_S(x) ] - \sum_{S \subsetneq T \subseteq [n]} (-1)^{|T| - |S|-1}\mathbb{P}\left[A_{T}(x)\right] \right) \notag \\
& = \sum_{\substack{S \subseteq [n]\\ |S| = k}} \sum_{S \subseteq T \subseteq [n]} (-1)^{|T| - |S|} \mathbb{P}\left[A_{T}(x)\right] \notag\\
& = \sum_{\ell = k}^n \sum_{\substack{S \subseteq [n] \\ |S| = k} }\sum_{\substack{S \subseteq T \subseteq [n] \\ |T|= \ell} } (-1)^{|T| - |S|} \mathbb{P}\left[A_{T}(x)\right]  \\
 & = \sum_{\ell = k}^n (-1)^{\ell - k} \binom{\ell}{k}  h_{\ell}(x), \label{eq: distribution}
\end{align}
\end{small}
where $h_\ell(x) =\sum_{\substack{T \subseteq [n]\\ |T| = \ell}} \mathbb{P}\left[A_T(x)\right] = \sum_{\substack{T \subseteq [n]\\ |T| = \ell}}\prod_{j \in T} x_j p_j,$ for all $\ell = k, \ldots, n.$ And the last equality holds because for each $\ell$, $h_\ell(x)$ is counted $\binom{\ell}{k}$ times. Then we derive the closed-form expression for the probability $\mathbb{P}\left[\sum_{j=1}^n \tilde{a}_j x_j  \ge k \right]$, that is, the probability of at least $k$ out of $n$ events $A_1(x), \ldots, A_n(x)$ occurring:
\begin{align*}
\mathbb{P} [\sum_{j=1}^n \tilde{a}_j x_j  \ge k ] & = \sum_{d=k}^n\mathbb{P} [\sum_{j=1}^n \tilde{a}_j x_j  = d ] 
 = \sum_{d= k}^n\sum_{\ell = d}^n (-1)^{\ell - d} \binom{\ell}{d}  h_{\ell}(x) \\ 
& = {\sum_{\ell= k}^n\sum_{d = k}^\ell (-1)^{\ell - d} \binom{\ell}{d}  h_{\ell}(x)
 = \sum_{\ell= k}^n (-1)^{\ell -k} \sum_{d = k}^\ell (-1)^{k - d} \binom{\ell}{d}  h_{\ell}(x)} \\
&  = \sum_{\ell= k}^n (-1)^{\ell -k} \sum_{d = k}^\ell (-1)^{d - k} \binom{\ell}{d}  h_{\ell}(x) = \sum_{\ell= k}^n (-1)^{\ell -k}  \binom{\ell -1}{\ell -k}  h_{\ell}(x),
\end{align*}
where the {last two identities follows from $(-1)^{k-d} = (-1)^{d-k}$ and the combinatorial identity $(-1)^{k-1}\binom{\ell -1}{k-1}=\sum_{d = 0}^{k-1} (-1)^d \binom{\ell}{d} = -\sum_{d = k}^{\ell} (-1)^d \binom{\ell}{d}$.}

\subsection{Proof of Theorem \ref{thm: bounds} }
 From the proof of Lemma \ref{lem: 1}, for any fixed $S \subseteq [n]$ with $|S| = d$, we have 
    \begin{equation*}
    \mathbb{P}[A_S(x)\bigcap_{j \in S^c} A_{j}^c(x)] = \mathbb{P}\left[A_S(x)\right] - \mathbb{P}[\bigcup_{j \in S^c} A_{S \cup \{j\}}(x)]=\sum_{S \subseteq T \subseteq [n]} (-1)^{|T| - |S|} \mathbb{P}\left[A_{T}(x)\right].
    \end{equation*}
Applying Bonferroni inequalities to the probability $\mathbb{P}\left[\bigcup_{j \in S^c} A_{S \cup \{j\}}(x)\right]$ yields 
{
\begin{align*}
\mathbb{P}\!\left[\bigcup_{j\in S^c} A_{S\cup\{j\}}(x)\right]
&=\sum_{S \subsetneq T \subseteq [n]} (-1)^{|T|-|S|-1}\,\mathbb{P}\!\left[A_T(x)\right] \\
&\ge \sum_{\ell=d}^{\min\{d+t,n\} } \;
\sum_{\substack{S \subsetneq T \subseteq [n]\\ |S|=d,\; |T|=\ell}}
(-1)^{|T|-|S|-1}\,\mathbb{P}\!\left[A_T(x)\right],
\text{for even $t$.}
\end{align*}
Therefore,
}
\begin{align*}
\mathbb{P}[A_S(x)\bigcap_{j \in S^c} A_{j}^c(x)] & =  \sum_{S \subseteq T \subseteq [n] } (-1)^{|T| - |S|} \mathbb{P}\left[A_{T}(x)\right]  = \sum_{\ell=d}^{n} \sum_{\substack{S \subseteq T \subseteq [n] \\|S| =d, |T|= \ell} } (-1)^{|T| - |S|} \mathbb{P}\left[A_{T}(x)\right] \\
& {\le \sum_{\ell=d}^{\min\{d+t,n\}} \sum_{\substack{S \subseteq T \subseteq [n] \\ |S| =d, |T|= \ell} } (-1)^{|T| - |S|} \mathbb{P}\left[A_{T}(x)\right], \text{for even $t$.}} 
\end{align*}
{Similarly, we have
\begin{align*}
\mathbb{P}[A_S(x)\bigcap_{j \in S^c} A_{j}^c(x)]  & =  \sum_{S \subseteq T \subseteq [n] } (-1)^{|T| - |S|} \mathbb{P}\left[A_{T}(x)\right]  = \sum_{\ell=d}^{n} \sum_{\substack{S \subseteq T \subseteq [n] \\|S| =d, |T|= \ell} } (-1)^{|T| - |S|} \mathbb{P}\left[A_{T}(x)\right]  \\
& \ge \sum_{\ell=d}^{\min\{d+t,n\}} \sum_{\substack{S \subseteq T \subseteq [n] \\ |S| =d, |T|= \ell} } (-1)^{|T| - |S|} \mathbb{P}\left[A_{T}(x)\right], \text{for odd $t$.}
\end{align*} }
Moreover, if the inequality
\begin{equation*}
\sum_{\substack{S \subseteq T \subseteq [n]\\|S| = d, |T| = \ell}} \mathbb{P}\left[A_T(x)\right] \ge \sum_{\substack{S \subseteq T \subseteq [n]\\  |S| = d, |T| = \ell + 1}} \mathbb{P}\left[A_T(x)\right],
\end{equation*}
holds for any fixed $S \subseteq [n]$ with $|S| = d \le \ell<n$, then we have
\begin{small}
\begin{subequations}
\begin{align}
\mathbb{P}[A_S(x)\bigcap_{j \in S^c} A_{j}^c(x)] & = \sum_{\ell=d}^{n} \sum_{\substack{S \subseteq T \subseteq [n] \\ |S| = d, |T|= \ell} } (-1)^{|T| - |S|} \mathbb{P}\left[A_{T}(x)\right] \le \cdots \le \sum_{\ell=d}^{d+2} \sum_{\substack{S \subseteq T \subseteq [n] \\ |S| = d, |T|= \ell} } (-1)^{|T| - |S|} \mathbb{P}\left[A_{T}(x)\right] \\
& \le \sum_{\ell=d}^{d} \sum_{\substack{S \subseteq T \subseteq [n] \\ |S| = d, |T|= \ell} } (-1)^{|T| - |S|} \mathbb{P}\left[A_{T}(x)\right] 
\end{align} \label{eq: upper_bounds}
\end{subequations}
\end{small}
and 
\begin{small}
\begin{subequations}
\begin{align*}
\mathbb{P}[A_S(x)\bigcap_{j \in S^c} A_{j}^c(x)] & = \sum_{\ell=d}^{n} \sum_{\substack{S \subseteq T \subseteq [n] \\|S| = d,  |T|= \ell} } (-1)^{|T| - |S|} \mathbb{P}\left[A_{T}(x)\right] \ge \cdots \ge \sum_{\ell=d}^{d+3} \sum_{\substack{S \subseteq T \subseteq [n] \\|S| = d,  |T|= \ell} } (-1)^{|T| - |S|} \mathbb{P}\left[A_{T}(x)\right] \\
& \ge \sum_{\ell=d}^{d+1} \sum_{\substack{S \subseteq T \subseteq [n] \\ |S| = d, |T|= \ell} } (-1)^{|T| - |S|} \mathbb{P}\left[A_{T}(x)\right] 
\end{align*}
\end{subequations}
\end{small}
Taking the summations on both sides of inequalities \eqref{eq: upper_bounds}, we obtain
\begin{footnotesize}
\begin{subequations}
\begin{align*}
\mathbb{P}[\sum_{j=1}^n \tilde{a}_j x_j  \le k -1 ] &= \sum_{d = 0}^{k -1} \mathbb{P}[\sum_{j=1}^n \tilde{a}_j x_j = d ] = \sum_{d = 0}^{k -1}  \sum_{\substack{S \subseteq [n]\\ |S| = d}}\mathbb{P}[A_S(x)\bigcap_{j \in S^c} A_{j}^c(x)] \\
 & =\sum_{d = 0}^{k -1} \sum_{\substack{S \subseteq [n]\\ |S| = d}}  \sum_{\ell=d}^{n} \sum_{\substack{S \subseteq T \subseteq [n] \\ |T|= \ell} } (-1)^{|T| - |S|} \mathbb{P}\left[A_{T}(x)\right] \le \cdots \le \sum_{d = 0}^{k -1} \sum_{\substack{S \subseteq [n]\\ |S| = d}} \sum_{\ell=d}^{d+2} \sum_{\substack{S \subseteq T \subseteq [n] \\ |T|= \ell} } (-1)^{|T| - |S|} \mathbb{P}\left[A_{T}(x)\right] \\
& \le \sum_{d = 0}^{k -1} \sum_{\substack{S \subseteq [n]\\ |S| = d}} \sum_{\ell=d}^{d} \sum_{\substack{S \subseteq T \subseteq [n] \\ |T|= \ell} } (-1)^{|T| - |S|} \mathbb{P}\left[A_{T}(x)\right] 
\end{align*}
\end{subequations}
\end{footnotesize}
{From eq. \eqref{eq: distribution} of Lemma \ref{lem: 1}, we have  
\begin{align*}
 \sum_{\substack{S \subseteq [n] \\ |S| = k} }\sum_{\substack{S \subseteq T \subseteq [n] \\ |T|= \ell} } (-1)^{|T| - |S|} \mathbb{P}\left[A_{T}(x)\right]  = (-1)^{\ell - k} \binom{\ell}{k}  h_{\ell}(x). 
\end{align*}
Therefore, we conclude that
}
\begin{subequations}
\begin{align*}
\mathbb{P}\left[\sum_{j=1}^n \tilde{a}_j x_j  \le k -1 \right] & =\sum_{d = 0}^{k-1}\sum_{\ell^=d}^{n}(-1)^{\ell - d}\binom{\ell}{d}  h_\ell(x) \le \cdots \le \sum_{d = 0}^{k-1}\sum_{\ell^=d}^{d+2}(-1)^{\ell - d}\binom{\ell}{d}  h_\ell(x) \\
& \le \sum_{d = 0}^{k-1}\sum_{\ell^=d}^{d}(-1)^{\ell - d}\binom{\ell}{d}  h_\ell(x), 
\end{align*}
\end{subequations}
which is equivalent to 
\begin{equation*}
\mathbb{P}\left[\sum_{j=1}^n \tilde{a}_j x_j  \le k -1 \right]=g_n(x) \le \cdots \le g_4(x) \le g_2(x) \le g_0(x), 
\end{equation*}
and the sequence of inequalities
\begin{equation*}
\mathbb{P}\left[\sum_{j=1}^n \tilde{a}_j x_j  \le k -1 \right]=g_n(x) \ge \cdots \ge g_5(x) \ge g_3(x) \ge g_1(x)
\end{equation*}
holds similarly.

\subsection{Proof of Lemma \ref{lem:monotone}}

{
  Recall $h_\ell(x) = \sum_{\substack{T \subseteq [n]\\ |T| = \ell}}\prod_{j \in T} x_j p_j$ and 
$g_{t,d}(x) = \ssum\limits_{\ell=d}^{t+d}(-1)^{\ell - d}\binom{\ell}{d}h_\ell(x)$, where $t+d \le n$. 
For any fixed $d$ such that $0 \le d \le k-1$, we have
\begin{align}
g_{t+2,d}(x)-g_{t,d}(x)
&=\sum_{\ell=t+d+1}^{t+d+2}(-1)^{\ell-d}\binom{\ell}{d}\,h_\ell(x) \notag\\
&=(-1)^{t+1}\Big(\binom{q}{d}h_q(x)-\binom{q+1}{d}h_{q+1}(x)\Big), \label{eq:gt-diff}
\end{align}
where $q:=t+d+1$. Define
\[
\Phi_q(x):=\binom{q}{d}h_q(x)-\binom{q+1}{d}h_{q+1}(x),
\]
so \eqref{eq:gt-diff} can be written as
\begin{equation}\label{eq:star}
g_{t+2,d}(x)-g_{t,d}(x)=(-1)^{t+1}\,\Phi_{q}(x).
\end{equation}

If $h_q(x)=0$ for some $q$, then $h_{q'}(x)=0$ for all $q'\ge q$ (since each term in $h_{q'}$ is a
product of $q'$ factors from $\{x_jp_j\}$, and $h_q=0$ implies any $q$ factors are zero). Consequently, $\Phi_{q'}(x)=0$ for all
$q'\ge q$, and the desired monotonicity holds for all $q'\ge q$. Hence, it suffices to consider the case
where $h_q(x)>0,$ for all $q \le n$.

Note that $\Phi_q(x)\ge 0$ is equivalent to
\begin{equation}\label{eq:Phi-ineq}
\frac{h_q(x)}{h_{q+1}(x)}
\ge
\frac{\binom{q+1}{d}}{\binom{q}{d}}
=
\frac{q+1}{q+1-d}.
\end{equation}
Moreover, $h_\ell(x)$ is the $\ell$-th elementary symmetric polynomial in $(x_1p_1,\ldots,x_np_n)$.
By Newton's inequalities for elementary symmetric polynomials,
\[
\left(\frac{h_\ell(x)}{\binom{n}{\ell}}\right)^2 \ge 
\left(\frac{h_{\ell-1}(x)}{\binom{n}{\ell-1}}\right)
\left(\frac{h_{\ell+1}(x)}{\binom{n}{\ell+1}}\right),
\]
which implies
\[
\frac{h_\ell(x)}{h_{\ell+1}(x)} \ge 
\frac{h_{\ell-1}(x)}{h_\ell(x)}\cdot
\frac{\binom{n}{\ell}^2}{\binom{n}{\ell-1}\binom{n}{\ell+1}}
\ge
\frac{h_{\ell-1}(x)}{h_\ell(x)},
\]
where the last inequality follows from $\binom{n}{\ell}^2 \ge \binom{n}{\ell-1}\binom{n}{\ell+1}$.
Therefore, the ratio $R_\ell:=h_\ell(x)/h_{\ell+1}(x)$ is nondecreasing in $\ell$.
Moreover, the ratio $(\ell+1)/(\ell+1-d)$ is strictly decreasing in $\ell$. It follows from \eqref{eq:Phi-ineq} that if
$\Phi_{q}(x)\ge 0$ holds for some $q$, then $\Phi_{q'}(x)\ge 0$ holds for all $q'\ge q$.

Now assume $t$ is odd and $g_{t,d}(x)\le g_{t+2,d}(x)$. For any odd $\tau$ with
$t\le \tau\le n-2$, we have $(-1)^{\tau+1}=1$, so \eqref{eq:star} yields
\[
g_{\tau+2,d}(x)-g_{\tau,d}(x)=\Phi_{\tau+d+1}(x)\ge 0,
\]
and therefore $g_{\tau,d}(x)\le g_{\tau+2,d}(x)$.

Similarly, assume $t$ is even and $g_{t,d}(x)\ge g_{t+2,d}(x)$. \eqref{eq:star} yields
\[
g_{\tau+2,d}(x)-g_{\tau,d}(x)=-\Phi_{\tau+d+1}(x)\le 0,
\]
for any even $\tau$ with $t\le \tau\le n-2$.
}

\subsection{Proof of Proposition \ref{thm:monotone}}

From Lemma \ref{lem:monotone}, it follows that such $t_{\min}^{odd}$ and $t_{\min}^{even}$ exist. Therefore, for any odd $\tau$ such that $t_{\min}^{odd} \le \tau \le n-2$, we have $g_{\tau+2}(x) -g_\tau(x) = \sum_{d = 0}^{k-1} \left(g_{\tau+2,d}(x) -g_{\tau,d}(x)\right) \ge 0$. For even $\tau$, the statement holds similarly. 

\subsection{Proof of Lemma \ref{lem: k-e} }
   $X_{i_1} \subseteq X_{i_2}$ follows from  
    \begin{equation*}
    \mathbb{P}\left[\sum_{j \in [n]} \tilde{a}_{i_2 j} x_j \ge k_{i_2}\right] =    \mathbb{P}\left[\sum_{j \in [n]} \tilde{a}_{i_1 j} x_j \ge k_{i_2}\right]  \ge \mathbb{P}\left[\sum_{j \in [n]} \tilde{a}_{i_1 j} x_j \ge k_{i_1}\right] \ge 1-\epsilon_{i_1} \ge  1-\epsilon_{i_2},
    \end{equation*}
for any $x \in X_{i_1}$. 

\subsection{Proof of Lemma \ref{lem: relation} }
Similar to the proof of Lemma \ref{lem: k-e}, we have 
    \begin{equation*}
    \mathbb{P}\left[\sum_{j \in [n]} \tilde{a}_{i_2 j} x_j \ge k_{i_2}\right] \ge    \mathbb{P}\left[\sum_{j \in [n]} \tilde{a}_{i_1 j} x_j \ge k_{i_2}\right]  \ge \mathbb{P}\left[\sum_{j \in [n]} \tilde{a}_{i_1 j} x_j \ge k_{i_1}\right] \ge 1-\epsilon_{i_1} \ge  1-\epsilon_{i_2},
    \end{equation*}
for any $x \in X_{i_1}$. Thus, $X_{i_1} \subseteq X_{i_2}$.

\subsection{Proof of Lemma \ref{lem: covering-1} }
We denote the indices of sets in the corresponding candidate solution $\bar{x}$ by $D := \{j : \bar{x}_j = 1, j \in [n]\}$. Then 
    \begin{align*}
      \mathbb{P}\left[\sum_{j=1}^n \tilde{a}_j \bar{x}_j  \ge k \right] = \mathbb{P}\left[\sum_{j \in D} \tilde{a}_j  \ge k \right].
    \end{align*}
Note that the random variable $\sum_{j \in D} \tilde{a}_j \sim \text{Binomial}(d,p)$ since $a_j$'s are independent Bernoulli random variables with the same success probability $p$. Thus,  
\begin{equation*}
\mathbb{P}\left[\sum_{j=1}^n \tilde{a}_j \bar{x}_j  \ge k \right]=\sum_{\ell = k}^d\mathbb{P}\left[\sum_{j \in D} \tilde{a}_j  = \ell \right] = \sum_{\ell = k}^d\binom{d}{\ell}p^{\ell}(1-p)^{d-\ell}. 
\end{equation*}
\subsection{Proof of Lemma \ref{lem: covering-2} }
We consider the extended binomial coefficients by allowing $\alpha$ to be an arbitrary number (negative, rational, real), expressed as
\begin{equation*}
{\sbinom {\alpha }{k}}:=\textstyle\frac{\alpha (\alpha -1)(\alpha -2)\cdots (\alpha -k+1)}{k(k-1)(k-2)\cdots 1} \quad {\text{for }}k\in \mathbb {N} {\text{ and }} \alpha\in\mathbb{R}.
\end{equation*}
With this definition one has a generalization of the binomial formula 
\begin{equation*}
(x+y)^{\alpha}=\ssum _{\ell=0}^{+\infty}{ \binom{\alpha}{\ell}}x^{\ell}y^{\alpha-\ell}, \ \text{for any $x,y \in (0,1)$ and $\alpha \in \mathbb{R}$}. 
\end{equation*}
We still use $D := \{j : \bar{x}_j = 1, j \in [n]\}$ to denote the indices of sets in the corresponding candidate solution $\bar{x}$, and use the rising Pochhammer symbol $q^{(n)} = \binom{q+n-1}{n}n!$ in the proof. By Lemma \ref{lem: 1}, we have 
    \begin{align}
    \mathbb{P}\left[\ssum_{j=1}^n \tilde{a}_j \bar{x}_j  \ge k \right] &=\ssum_{\ell^=k}^n(-1)^{\ell - k}\binom{\ell-1}{\ell -k}  \ssum_{\substack{S \subseteq [n]\\ |S| = \ell}}\prod_{j \in S} \bar{x}_j p_j \\
    & =\ssum_{\ell^=k}^d(-1)^{\ell - k}\binom{\ell-1}{\ell -k}  \ssum_{\substack{S \subseteq D\\ |S| = \ell}}\prod_{j \in S} \bar{x}_j p_j \label{eq: thm5-0}\\
    &=\ssum_{\ell^=k}^d(-1)^{\ell - k}\binom{\ell-1}{\ell -k}  \binom{d}{\ell}p^\ell = \ssum_{\ell^=k}^d(-1)^{\ell - k}\binom{\ell-1}{k -1}  \binom{d}{\ell -1} \frac{d-\ell +1}{\ell} p^\ell \label{eq: thm5-1}\\
    & = \ssum_{\ell^=k}^d(-1)^{\ell - k}\binom{d}{k -1}  \binom{d-k+1}{d- \ell +1} \frac{d-\ell +1}{\ell} p^\ell \label{eq: thm5-2}\\
    & =k \binom{d}{k} \ssum_{\ell^=k}^d(-1)^{\ell - k}\binom{d-k}{\ell -k} \frac{p^\ell}{\ell}  = k \binom{d}{k} \ssum_{t^=0}^{d+k}(-1)^{t}\binom{d-k}{t} \frac{p^{k+t}}{k+t} \tag{*} \label{eq: thm5-3} \\ 
    & = \binom{d}{k}p^k \ssum_{t^=0}^{+\infty}\binom{t+k-d-1}{t} \frac{kp^{t}}{k+t} = \binom{d}{k}p^k \ssum_{t^=0}^{+\infty}(k-d)^{(t)} \cdot \frac{k}{k+t} \cdot \frac{p^{t}}{t!} \label{eq: thm5-4-5} \\ 
    & = \binom{d}{k}p^k \ssum_{t^=0}^{+\infty}(k-d)^{(t)} \cdot \frac{k^{(t)}}{(k+1)^{(t)}} \cdot \frac{p^{t}}{t!} = \binom{d}{k} p^k {}_{2}F_{1}(k-d,k;k+1;p) \tag{**} \label{eq: thm5-6} \\
    & = I_{p}(k,d-k+1) \tag{***} \label{eq: thm5-7}
    \end{align}
where the equation \eqref{eq: thm5-0} follows from the definition of $D$, the equation \eqref{eq: thm5-1} follows from the combinatorial identity $\binom {n}{k}={\frac {n}{k}}{\binom {n-1}{k-1}}$, the equation \eqref{eq: thm5-2} follows from $ {\binom {n}{h}}{\binom {n-h}{k}}={\binom {n}{k}}{\binom {n-k}{h}}$, the equation \eqref{eq: thm5-4-5} follows from $(-1)^k \binom{n}{k} = \binom{k-n-1}{k}$ and the definition of the rising Pochhammer symbol $q^{(n)}$, the equation \eqref{eq: thm5-7} follows from the fact $I_x(a,b)= x^a {}_{2}F_{1}(a,1-b;a+1,x)/(a \mathrm{B}(a,b)) $ and the equations \eqref{eq: thm5-3}, \eqref{eq: thm5-6} and \eqref{eq: thm5-7} indicate the conclusion.

Moreover, letting
\begin{equation*}
h(u) := k \binom{d}{k} \sum_{\ell^=k}^d(-1)^{\ell - k}\binom{d-k}{\ell -k} \frac{u^\ell}{\ell}= k \binom{d}{k} \sum_{t^=0}^{+\infty}(-1)^{t}\binom{d-k}{t} \frac{u^{k+t}}{k+t},
\end{equation*}
we have $h(0) =0$ and note that
\begin{equation*}
\frac{dh}{du} = k \binom{d}{k} u^{k-1} \sum_{t = 0}^{+\infty} (-1)^{t} \binom{d-k}{t} u^t = k\binom{d}{k}u^{k-1}(1-u)^{d-k}.
\end{equation*}
Hence,
\begin{equation*}
k\binom{d}{k}\int_{0}^pu^{k-1}(1-u)^{d-k}du = h(p) - h(0) = h(p),
\end{equation*}
which completes our proof.

\subsection{Proof of Lemma \ref{lem: feasible-sol} }

Let $D := \{j : \bar{x}_j = 1, j \in [n]\}$ and $\hat{D} := \{j : \hat{x}_j = 1, j \in [n]\}$. It is equivalent to prove $\mathbb{P}\left[ \sum_{j=1}^n \tilde{a}_j \hat{x}_j \ge k\right] \ge 1 - \epsilon$. Let $F(k;n,p)$ be the CDF of a binomial distribution. {It is well known that, for any fixed $k$ and $p$, the function $F(k; n, p)$ is monotone nonincreasing with respect to $n$, where $n \in \mathbb{N}$.} Note that the cover probability 
\begin{align*}
    \mathbb{P}[\ssum_{j=1}^n \tilde{a}_j \hat{x}_j  \ge k ] &= \mathbb{P}[\ssum_{j \in \hat{D}} \tilde{a}_j  \ge k ]  = \mathbb{P}[\ssum_{j = 1}^{\hat{d}} \tilde{a}_j  \ge k ]  \quad (\tilde{a}_j\text{'s are i.i.d random variables} ) \\ 
    & \ge \mathbb{P}[\ssum_{j = 1}^{d} \tilde{a}_j  \ge k ] = \mathbb{P}[\ssum_{j \in D} \tilde{a}_j  \ge k ]  \quad (\text{by monotone property of $F(k;n,p)$} ) \\ 
    & =  \mathbb{P}[\ssum_{j=1}^n \tilde{a}_j \bar{x}_j  \ge k ] \ge 1 -\epsilon. \quad (\text{$\bar{x}$ is a feasible solution})
\end{align*}
This completes our proof.
\subsection{Proof of Proposition \ref{prop: markov} }
First notice that $d \ge k$ since we cannot use less than $k$ sets to cover the item at least $k$ times with a predefined probability. Furthermore, by the Markov's inequality, we have 
    \begin{align*}
    1-\epsilon  \le \mathbb{P}\left[\sum_{j=1}^n \tilde{a}_j \bar{x}_j  \ge k \right]  \le \frac{\mathbb{E}\left[\sum_{j=1}^n \tilde{a}_j \bar{x}_j \right]}{k} = \frac{\sum_{j=1}^n \bar{x}_j \mathbb{E}\left[\tilde{a}_j  \right]}{k} = \frac{dp}{k}.
    \end{align*}
Therefore, $d : = \sum_{j \in [n]}\bar{x}_j  \ge \max\left\{k, \left\lceil\frac{k(1-\epsilon)}{p}\right\rceil \right\}.$

\subsection{Proof of Theorem \ref{thm: lower_bound} }
Let $x^*$ be an optimal solution to CC-SMCP. Then $q_i(x^*) := \mathbb{P}[\tilde{A}_{i}(\omega) x^* \le k_i - 1] \le \epsilon_i$ for each $i \in [m]$. By the definition of $X_\alpha^N$, $x^* \in X_\alpha^N$ if and only if no more than $\lfloor \alpha_i N \rfloor$ times the event $\{\tilde{A}_{i}(\omega) x^* \le k_i - 1\}$ happens in $N$ trials for each $ i \in [m]$. That is, for any $i \in [m]$, 
    \begin{equation*}
     \sum_{\omega \in \Omega}\mathbb{I}\left( \tilde{A}_{i}(\omega) x \ge k_i \right) \ge (1- \alpha_i)N \Leftrightarrow \sum_{\omega \in \Omega}\mathbb{I}\left( \tilde{A}_{i}(\omega) x \le k_i -1 \right) \le \lfloor \alpha_iN \rfloor. 
    \end{equation*}
Also, note that if $x^* \in  X_\alpha^N$, then $\nu_\alpha^N \le \nu^*$. Let $E_i$ be the event that $\{ \tilde{A}_{i}(\omega) x^* \le k_i - 1 \}$  happens at most $ \lfloor \alpha_i N \rfloor$ times. By Boole's inequality (union bound), we have 
\begin{equation*}
\mathbb{P}\left[\nu_\alpha^N \le \nu^*_\epsilon \right]  \ge \mathbb{P}\left[x^* \in X_\alpha^N\right] = \mathbb{P}\left[\cap_{i = 1}^m E_i\right] \ge 1- \sum_{i = 1}^m (1 - \mathbb{P}\left[E_i\right]).
\end{equation*}
Moreover, for any $i \in [m]$, 
\begin{equation*}
\mathbb{P}\left[E_i\right] =  \sum_{\ell = 0}^{\lfloor\alpha_i N \rfloor} \binom{N}{\ell} \left[q_i(x^*)\right]^\ell \left[1- q_i(x^*)\right]^{N- \ell} \ge \sum_{\ell = 0}^{\lfloor\alpha_i N \rfloor} \binom{N}{\ell} \epsilon_i^\ell (1- \epsilon_i)^{N- \ell}
\end{equation*}
where the inequality follows from the fact that the CDF of the binomial distribution, $F(k; n, p)$, is non-increasing with respect to $p$. Therefore,
\begin{equation*}
\mathbb{P}\left[\nu_\alpha^N \le \nu^*_\epsilon \right]  \ge 1 - \sum_{i = 1}^m \sum_{\ell = \lfloor\alpha_i N \rfloor +1 }^N \binom{N}{\ell} \epsilon_i^\ell (1- \epsilon_i)^{N- \ell} = 1 - \sum_{i = 1}^m I_{\epsilon_i}(\lfloor \alpha_i N \rfloor + 1, N - \lfloor \alpha_i N \rfloor).
\end{equation*}
Moreover, we define $Y_i$ to be a binomial random variable with parameters $(N, \epsilon_i)$. Using the Chernoff Bound in Theorem \ref{lem: chernoff}, we obtain if $\epsilon_i \le \alpha_i$ for some $i \in [m]$, then
\begin{equation*}
\mathbb{P}\left[Y_i \ge \alpha_i N \right] =  \mathbb{P}\left[Y_i \ge \epsilon_i N \left(1+ \left(\frac{\alpha_i}{\epsilon_i} - 1 \right) \right) \right] \le \exp \left(\frac{- (\alpha_i - \epsilon_i)^2 N}{ \alpha_i + \epsilon_i} \right) \le \exp(-\kappa_1 N).
\end{equation*}
where $\kappa_1:= \min_{i \in [m]} \left\{(\alpha_i -\epsilon_i)^2 /(\alpha_i + \epsilon_i)\right\} $. 
Therefore, if $\epsilon_i \le \alpha_i$ for each $i \in [m]$, then
\begin{equation*}
\mathbb{P}\left[\nu_\alpha^N \le \nu^*_\epsilon \right] \ge 1- \sum_{i = 1}^m\mathbb{P}\left[Y_i \ge \alpha_i N \right] \ge 1- \sum_{i = 1}^m \exp \left(\frac{- (\alpha_i - \epsilon_i)^2 N}{  \alpha_i + \epsilon_i} \right)  \ge 1-m \exp\left(-\kappa_1 N\right).
\end{equation*}

\subsection{Proof of Theorem \ref{thm: feasible_region} }
We define the indicator random variable $Y_{i, \omega}$ by $Y_{i, \omega} = 1$ if $\tilde{A}_{i}(\omega) x \ge k_i $ and $Y_{i, \omega} = 0$ otherwise. Recall that $X_\epsilon := \cap_{i \in [m]} X_{\epsilon, i}$ where $X_{\epsilon, i}:= \left\{x \in B: \mathbb{P}[\tilde{A}_{i} x \ge k_i] \ge 1 - \epsilon_i \right\}$. Then for any $x \in B \backslash X_{\epsilon, i}$, we obtain $\mathbb{E}[Y_{i, \omega}] = \mathbb{P}[\tilde{A}_{i}(\omega)x \ge k_i] < 1- \epsilon_i$. Let 
        \begin{equation*}
         X_{\alpha, i}^N := \left\{x \in B: \frac{1}{N} \sum_{\omega \in \Omega} \mathbb{I}\left( \tilde{A}_{i}(\omega) x \ge k_i \right) \geq 1-\alpha_i\right\}, \forall i \in [m].
        \end{equation*}
Note that $X_\alpha^N = \cap_{i \in [m]} X_{\alpha, i}^N$ and $x \in  X_{\alpha, i}^N$ if and only if $(1/N) \sum_{\omega \in \Omega} Y_{i, \omega} \ge 1 - \alpha_i $. Hoeffding's inequality yields, for any $x \in B \backslash X_{\epsilon, i}$, 
\begin{align*}
\mathbb{P}[x \in X_{\alpha, i}^N] & = \mathbb{P}\left[ \frac{1}{N} \sum_{\omega \in \Omega} Y_{i, \omega}] \ge 1 - \alpha_i\right] \\ 
& \le \mathbb{P}\left[\sum_{\omega \in \Omega} \left( Y_{i, \omega} - \mathbb{E}\left[Y_{i, \omega}\right] \right)  \ge N(\epsilon_i - \alpha_i) \right] \le \exp \left\{- 2 N (\epsilon_i - \alpha_i)^2  \right\}.
\end{align*}
Therefore, we can obtain an upper bound for the probability 
\begin{align*}
\mathbb{P}\left[X_\alpha^N \nsubseteq X_\epsilon \right] & = \mathbb{P}\left[\exists x \in X_\alpha^N \text{ where } x \in B \backslash X_\epsilon \right] = \mathbb{P}\left[\bigcup_{x \in B \backslash X_\epsilon} \left\{ x \in X_\alpha^N\right\}\right] \\
& = \mathbb{P}\left[\bigcup_{i \in [m]}\bigcup_{x \in B \backslash X_{\epsilon, i}} \left\{ x \in \bigcap_{i \in [m]} X_{\alpha, i}^N\right\}\right]  \le \sum_{i \in [m]} \sum_{x \in B \backslash X_{\epsilon, i}}\mathbb{P}\left[x \in X_{\alpha, i}^N\right] \\
& \le \sum_{i \in [m]} \sum_{x \in B \backslash X_{\epsilon, i}} \exp \left\{- 2 N (\epsilon_i - \alpha_i)^2  \right\} \le m |B \backslash X_\epsilon| \exp\left(-\kappa_2 N\right)
\end{align*}
where $\kappa_2 :=2 \min_{i \in [m]} (\epsilon_i - \alpha_i)^2$. This completes our proof.

\subsection{Proof of Theorem \ref{thm: SAA_opt} }
We will use the following lemma in the proof:  
\begin{lem} \label{lem: equivalent}
Let $X_{\epsilon}$ be the feasible region of CC-SMCP. Then
    \begin{equation*}
    X_{\epsilon} = X_{\overline{\alpha}}',
    \end{equation*}
where $X_{\overline{\alpha}}':= \cap_{i \in [m]} X_{\overline{\alpha},i}'$ and $X_{\overline{\alpha},i}' : = \left\{x \in B: \mathbb{P}[\tilde{A}_{i} x \le k_i -1] < \overline{\alpha}_i \right\}$ for each $i \in [m]$.
\end{lem}

\textit{Proof (of Lemma \ref{lem: equivalent}).} We prove that $X_{\epsilon,i} = X_{\overline{\alpha},i}'$ for each $i \in [m]$. Clearly, $X_{\epsilon,i} \subseteq X_{\overline{\alpha},i}'$ for any $i \in [m]$ since $\mathbb{P}[\tilde{A}_{i} x \le k_i -1] \le \epsilon_i < \overline{\alpha}_i$ for every $x \in X_{\epsilon,i}$. Conversely, if $x \in X_{\alpha,i}'$ for a fixed $i$, then by the definition of $\overline{\alpha}_i$, we have $x \in X_{\epsilon, i}$. \hfill $\blacksquare$ 

\textit{Proof (of Theorem \ref{thm: SAA_opt}).} By the definition of $\underline{\alpha}_i$, we have $\nu_{\underline{\alpha}}^* = \nu^*_\epsilon$. Theorem \ref{thm: lower_bound} implies 
    \begin{equation*}
     \mathbb{P}\left[\nu_\epsilon^N > \nu_\epsilon^*\right] = \mathbb{P}\left[\nu_\epsilon^N > \nu_{\underline{\alpha}}^*\right] \le m \exp\left(-\min_{i \in [m]} \left\{ \frac{(\underline{\alpha}_i -\epsilon_i)^2 N}{(\underline{\alpha}_i + \epsilon_i)} \right\} \right).
    \end{equation*}
    Further, observe that the proof of Theorem \ref{thm: feasible_region} can be modified to show the slightly stronger result:
    \begin{equation*}
     \mathbb{P}\left[X_\alpha^N \subseteq X_\epsilon' \right] \ge 1 - m|B \backslash X_\epsilon'| \exp \left\{ -2N \min_{i \in [m]} (\epsilon_i - \alpha_i)^2 \right\},
    \end{equation*}
where $X_{\epsilon}':= \cap_{i \in [m]} X_{\epsilon,i}'$ and $X_{\epsilon,i}' : = \left\{x \in B: \mathbb{P}[\tilde{A}_{i} x \le k_i -1] < \epsilon_i \right\}$ for each $i \in [m]$. (In the proof, we consider each $x \in X \backslash X_{\epsilon}'$ and replace $X_\epsilon$ with $X_{\epsilon}'$. The remainder of the proof is identical to the original proof.) Applying this result yields 
    \begin{equation*}
     \mathbb{P}\left[X_\epsilon^N \subseteq X_{\overline{\alpha}}' \right] \ge 1 - m|B \backslash X_{\overline{\alpha}}'| \exp \left\{ -2N \min_{i \in [m]} (\overline{\alpha}_i  - \epsilon_i)^2 \right\}.
    \end{equation*}
Besides, if $X_\epsilon^N \subseteq X_\epsilon $ then we have $\nu_{\epsilon}^N \ge \nu^*_\epsilon$. By Lemma \ref{lem: equivalent}, it follows that 
\begin{equation*}
\mathbb{P}\left[\nu_\epsilon^N < \nu_\epsilon^*\right] \le \mathbb{P}\left[X_\epsilon^N \nsubseteq X_\epsilon \right] = \mathbb{P}\left[X_\epsilon^N \nsubseteq X_{\overline{\alpha}}' \right] \le m|B \backslash X_{\overline{\alpha}}'| \exp \left\{ -2N \min_{i \in [m]} (\overline{\alpha}_i  - \epsilon_i)^2 \right\}.
\end{equation*}
Therefore,
\begin{align*}
\mathbb{P}\left[\nu_\epsilon^N \neq \nu_\epsilon^*\right] & \le \mathbb{P}\left[\nu_\epsilon^N > \nu_\epsilon^*\right] + \mathbb{P}\left[\nu_\epsilon^N < \nu_\epsilon^*\right] \\
&\le m \exp\left(-\min_{i \in [m]} \left\{ \frac{(\underline{\alpha}_i -\epsilon_i)^2 N}{(\underline{\alpha}_i + \epsilon_i)} \right\} \right) + m |B \backslash X_\epsilon| \exp \left\{ -2N \min_{i \in [m]} (\overline{\alpha}_i  - \epsilon_i)^2 \right\} \\
& \le  m\left(|B \backslash X_\epsilon| + 1\right)\exp\{-\kappa_3 N\}.
\end{align*}

\subsection{Proof of Theorem \ref{thm: IS_OPT} }
To prove the theorem, we need the following lemma:  
\begin{lem}\label{lem: transform}
For $n \ge 1$, let $p_j, j = 1, \ldots, n$ be numbers such that $1>p_1 \ge \cdots \ge p_n >0$. Given an integer $w$ with $1 \le w \le n$, consider the following problem:
\begin{equation}
\min _{\lambda \in \mathbb{R}_{+}^n} \max _{\substack{z_j \in\{0,1\}^n \\ \sum_j z_j=w}}\sum_{j=1}^n z_j \lambda_j+\sum_{j=1}^n \log \left(e^{-\lambda_j} p_j+\left(1-p_j\right)\right) .
\end{equation}
Then, there exists an optimal solution to the above problem that satisfies $\lambda_1 \ge \cdots \ge \lambda_n \ge 0$.
\end{lem}  \vspace{-1.2em}

\textit{Proof (of Lemma \ref{lem: transform})}. Suppose that $\lambda := (\lambda_1,\ldots,\lambda_n)$ is an optimal solution to the problem and there exists some $j < n$ such that $\lambda_j < \lambda_{j+1}$. We define $\bar{\lambda}$ as $\bar{\lambda}_j:= \lambda_{j+1}$, $\bar{\lambda}_{j+1}:= \lambda_{j}$ and $\bar{\lambda}_\ell:= \lambda_\ell$ for $\ell \neq \{j,j+1\}$. We will show that $\bar{\lambda}$ has no worse function value than $\lambda$. First notice that 
\begin{equation*}
\max_{\substack{z_j \in\{0,1\}^n \\ \sum_j z_j=w}}\sum_{j=1}^n z_j \lambda_j=\max_{\substack{z_j \in\{0,1\}^n \\ \sum_j z_j=w}}\sum_{j=1}^n z_j \bar{\lambda}_j,
\end{equation*}
which is equal to the sum of the largest $w$ components of the vector $\lambda$. Therefore, we only need to compare the remaining part of the objective function. Define $\Delta$ as the difference in the objective function between $\lambda$ and $\bar{\lambda}$, i.e., 
\begin{align*}
\Delta &= \sum_{j=1}^n z_j \lambda_j+\sum_{j=1}^n \log \left(e^{-\lambda_j} p_j+\left(1-p_j\right)\right) - \left(\sum_{j=1}^n z_j \bar{\lambda}_j+\sum_{j=1}^n \log \left(e^{-\bar{\lambda}_j} p_j+\left(1-p_j\right)\right)\right) \\
&= \sum_{j=1}^n \log \left(e^{-\lambda_j} p_j+\left(1-p_j\right)\right) - \sum_{j=1}^n \log \left(e^{-\bar{\lambda}_j} p_j+\left(1-p_j\right)\right) \\
& = \log \left(e^{-\lambda_j} p_j+\left(1-p_j\right)\right)+ \log \left(e^{-\lambda_{j+1}} p_{j+1}+\left(1-p_{j+1}\right)\right) \\
& \quad - \log \left(e^{-\bar{\lambda}_j} p_j + \left(1-p_j\right)\right) - \log \left(e^{-\bar{\lambda}_{j+1}} p_{j+1}+\left(1-p_{j+1}\right)\right) 
\end{align*}
Since $\bar{\lambda}_j = \lambda_{j+1}$ and $\bar{\lambda}_{j+1} = \lambda_j$, it follows that
\begin{align*}
    \Delta & = \log \left(\frac{e^{-\lambda_j} p_j+\left(1-p_j\right)}{e^{-\lambda_{j+1}} p_j+\left(1-p_j\right)} \right) - \log \left(\frac{e^{-\lambda_j} p_{j+1}+\left(1-p_{j+1}\right)}{e^{-\lambda_{j+1}} p_{j+1}+\left(1-p_{j+1}\right)} \right) \\
    & = \log \left(\frac{e^{-\lambda_j} - e^{-\lambda_{j+1}}}{e^{-\lambda_{j+1}} + \frac{1}{p_j}-1}+1 \right) - \log \left(\frac{e^{-\lambda_j} - e^{-\lambda_{j+1}}}{e^{-\lambda_{j+1}} + \frac{1}{p_{j+1}}-1}+1 \right).   
\end{align*}
Note that the term inside the log is positive, since $-\lambda_j>-\lambda_{j+1}$. Further, since $p_j \ge p_{j+1}$, it follows that $1/p_j -1 \le 1/p_{j+1} -1 $. Therefore, we conclude that $\Delta \ge 0$. \hfill $\blacksquare$

\textit{Proof (of Theorem \ref{thm: IS_OPT}).} We consider obtaining an optimal solution by minimizing 
\begin{equation*}
\log (B_x(\lambda)) = \max_{\tilde{a}: \sum_{j = 1}^n \tilde{a}_jx_j < k} \sum_{j = 1}^n \lambda_j \tilde{a}_j + \sum_{j =1}^n \log \left(e^{-\lambda_j}p_j + (1 - p_j) \right).
\end{equation*}
We consider a vector $\tilde{a}:=(z_1,z_2,\ldots,z_n) \in \{0,1\}^n$ such that $z_j =1$ if $x_j = 0$, and $\sum_{j=1}^n z_j = n-u +k -1$. Without loss of generality, we assume the index set $\{j: z_j = 1\}$ corresponds to $\{1,\ldots,n-u+k-1\}$. By lemma \ref{lem: transform}, it follows that minimizing $\log(B_x(\lambda))$ over $\lambda \ge 0$ amounts to solving the following problem:
\begin{subequations}
\begin{align}
\min _{\lambda \in \mathbb{R}^n_+} \phi(\lambda):=& \sum_{j=1}^{n - u + k - 1} \lambda_j+\sum_{j=1}^n \log \left(e^{-\lambda_j} p_j+\left(1-p_j\right)\right) \\
&\lambda_j \ge \lambda_{j+1} \quad j=1 \ldots n-1  \label{eq: mu1}\\
&\lambda_n \geq 0 \label{eq: mu2}
\end{align} \label{eq: phi}
\end{subequations}
Note that the objective function $\phi(\lambda)$ is strictly convex in $\lambda$. In fact, its second derivatives are
\begin{equation*}
\frac{\partial^2 \phi}{\partial \lambda_j^2}=\frac{e^{-\lambda_j} p_j\left(1-p_j\right)}{\left(e^{-\lambda_j} p_j+\left(1-p_j\right)\right)^2}>0, \quad \frac{\partial^2 \phi}{\partial \lambda_{j_1} \partial \lambda_{j_2}}=0.
\end{equation*}
which implies the above problem \eqref{eq: phi} has a unique optimal solution that can be found by using Karush-Kuhn-Tucker (KKT) conditions:
\begin{subequations}
\begin{align}
1_{(1 \le n-u+k -1)} - \frac{p_{j}e^{-\lambda_{j}}}{e^{-\lambda_j} p_j+\left(1-p_j\right)}  - \mu_1 &= 0 \label{eq: kkt-1}\\
1_{(j \le n-u+k -1)} - \frac{p_{j}e^{-\lambda_{j}}}{e^{-\lambda_j} p_j+\left(1-p_j\right)} + \mu_{j-1} - \mu_j &= 0, ~~  j = 2,\ldots, n \label{eq: kkt-2}\\
\mu_j (\lambda_{j+1} - \lambda_j) &=0, ~~ j = 1,\ldots, n - 1 \label{eq: kkt-3}\\
\mu_n\lambda_n &= 0 \label{eq: kkt-4}\\ 
\mu_j &\ge 0, ~~ j = 1, \ldots, n  \label{eq: kkt-5}
\end{align}
\end{subequations}
where $(\mu_j)_{j \in [n-1]}$ is the Lagrangian multipliers of constraints \eqref{eq: mu1} and $\mu_n$ is the Lagrangian multiplier of constraint \eqref{eq: mu2}. 

We consider a particular choice of vectors \textbf{$\mu$} and $\mathbf{\lambda}$ defined as follows:
\begin{align}
\mu_n &:= 0, \label{eq: mu_n}\\
\mu_j &:= \min\{j, n-u+k-1\} - \ssum_{\ell = 1}^j \frac{e^{-\lambda^*}}{e^{-\lambda^*}p_\ell +(1- p_\ell)}, ~~ j = 1\ldots, n-1 \label{eq: mu_j}\\
\lambda_1 &= \cdots = \lambda_n := \lambda^* \label{eq: lambda}
\end{align}
where $\lambda^* \in \mathbb{R}_{+}$ solves the following equation
\begin{equation}
\eta(\lambda^*): = \sum_{j=1}^n\frac{e^{-\lambda^*}p_j}{e^{-\lambda^*}p_j+(1-p_j)}= n-u+k-1. \label{eq: lambda_star}
\end{equation}
Note that we can always find such $\lambda^*$, since the function $\eta(\lambda)$ is continuous and decreasing, and
\begin{align*}
\eta(0) &= \sum_{j = 1}^n{p}_j \ge \sum_{j = 1}^n p_jx_j >n-u+k -1, \\
\lim_{\lambda \to +\infty}\eta(\lambda) & = 0 < n-u+k -1.
\end{align*}
Then we claim that $\mathbf{\mu}$ and $\mathbf{\lambda}$ satisfy the KKT conditions \eqref{eq: kkt-1} -- \eqref{eq: kkt-5}. First note that equations \eqref{eq: mu_j} imply \eqref{eq: kkt-1} and \eqref{eq: kkt-2}. Equations \eqref{eq: kkt-3} follow from \eqref{eq: lambda}, Equation \eqref{eq: kkt-4} follows from \eqref{eq: mu_n}. Finally, \eqref{eq: lambda_star} implies that  
\begin{equation*}
\sum_{\ell=1}^j\frac{e^{-\lambda^*}p_j}{e^{-\lambda^*}p_j+(1-p_j)} <n-u+k-1, ~~ j = 1,\ldots, n -1. 
\end{equation*}
Additionally, since each term on the left-hand side of the above summation is less than $1$, we have 
\begin{equation*}
\sum_{\ell=1}^j\frac{e^{-\lambda^*}p_j}{e^{-\lambda^*}p_j+(1-p_j)} <j, ~~ j = 1, \ldots, n.
\end{equation*}
The above two observations yield inequalities \eqref{eq: kkt-5}.

\section{Additional Experiments} \label{sec: additional_experiments}

\subsection{The performance of checking infeasibility} 

In this section, we investigate computational performance of four approaches for checking infeasibility by reducing the value of the success probability $p_{ij}$. Similar to the settings in Section \ref{sec: overall_performance}, we select up to $12$ columns at random for each row  $i$ where $k_i \ge 2$. However, this time we uniformly sample the success probability $p_{ij} = P[\tilde{a}_{ij} = 1]$ from the range $[0.2, 0.6]$. For the remaining $n-n'$ columns in row $i$, we still assign a success probability $p_{ij}$ of $0$. All other settings remain consistent with those in Section \ref{sec: overall_performance}. By reducing the success probability $p_{ij}$, we can generate infeasible instances with a high probability. If an infeasible instance cannot be obtained for fixed $n,m$, and $\epsilon$, then this process should be repeated until an infeasible instance is found.

\begin{table}[htbp]
  \centering
  \caption{Computational performance of methods for checking infeasibility ($p_{ij} \in [0.2, 0.6], N = 200$)}
  \scriptsize
       \begin{tabular}{ccccccccc}
\cmidrule{1-9}    \multicolumn{3}{c}{Parameters} & \multicolumn{2}{c}{OA-I} & \multicolumn{2}{c}{OA-II} & SAA   & IS \\
    $n$   & $m$   & $\epsilon$ & Time  & Iterations & Time  & Iterations & Time  & Time \\
    \midrule
    30    & 10    & 0.05  & 0.26  & 2     & 0.04  & 2     & 0.00  & 0.00 \\
    30    & 10    & 0.1   & 0.13  & 2     & 0.05  & 2     & 0.00  & 0.00 \\
    30    & 20    & 0.05  & 1.30  & 2     & 0.88  & 2     & 0.01  & 0.01 \\
    30    & 20    & 0.1   & 0.81  & 2     & 0.93  & 2     & 0.01  & 0.01 \\
    30    & 30    & 0.05  & 0.37  & 2     & 0.18  & 2     & 0.01  & 0.01 \\
    30    & 30    & 0.1   & 0.46  & 2     & 0.20  & 2     & 0.01  & 0.01 \\
    30    & 50    & 0.05  & 2.90  & 2     & 1.62  & 2     & 0.02  & 0.02 \\
    30    & 50    & 0.1   & 2.66  & 2     & 1.24  & 2     & 0.02  & 0.02 \\
    30    & 100   & 0.05  & 0.00  & 1     & 0.01  & 1     & 0.04  & 0.03 \\
    30    & 100   & 0.1   & 0.00  & 1     & 0.00  & 1     & 0.04  & 0.03 \\
    30    & 150   & 0.05  & 0.01  & 1     & 0.01  & 1     & 0.05  & 0.05 \\
    30    & 150   & 0.1   & 0.01  & 1     & 0.01  & 1     & 0.05  & 0.05 \\
    50    & 50    & 0.05  & 1.18  & 2     & 0.35  & 2     & 0.02  & 0.02 \\
    50    & 50    & 0.1   & 0.32  & 2     & 0.35  & 2     & 0.02  & 0.02 \\
    50    & 100   & 0.05  & 0.00  & 1     & 0.01  & 1     & 0.05  & 0.07 \\
    50    & 100   & 0.1   & 0.00  & 1     & 0.01  & 1     & 0.05  & 0.04 \\
    50    & 150   & 0.05  & 1.69  & 2     & 1.03  & 2     & 0.08  & 0.07 \\
    50    & 150   & 0.1   & 1.60  & 2     & 1.09  & 2     & 0.07  & 0.07 \\
    50    & 200   & 0.05  & 0.01  & 1     & 0.01  & 1     & 0.11  & 0.10 \\
    50    & 200   & 0.1   & 0.01  & 1     & 0.01  & 1     & 0.11  & 0.10 \\
    50    & 250   & 0.05  & 0.01  & 1     & 0.01  & 1     & 0.15  & 0.13 \\
    50    & 250   & 0.1   & 0.01  & 1     & 0.01  & 1     & 0.15  & 0.13 \\
    50    & 300   & 0.05  & 0.01  & 1     & 0.02  & 1     & 0.18  & 0.19 \\
    50    & 300   & 0.1   & 0.01  & 1     & 0.02  & 1     & 0.18  & 0.17 \\
    100   & 50    & 0.05  & 0.00  & 1     & 0.00  & 1     & 0.04  & 0.03 \\
    100   & 50    & 0.1   & 0.00  & 1     & 0.00  & 1     & 0.04  & 0.04 \\
    100   & 100   & 0.05  & 0.01  & 1     & 0.01  & 1     & 0.09  & 0.07 \\
    100   & 100   & 0.1   & 0.01  & 1     & 0.01  & 1     & 0.08  & 0.08 \\
    100   & 150   & 0.05  & 13.27 & 2     & 3.17  & 2     & 0.14  & 0.11 \\
    100   & 150   & 0.1   & 9.34  & 2     & 3.18  & 2     & 0.14  & 0.11 \\
    300   & 50    & 0.05  & 0.00  & 1     & 0.00  & 1     & 0.08  & 0.07 \\
    300   & 50    & 0.1   & 0.00  & 1     & 0.00  & 1     & 0.08  & 0.07 \\
    300   & 100   & 0.05  & 740.40 & 2     & 955.05 & 2     & 0.18  & 0.15 \\
    300   & 100   & 0.1   & 747.88 & 2     & 162.76 & 2     & 0.18  & 0.15 \\
    300   & 150   & 0.05  & 0.01  & 1     & 0.01  & 1     & 0.27  & 0.23 \\
    300   & 150   & 0.1   & 487.47 & 2     & 1406.14 & 2     & 0.25  & 0.25 \\
    300   & 200   & 0.05  & 0.01  & 1     & 0.02  & 1     & 0.36  & 0.31 \\
    300   & 200   & 0.1   & 0.01  & 1     & 0.02  & 1     & 0.35  & 0.31 \\
    300   & 250   & 0.05  & 0.02  & 1     & 0.02  & 1     & 0.44  & 0.45 \\
    300   & 250   & 0.1   & 0.02  & 1     & 0.02  & 1     & 0.45  & 0.42 \\
    300   & 300   & 0.05  & 0.02  & 1     & 0.03  & 1     & 0.56  & 0.81 \\
    300   & 300   & 0.1   & 0.02  & 1     & 0.03  & 1     & 0.54  & 0.47 \\
    \bottomrule
    \end{tabular}%
  \label{tab:infeasibility}%
\end{table}

Table \ref{tab:infeasibility} presents the computational performance of four approaches for detecting infeasibility of instances. Observe that SAA and IS are efficient in detecting infeasibility, while OA methods perform well in all cases except for some instances when $n = 300$. It is worth noting that even if the corresponding SAA or IS problem is infeasible, this does not necessarily mean that the original CC-SMCP problem is infeasible, as reflected in Table \ref{tab: comparison}.

\subsection{Experiments on special cases}

In this section, we conduct a series of experiments on the special cases outlined in Section \ref{sec: special_cases}. For the PSC problem, extensive experiments have already been conducted in the previous literature \citep{haight2000integer,fischetti2012cutting,ahmed2013probabilistic}. Therefore, our main focus is on the case when cover probabilities of each item are equal, as discussed in Section \ref{sec: case 2}.

In the following experiments, the covering parameter $k_i$ is randomly selected from \{2,3\} and all other experimental settings remain consistent with those in Section \ref{sec: overall_performance}. Given that we are working with a \textit{sparse} probability matrix $P$ in our experimental setup, this implies that for a given item $i$, it can either be covered with a set with an equal probability $p_i$ or not be covered at all by that set. This assumption extends the case discussed in Section \ref{sec: case 2} and motivates us to avoid computing $\bar{d}$. 
Instead, we aim to address the following problem: 
\begin{equation*}
\min\left\{\ssum_{j \in [n]} c_j x_j : \ssum_{j \in J_i}x_j \ge \bar{d}_i, \forall i \in [m], x \in B \subseteq \{0,1\}^n\right\},
\end{equation*}
where $J_i = \{j \in [n]\ :\ p_{ij} \neq 0\}$. This problem takes into account the nonzero coefficients of the probability matrix for each item $i$ and integrates the corresponding subproblems together.

Table \ref{tab: special_case} presents the computational performance of three methods when applied to special cases. In this table, ``SC'' refers to the method described in Section \ref{sec: case 2} that is used to address these special cases. ``BS'', ``Solving'' and ``Total'' indicate the time spent on binary search, the time required to solve the reformulated problem, and the overall time (BS+ Solving) for the SC method, respectively.

From Table \ref{tab: special_case}, we can see that the SC method successfully solves all instances very quickly, indicating that the reformulation technique introduced in Section \ref{sec: case 2} can significantly reduce solution time compared to the general case in Section \ref{sec: overall_performance}. Additionally, the SC method outperforms sampling-based methods in almost all cases. This may imply that sampling-based methods are not very effective for special cases when $n$ or $m$ increases.  

\begin{table}[htbp]
  \centering
  \caption{Computational performance of methods on special cases ($p_{ij} \in [0.9,1], N = 200$)}
  \resizebox{\textwidth}{!}{%
     \begin{tabular}{ccccccccccccc}
    \toprule
    \multicolumn{3}{c}{Parameters} & \multicolumn{4}{c}{SC}        & \multicolumn{3}{c}{SAA} & \multicolumn{3}{c}{IS} \\
    $n$   & $m$   & $\epsilon$ & BS & Solving & Total & Val   & Time  & Val   & Gap(\%) & Time  & Val   & Gap(\%) \\
    \midrule
    30    & 10    & 0.05  & 0.00  & 0.00  & 0.00  & 9     & 0.41  & 9     & 0.00  & 1.04  & 9     & 0.00 \\
    30    & 10    & 0.1   & 0.00  & 0.00  & 0.00  & 8     & 0.29  & 7     & -12.50 & 0.10  & 7     & -12.50 \\
    30    & 20    & 0.05  & 0.00  & 0.00  & 0.00  & 15    & 3.28  & 16    & 6.67  & 1.77  & 15    & 0.00 \\
    30    & 20    & 0.1   & 0.00  & 0.00  & 0.00  & 14    & 5.41  & 13    & -7.14 & 3.34  & 15    & 7.14 \\
    30    & 30    & 0.05  & 0.00  & 0.00  & 0.00  & 17    & 1.48  & 18    & 5.88  & 0.71  & 17    & 0.00 \\
    30    & 30    & 0.1   & 0.00  & 0.00  & 0.00  & 17    & 25.02 & 15    & -11.76 & 1.77  & 17    & 0.00 \\
    30    & 50    & 0.05  & 0.00  & 0.00  & 0.00  & 21    & 2.01  & 21    & 0.00  & 1.84  & 22    & 4.76 \\
    30    & 50    & 0.1   & 0.00  & 0.00  & 0.00  & 20    & 7.75  & 19    & -5.00 & 2.34  & 19    & -5.00 \\
    30    & 100   & 0.05  & 0.01  & 0.00  & 0.01  & 24    & 0.78  & 25    & 4.17  & 1.80  & 25    & 4.17 \\
    30    & 100   & 0.1   & 0.01  & 0.00  & 0.01  & 24    & 2.09  & 24    & 0.00  & 2.90  & 24    & 0.00 \\
    30    & 150   & 0.05  & 0.03  & 0.00  & 0.03  & 25    & 3.26  & 26    & 4.00  & 1.37  & 24    & -4.00 \\
    30    & 150   & 0.1   & 0.01  & 0.00  & 0.01  & 24    & 2.45  & 23    & -4.17 & 2.58  & 24    & 0.00 \\
    50    & 50    & 0.05  & 0.00  & 0.00  & 0.00  & 27    & 7.71  & 27    & 0.00  & 2.51  & 27    & 0.00 \\
    50    & 50    & 0.1   & 0.00  & 0.00  & 0.00  & 25    & 11.00 & 25    & 0.00  & 3.22  & 26    & 4.00 \\
    50    & 100   & 0.05  & 0.01  & 0.00  & 0.01  & 35    & 0.03  & INF   & INF   & 0.02  & INF   & INF \\
    50    & 100   & 0.1   & 0.01  & 0.00  & 0.01  & 33    & 7.11  & 34    & 3.03  & 2.88  & 33    & 0.00 \\
    50    & 150   & 0.05  & 0.01  & 0.00  & 0.01  & 39    & 0.03  & INF   & INF   & 0.04  & INF   & INF \\
    50    & 150   & 0.1   & 0.01  & 0.00  & 0.01  & 38    & 14.95 & 35    & -7.89 & 3.15  & 38    & 0.00 \\
    50    & 200   & 0.05  & 0.01  & 0.00  & 0.01  & 39    & 13.97 & 37    & -5.13 & 4.10  & 39    & 0.00 \\
    50    & 200   & 0.1   & 0.01  & 0.05  & 0.06  & 35    & 152.25 & 33    & -5.71 & 11.23 & 35    & 0.00 \\
    50    & 250   & 0.05  & 0.01  & 0.00  & 0.01  & INF   & 0.12  & INF   & 0.00  & 0.07  & INF   & INF \\
    50    & 250   & 0.1   & 0.01  & 0.06  & 0.07  & 40    & 30.12 & 40    & 0.00  & 6.87  & 40    & 0.00 \\
    50    & 300   & 0.05  & 0.02  & 0.00  & 0.02  & INF   & 0.20  & INF   & 0.00  & 0.05  & INF   & INF \\
    50    & 300   & 0.1   & 0.02  & 0.00  & 0.02  & 42    & 20.48 & 41    & -2.38 & 6.90  & 42    & 0.00 \\
    100   & 50    & 0.05  & 0.00  & 0.00  & 0.01  & 36    & 49.49 & 35    & -2.78 & 6.87  & 35    & -2.78 \\
    100   & 50    & 0.1   & 0.00  & 0.09  & 0.09  & 33    & 155.75 & 31    & -6.06 & 27.11 & 33    & 0.00 \\
    100   & 100   & 0.05  & 0.01  & 0.00  & 0.01  & 58    & 10.91 & 55    & -5.17 & 4.62  & 56    & -3.45 \\
    100   & 100   & 0.1   & 0.01  & 0.00  & 0.01  & 51    & 56.68 & 50    & -1.96 & 14.32 & 50    & -1.96 \\
    100   & 150   & 0.05  & 0.01  & 0.05  & 0.06  & 59    & 103.79 & 58    & -1.69 & 32.99 & 58    & -1.69 \\
    100   & 150   & 0.1   & 0.01  & 0.01  & 0.01  & 53    & 1040.21 & 53    & 0.00  & 198.95 & 54    & 1.89 \\
    300   & 50    & 0.05  & 0.00  & 0.01  & 0.01  & 71    & 836.67 & 70    & -1.41 & 7.28  & 67    & -5.63 \\
    300   & 50    & 0.1   & 0.00  & 0.06  & 0.06  & 66    & 948.75 & 64    & -3.03 & 102.95 & 64    & -3.03 \\
    300   & 100   & 0.05  & 0.00  & 0.01  & 0.01  & 105   & 0.03  & INF   & INF   & 190.44 & 105   & 0.00 \\
    300   & 100   & 0.1   & 0.01  & 0.05  & 0.06  & 100   & -     & 97    & -3.00 & - & 98    & -2.00 \\
    300   & 150   & 0.05  & 0.01  & 0.07  & 0.08  & 135   & 0.06  & INF   & INF   & 0.02  & INF   & INF \\
    300   & 150   & 0.1   & 0.01  & 0.03  & 0.05  & 126   & 0.03  & INF   & INF   & 267.93 & 123   & -2.38 \\
    300   & 200   & 0.05  & 0.01  & 0.01  & 0.02  & 148   & 0.08  & INF   & INF   & 256.41 & 147   & -0.68 \\
    300   & 200   & 0.1   & 0.01  & 0.01  & 0.02  & 138   & -     & 139   & 0.72  & 1068.16 & 139   & 0.72 \\
    300   & 250   & 0.05  & 0.01  & 0.00  & 0.01  & INF   & 0.09  & INF   & 0.00  & 0.05  & INF   & INF \\
    300   & 250   & 0.1   & 0.01  & 0.10  & 0.12  & 135   & -     & 134   & -0.74 & - & 132   & -2.22 \\
    300   & 300   & 0.05  & 0.01  & 0.09  & 0.10  & 159   & 0.09  & INF   & INF   & 0.05  & INF   & INF \\
    300   & 300   & 0.1   & 0.02  & 0.09  & 0.11  & 149   & -     & 146   & -2.01 & - & 147   & -1.34 \\
    \bottomrule
    \end{tabular}%
    }
  \label{tab: special_case}%
\end{table}

\subsection{The effect of risk parameter $\epsilon_i$}

\begin{figure}[htp]
\centering
\begin{tabular}{cc}
\includegraphics[width=.4\textwidth]{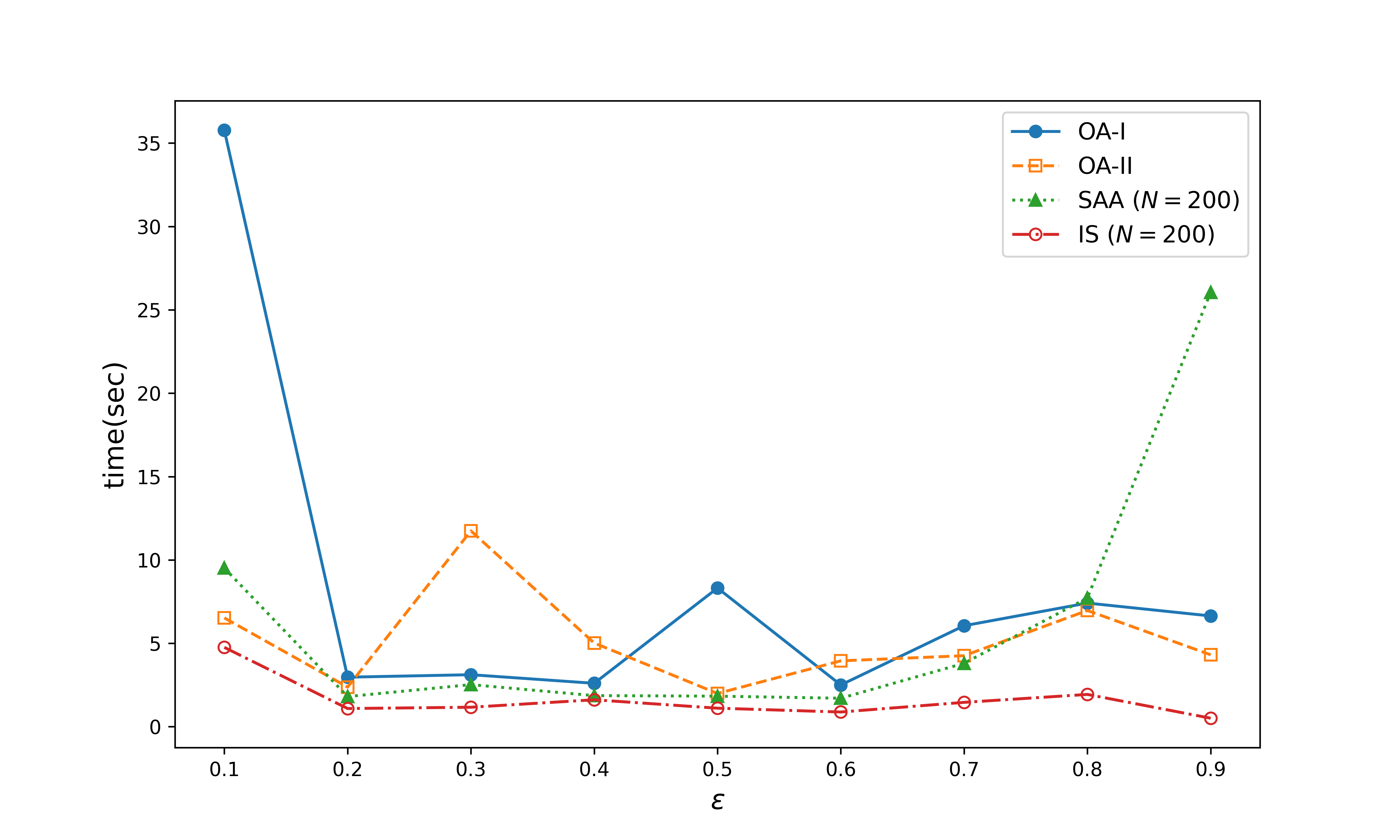} & \qquad \includegraphics[width=.4\textwidth]{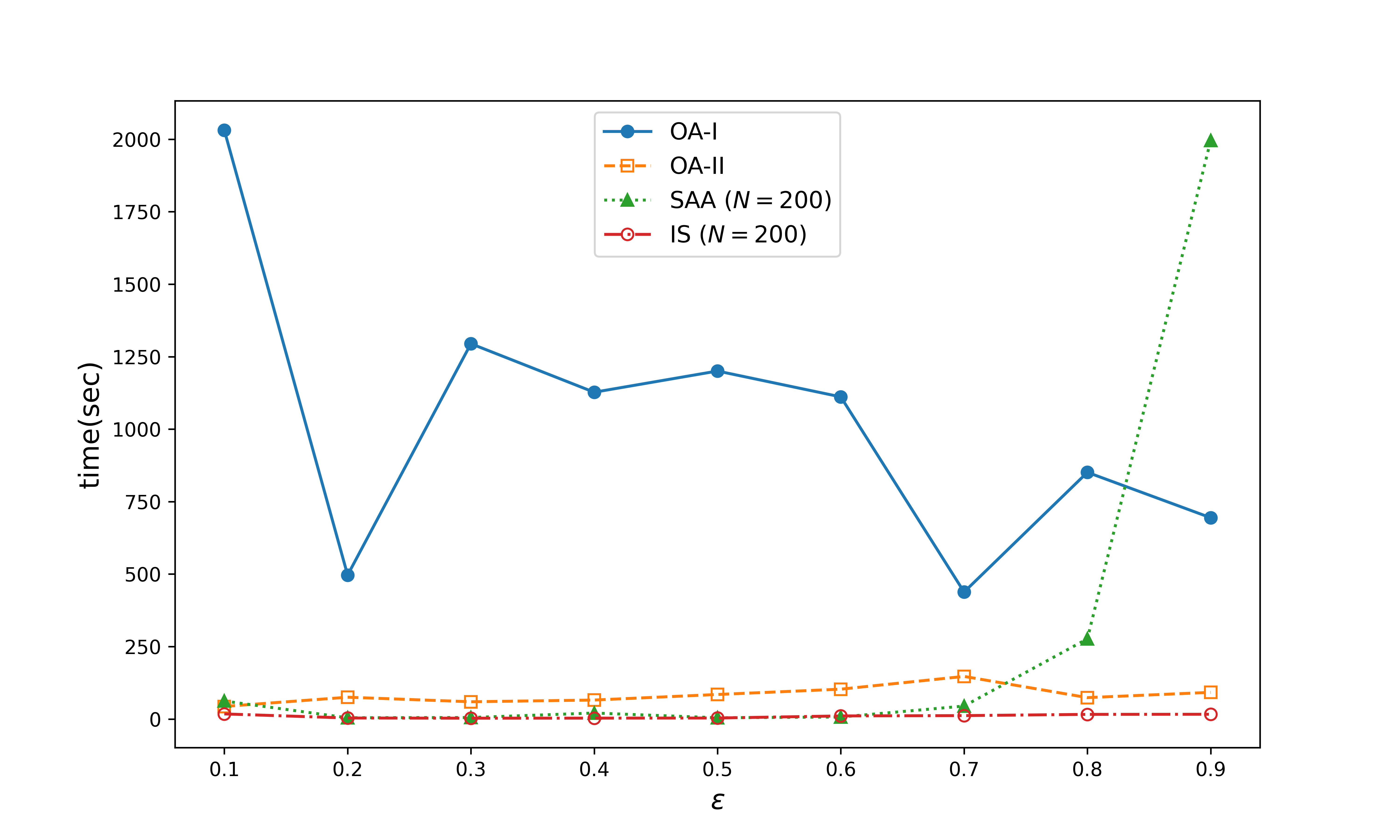} \vspace*{-0.5em} \\
{\small (a) $n=50$ and $m = 30$} & \qquad {\small (b) $n=100$ and $m = 50$}  \\
\end{tabular}
\caption{The average solution time of four approaches as a function of $\epsilon$ on two instances} 
\label{fig: epsilon}
\end{figure}

The results from Table \ref{tab: comparison} indicate that OA approaches can solve the instances to optimality when $\epsilon_i := \epsilon$ is small, which is commonly the case. However, it remains an intriguing question to investigate how well these approaches work for larger values of $\epsilon$. In this section, we calculate the average solution time for four approaches as a function of $\epsilon$, considering cases where $n = 50, m = 30$ and $n = 100, m = 50$. For each approach, we conduct five replications at each $\epsilon$ to compute the average time. All other settings remain consistent with those in Section \ref{sec: overall_performance}. Figure \ref{fig: epsilon} illustrates the average solution time of four different approaches as a function of $\epsilon$ on two different instances.

From Figure \ref{fig: epsilon}, we observe that OA-I exhibits high variability, with significant spikes around $\epsilon=0.1$ and $\epsilon=0.5$. SAA displays a sudden large spike at $\epsilon=0.9$, but maintains a low solution time in other cases. This could be attributed to the fact that the LP relaxation of SAA reformulation \eqref{eq: SAA-compact} is too weak to provide a strong lower bound at $\epsilon=0.9$. In contrast, the OA-II and IS approaches are more stable, generally maintaining low and consistent solution times across different values of $\epsilon$. OA-II improves upon OA-I, possibly by incorporating stronger reformulations that reduce computational time, resulting in more stable performance across various $\epsilon$ values.

\end{document}